\documentclass{gtart_h}

\def\ifplaintex{\expandafter\ifx\csname documentclass\endcsname\relax}

\ifplaintex 
\else
\fi

\def\gt{{\mathsurround=0pt\it $\cal G\mskip-2mu$eometry \&\ 
$\cal T\!\!$opology}}        %  journal title in recommended style

\def\gtm{{\mathsurround=0pt\it $\cal G\mskip-2mu$eometry \&\ 
$\cal T\!\!$opology $\cal M\mskip-1mu$onographs}}    %  for monographs

\def\gtp{{\mathsurround=0pt\it $\cal G\mskip-2mu$eometry \&\ 
$\cal T\!\!$opology $\cal P\!$ublications}}  % GT publications

\def\recd{{\small Received:\qua\receiveddate\ifx\reviseddate\relax
\else\qquad Revised:\qua\reviseddate\fi\par}} 

\def\volumenumber#1{\def\thevolumenumber{#1}}
\def\volumeyear#1{\def\thevolumeyear{#1}}
\def\volumename#1{\def\thevolumename{#1}}
\def\papernumber#1{\def\thepapernumber{#1}}
\def\pagenumbers#1#2{\def\startpage{#1}\def\finishpage{#2}}
\def\published#1{\def\publishdate{#1}}
\def\received#1{\def\receiveddate{#1}}
\def\revised#1{\def\reviseddate{#1}}
\def\accepted#1{\def\accepteddate{#1}}

\def\asciiemail#1{\def\theasciiemail{#1}}

\long\def\asciiabstract#1{\long\def\theasciiabstract{#1}}
\def\asciikeywords#1{\def\theasciikeywords{#1}}

\let\\\par
\let\thevolumenumber\relax\let\thepapernumber\relax
\let\thevolumeyear\relax\let\startpage\relax
\let\finishpage\relax\let\publishdate\relax\let\receiveddate\relax
\let\reviseddate\relax\let\accepteddate\relax\let\theasciititle\relax
\let\theasciiauthors\relax
\let\theasciiabstract\relax\let\theasciikeywords\relax

\let\theerratum\relax\let\theasciiemail\relax
\let\theshortauthors\relax\let\theshorttitle\relax

\def\startpage{1}\def\finishpage{15}\def\thepapernumber{77}

\volumenumber{2}
\volumename{Proceedings of the Kirbyfest}
\volumeyear{1999}

\long\def\maketitlep{   % start of definition of \maketitlep

\count0=\startpage

\gtm\nl        %   GT mongraphs (top left) 
{\small Volume \thevolumenumber: \thevolumename\nl 
\ifx\theerratum\relax\else Erratum \erratumnumber\nl\fi
Pages \startpage--\finishpage\nl}

\vglue 0.1truein   % top margin

{\parskip=0pt\leftskip 0pt plus 1fil\def\\{\par\smallskip}{\ifplaintex\large
\else\Large\fi\bf\thetitle}\par\medskip}   
\vglue 0.05truein 

{\parskip=0pt\leftskip 0pt plus 1fil\def\\{\par}{\sc\theauthors}
\par\medskip}%
 
\vglue 0.03truein 

{\small\leftskip 25pt\rightskip 25pt{\bf Abstract}\stdspace\theabstract

{\bf AMS Classification}\stdspace\theprimaryclass
\ifx\thesecondaryclass\relax\else; \thesecondaryclass\fi\par
{\bf Keywords}\stdspace \thekeywords\par}\vglue 7pt

}   % end of definition of \maketitlep

\font\phead=cmsl9 scaled 950
\font=cmsl9 scaled 1050
\font\pnum=cmbx10 scaled 913
\font\lnum=cmbx10 
\font\pfoot=cmsl9 scaled 950
\font=cmsl9 scaled 1050
\ifplaintex
\headline{\vbox to 0pt{\vskip -4.5mm\line{\small\phead\ifnum
\count0=\startpage ISSN 1464-8997 (on line)
1464-8989 (printed) \hfill {\pnum\folio}\else\ifodd\count0\def\\{ }% 
\ifx\theshorttitle\relax\thetitle\else\theshorttitle\fi\hfill{\pnum\folio}
\else\def\\{ and }{\pnum\folio}\hfill\ifx\theshortauthors\relax\theauthors
\else\theshortauthors\fi\fi\fi}\vss}}
\footline{\vbox to 0pt{\vglue 0mm\line{\small\pfoot\ifnum\count0=\startpage
Published \publishdate:\qua\copyright\ \gtp\hfill\else
\gtm, Volume \thevolumenumber\ (\thevolumeyear)\hfill\fi}\vss
}}
\else
\makeatletter
\def\@oddhead{{\small\lhead\ifnum\count0=\startpage ISSN 1464-8997 (on line)
1464-8989 (printed) \hfill {\lnum\number\count0}\else\ifodd\count0
\def\\{ }\ifx\theshorttitle\relax \thetitle \else\theshorttitle\fi\hfill
{\lnum\number\count0}\else\def\\{ and }{\lnum\number\count0}
\hfill\ifx\theshortauthors\relax 
\theauthors\else\theshortauthors\fi\fi\fi}}\def\@evenhead{@oddhead}
\def\@oddfoot{\small\lfoot\ifnum\count0=\startpage Published \publishdate:\qua\copyright\ \gtp\hfill\else
\gtm, Volume \thevolumenumber\ (\thevolumeyear)\hfill\fi}
\def\@evenfoot{@oddfoot}
\makeatother
\fi

\let\maketitlepage\maketitlep
\let\makeshorttitle\maketitlepage
\let\maketitle\maketitlepage

\newwrite\gtoutfile
\long\gdef\makeheadfile{  %%% start of definition of \makeheadfile
{\def\\{, }\def\s{ }
\immediate\openout\gtoutfile head.xxx
\immediate\write\gtoutfile{Proxy-for: \ifx\theasciiauthors\relax
\theauthors\else\theasciiauthors\fi\s<\ifx\theasciiemail\relax\theemail\else\theasciiemail\fi>}
\immediate\write\gtoutfile{\noexpand\\}
\immediate\write\gtoutfile{Authors: \ifx\theasciiauthors\relax
\theauthors\else\theasciiauthors\fi}
{\def\\{ }\immediate\write\gtoutfile{Title: \ifx\theasciititle\relax
\thetitle\else\theasciititle\fi}}
\immediate\write\gtoutfile{Subj-class: GT or SG, GR etc}
\immediate\write\gtoutfile{MSC-class: \theprimaryclass\ifx\thesecondaryclass\relax\else, \thesecondaryclass\fi}
\immediate\write\gtoutfile{Journal-ref: Geom. Topol. Monogr. \thevolumenumber\s
(\thevolumeyear) \startpage-\finishpage}
\immediate\write\gtoutfile{Comments: Published by Geometry and Topology Monographs at}
\immediate\write\gtoutfile{\s\s\s  http://www.maths.warwick.ac.uk/gt/GTMon\thevolumenumber/paper\thepapernumber.abs.html}
\immediate\write\gtoutfile{\noexpand\\}
\immediate\write\gtoutfile{}
\ifx\theasciiabstract\relax
\immediate\write\gtoutfile{\theabstract}\else
\immediate\write\gtoutfile{\theasciiabstract}\fi
\immediate\write\gtoutfile{}
\immediate\write\gtoutfile{\noexpand\\}
\immediate\write\gtoutfile{}
\immediate\closeout\gtoutfile}}  %%% end of definition of \makeheadfile

\def\maketitlepage{\maketitlep\makeheadfile}
\let\makeshorttitle\maketitlepage
\let\maketitle\maketitlepage

\volumenumber{7}
\volumename{Proceedings of the Casson Fest}
\volumeyear{2004}
\papernumber{14}
\pagenumbers{335}{430}
\received{9 April 2004}
\revised{8 December 2004}
\accepted{30 August 2004}
\published{11 December 2004}

\usepackage{amssymb,amsmath}

\def\S{Section }

\newtheoremstyle{Para}{14pt plus6pt minus6pt}%
{6pt plus3pt minus3pt}{\rm}{}{\bf}{}{1em}%
{\thmname{#1}\thmnumber{#2}\thmnote{\bf\stdspace[#3]}}

\theoremstyle{Para}
\newtheorem{para}{}[section]

\theoremstyle{remark}
\newtheorem{remark}[para]{Remark}
\newtheorem{remarks}[para]{Remarks}
\newtheorem{definition}[para]{Definition}
\newtheorem{definitions}[para]{Definitions}
\newtheorem{notation}[para]{Notation}

\theoremstyle{plain}
\newtheorem{theorem}[para]{Theorem}
\newtheorem{lemma}[para]{Lemma}
\newtheorem{proposition}[para]{Proposition}
\newtheorem{corollary}[para]{Corollary}

\renewcommand{\tilde}{\widetilde}

\numberwithin{equation}{para}

\newcommand\Number{\begin{para}}
\newcommand\EndNumber{\end{para}}
\newcommand\Definition{\begin{definition}}
\newcommand\EndDefinition{\end{definition}}
\newcommand\Definitions{\begin{definitions}}
\newcommand\EndDefinitions{\end{definitions}}
\newcommand\Theorem{\begin{theorem}}
\newcommand\EndTheorem{\end{theorem}}
\newcommand\Remark{\begin{remark}}
\newcommand\EndRemark{\end{remark}}
\newcommand\Remarks{\begin{remarks}}
\newcommand\EndRemarks{\end{remarks}}
\newcommand\Notation{\begin{notation}}
\newcommand\EndNotation{\end{notation}}
\newcommand\Lemma{\begin{lemma}}
\newcommand\EndLemma{\end{lemma}}
\newcommand\Proposition{\begin{proposition}}
\newcommand\EndProposition{\end{proposition}}
\newcommand\Corollary{\begin{corollary}}
\newcommand\EndCorollary{\end{corollary}}
\newcommand\Proof{\begin{proof}}
\newcommand\EndProof{\end{proof}}
\newcommand\Equation{\begin{equation}}
\newcommand\EndEquation{\end{equation}}

\newdimen\partindent
\global\partindent=28.5pt
\newcommand\Parts{\begingroup\parskip 3pt plus 2pt}
\newcommand\Part[1]{\par\noindent\hangindent\partindent
     \hbox to \partindent{\hskip .5\partindent minus .5\partindent
     {\rm#1}\enspace\hfill}\ignorespaces}
\newcommand\Subpart[1]{\par\noindent\hangindent1.8\partindent
     \hbox to 1.8\partindent{\hskip \partindent minus \partindent
     {\rm#1}\enspace\hfill}\ignorespaces}
\newcommand\EndParts{\par\endgroup}

\newcommand\chiminus{\chi_{\_{}}}
\newcommand\emm{{\mathfrak m}}
\newcommand\II{{\mathcal I}}
\newcommand\HH{{\mathcal H}}
\newcommand\sigmasurface{$T$--surface}
\newcommand\sigmasurfaces{$T$--surfaces}
\newcommand\sigmaxisurface{$T_{x_i}$--surface}

\newcommand\id{{\rm id}}
\newcommand\four{\empty}
\newcommand\inter{\mathop{\rm int}}
\newcommand\frontier{\mathop{\rm frontier}}

\newcommand\height{\mathop{\rm height}}
\newcommand\minheight{\mathop{\rm minheight}}
\newcommand\trace{\mathop{\rm trace}}
\newcommand\image{\mathop{\rm im}}

\newcommand\sltwo{{\rm SL}_2}
\newcommand\gltwo{{\rm GL}_2}
\newcommand\bigirth{\mathop{\rm bigirth}}
\newcommand\length{\mathop{\rm length}}
\newcommand\howmany{seven }
\newcommand\init{\mathop{\rm init}}
\newcommand\term{\mathop{\rm term}}

\begin{document}

\title[Two-surface knots]{Knots with only two strict essential surfaces}

\authors{Marc Culler\\Peter B Shalen}
\address{Department of Mathematics, Statistics, 
and Computer Science (M/C 249)\\University of Illinois at Chicago\\
851 S Morgan St, Chicago, IL 60607-7045, USA}
\gtemail{\mailto{culler@math.uic.edu}, \mailto{shalen@math.uic.edu}}
\asciiemail{culler@math.uic.edu, shalen@math.uic.edu}

\begin{abstract}
We consider irreducible $3$--manifolds $M$ that arise as knot
complements in closed $3$--manifolds and that contain at most two
connected strict essential surfaces.  The results in the paper relate
the boundary slopes of the two surfaces to their genera and numbers of
boundary components.  Explicit quantitative relationships, with
interesting asymptotic properties, are obtained in the case that $M$
is a knot complement in a closed manifold with cyclic fundamental
group.
\end{abstract}
\asciiabstract{%
We consider irreducible 3-manifolds M that arise as knot
complements in closed 3-manifolds and that contain at most two
connected strict essential surfaces.  The results in the paper relate
the boundary slopes of the two surfaces to their genera and numbers of
boundary components.  Explicit quantitative relationships, with
interesting asymptotic properties, are obtained in the case that M
is a knot complement in a closed manifold with cyclic fundamental
group.}

\asciikeywords{Knot complement, hyperbolic 3-manifold, boundary slope,
strict essential surface, essential homotopy, cyclic fundamental
group, character variety}

\maketitle
\cl{\small\it We dedicate this paper to Andrew Casson, in honor of
his 60th birthday.}

\section*{Introduction}
\addcontentsline{toc}{section}{Introduction}

It is well known that if $K$ is a torus knot of type $(p,q)$ in $S^3$,
then the exterior of $M(K)$ contains exactly two connected essential
surfaces up to isotopy.  (See Section 1 for precise definitions of
``essential surface'' and other terms used in this introduction.)  One
of the connected essential surfaces in $M(K)$ is an essential annulus
whose numerical boundary slope with respect to the standard framing is
$pq$.  The other connected essential surface in $M(K)$ has boundary
slope $0$ (ie, is a spanning surface) and has genus
$(p-1)(q-1)/2$. Note that the genus of the surface with boundary slope
$0$ is less than the boundary slope of the other essential surface.

This relationship between two quantities that are computed in entirely
different ways may appear coincidental. However, one of the results of
this paper, Corollary \ref{easycorollary}, asserts that a similar,
although weaker, inequality holds for any knot $K$ in a homotopy
$3$--sphere $\Sigma$ such that $M(K)$ is irreducible and has only two
essential surfaces up to isotopy.  In this case one of the essential
surfaces has boundary slope $0$ with respect to the standard framing
and the other has boundary slope $r\not=0$.  If the surface with
boundary slope $0$ has genus $g\ge 2$ and if $r \not=\infty$ 
then
$$ \frac{g - 1}{4\log_2(2g-2)} \le r^2 .$$ Thus $g$ is bounded above by
a function of $r$ which grows only slightly faster than a quadratic
function.

This result illustrates the general theme of this paper. We consider
knots in a closed, orientable $3$--manifold $\Sigma$, whose complements
are irreducible and contain at most two connected essential surfaces
that are {\it strict} in the sense of \ref{strict}.  For such knots
our results relate the boundary slopes of the connected strict
essential surfaces in the knot exteriors to their intrinsic
topological invariants (genera and numbers of boundary curves). Our
deepest results concern the case in which $\pi_1(\Sigma)$ is cyclic.
(This includes the case $\Sigma=S^3$.)  However, we also obtain some
non-trivial results of this type for an arbitrary $\Sigma$; studying a
general knot exterior in an arbitrary $\Sigma$ is equivalent to
studying an arbitrary irreducible {\it knot manifold}\/, ie, a
compact, irreducible, orientable $3$--manifold whose boundary is a
torus.

In Section \ref{fewsurfaces} we give a general qualitative
description of irreducible knot manifolds that contain at most two
connected strict essential surfaces up to isotopy.  Theorem
\ref{AtMostTwoSES} includes a complete classification of the
irreducible knot manifolds that contain at most one isotopy class of
connected strict essential surfaces; in particular they are all
Seifert fibered spaces. This completes a partial result proved in
\cite{boundsep}. Theorem \ref{AtMostTwoSES}, with Proposition
\ref{exceptionalprop}, also provides a dichotomy among the
irreducible knot manifolds that contain exactly two connected strict
essential surfaces up to isotopy. One subclass of such knot manifolds,
called exceptional graph manifolds, are defined by an explicit
classification (\ref{exceptional}). The complementary subclass,
called non-exceptional two-surface knot manifolds, are homeomorphic to
compact cores of one-cusped finite-volume hyperbolic $3$--manifolds.
Furthermore, the two connected strict essential surfaces in a
non-exceptional two-surface knot manifold are bounded, have distinct
boundary slopes, and have strictly negative Euler characteristic.

The subsequent sections are devoted to studying the relationships
among the boundary slopes and intrinsic topological invariants of
connected strict essential surfaces in non-exceptional two-surface
knot manifolds. Theorem \ref{disksandannuli} applies to an arbitrary
non-exceptional two-surface knot manifold $M$ and asserts, roughly
speaking, that when the two connected strict essential surfaces $F_1$
and $F_2$ in $M$ are isotoped into standard position with respect to
each other, they cut each other up into disks and annuli. This easily
implies Corollary \ref{chibound}, which asserts that for $i=1,2$ we
have
$$|\chi(F_i)| \le \frac{m_1m_2\Delta}{2}, $$ where $m_i$ denotes
the number of boundary components of $F_i$, and $\Delta$ denotes
geometric intersection number of the boundary slopes of
$F_1$ and $F_2$. (Thus if $M$ is a knot exterior in a closed manifold,
and if with respect to some framing of the knot the numerical boundary
slope of $F_i$ is $s_i=p_i/q_i$, then $\Delta=|p_1q_2-p_2q_1|$.)

For the case where the non-exceptional two-surface knot manifold $M$
arises as the complement of a knot $K$ in a closed $3$--manifold
$\Sigma$ with cyclic fundamental group, our main results are Theorem
\ref{easyboundconsequence} and Theorem \ref{hardboundconsequence}.
If $F_i$, $m_i$ and $\Delta$ are defined as above, if $g_i$
denotes the genus of $F_i$ and $q_i$ the denominator of its numerical
boundary slope with respect to any framing of $K$, and if $g_2\ge 2$,
Theorem \ref{easyboundconsequence} asserts that
$$\left(\frac{q_1}{\Delta}\right)^2\le\frac{4 m_2^2\log_2(2g_2-2)}{g_2-1}.$$
For the case of a knot in a homotopy $3$--sphere whose exterior is
irreducible and has only two essential surfaces up to isotopy, Theorem
\ref{easyboundconsequence} specializes to Corollary
\ref{easycorollary}. However, Theorem \ref{easyboundconsequence}
applies more generally to knots in $S^3$ and other lens spaces whose
exteriors contain three non-isotopic essential surfaces of which only
two are strict.  The figure eight knot and its sister are well known
examples of such knots.  Theorem \ref{hardboundconsequence} gives a
somewhat different inequality under the same hypotheses as Theorem
\ref{easyboundconsequence}. An examination of the inequality in the
statement of \ref{hardboundconsequence} will reveal that it may be
written in the form
$$\frac{q_1^2}{\Delta}\le\frac{m_2|\chi_1|}{m_1|\chi_2|}f(|\chi_2|),$$ where
$\chi_i=\chi(F_i)=2-2g_i-m_i$ for $i=1,2$, and $f(x)$ is a function of
a positive variable $x$ which grows more slowly than any
positive power of $x$. Theorem \ref{hardboundconsequence} is more difficult to
prove than Theorem \ref{easyboundconsequence}, and in a sense that we
shall explain in \ref{qualitativediscussion} it is qualitatively
stronger than Theorem \ref{easyboundconsequence}, although it does
not imply the latter theorem.

This paper, like our earlier papers \cite{boundsep} and
\cite{slopediff}, is based on the idea of using character variety
techniques to study the essential surfaces in a knot exterior.
Sections \ref{essentialsurfacesection}--\ref{actionsfromidealpoints}
are foundational in nature. Much of the work in these sections
consists of refining and systematizing material that has its origins
in such papers as \cite{varieties}, \cite{cgls}, and
\cite{ccgls}, concerning character varieties, actions on trees,
essential surfaces, and the norm on the homology of the boundary of a
knot manifold that was first used in the proof of the Cyclic Surgery
Theorem.

The material in Section \ref{generalsection} is crucial to the proofs
of all the main results of the paper. This section centers around the
study of the norm on the plane $H_1(M;\mathbb{R})$ in the case where $M$
is a non-exceptional two-surface knot manifold $M$. In this case the
ball of any radius with respect to the norm is a parallelogram for
which the slopes of the diagonals are the boundary slopes of the
strict essential surfaces $F_1$ and $F_2$ in $M$. The extra
information needed to determine the shape of such a parallelogram is
the ratio of the norms of the boundary slopes. Unlike the slopes
themselves, this ratio has no obvious topological
interpretation. Theorem \ref{generalinequality} asserts that this
ratio is bounded above by a topologically defined quantity associated to
the surfaces $F_1$ and $F_2$. This quantity, denoted
$\kappa(F_1,F_2)$, is defined in \ref{kappadef}.

A key ingredient in the proof of Theorem \ref{generalinequality} is
the study of degrees of trace functions of non-peripheral elements of
$\pi_1(M)$.  This appears to be the first application of information
of this type to the topology of $3$--manifolds.

Theorem \ref{disksandannuli}, the result mentioned above which
asserts that when $F_1$ and $F_2$ are isotoped into standard position
they cut each other up into disks and annuli, is relatively easy to
derive from Theorem \ref{generalinequality}. The results about knot
exteriors in a manifold with cyclic fundamental group depend on
combining Theorem \ref{generalinequality} with a fundamental result
from \cite{cgls} relating the norm on $H_1(M;\mathbb{R})$ to cyclic
Dehn fillings of $M$. Combining these directly gives a purely
topological result, Theorem \ref{generalknotinequality}, which asserts
that $$\frac{q_1^2}{\Delta}\le{\four}2\kappa(F_1,F_2).$$
The deepest theorems in the paper, Theorems
\ref{easyboundconsequence} and \ref{hardboundconsequence}, which
relate the boundary slopes and intrinsic topological invariants of
$F_1$ and $F_2$, are proved by combining Theorem
\ref{generalknotinequality} with combinatorial results, Proposition
\ref{easykappabound} and Proposition \ref{explicithardkappabound},
which relate
$\kappa(F_1,F_2)$ to more familiar topologically defined quantities.
These results depend on combining graph-theoretical arguments, given
in Sections \ref{firstbigirthsection} and \ref{secondbigirthsection}
respectively, with material in $2$-- and $3$--manifold topology that is
presented in Sections \ref{firstgenussection} and
\ref{secondgenussection} respectively.

The theorems proved in this paper remain true if the condition that
$M$ is a non-exceptional two-surface manifold is replaced by a
condition that is weaker, but more technical.  This condition is
described in \ref{weakerhyp}.  Computational evidence suggests that
there are many examples of knot exteriors in lens spaces that
satisfy this condition.

The work presented in this paper was partially supported by NSF
grant DMS 0204142.  

\section{Conventions}

\Number
In this paper the results about manifolds may be interpreted in the
smooth or PL category, or in the category in which objects are
topological manifolds, embeddings are locally flat, and polyhedra
contained in manifolds are tame. All of these categories are
equivalent in dimensions $\le3$, and we will often implicitly choose
one for the proof of a particular result.

We shall generally denote the unit interval $[0,1]\subset\mathbb{R}$ by
$I$. The Euler characteristic of a compact polyhedron $P$ will be
denoted by $\chi(P)$, and the cardinality of a finite set $X$ by $\#(X)$.
\EndNumber

\Number\label{basepoints}
Base points will often be suppressed when the choice of a base point
does not affect the truth value of a statement; for example, if $f$ is
a map between path-connected spaces $X$ and $Y$, to say that
$f_\sharp\co\pi_1(X)\to\pi_1(Y)$ is injective means that for some, and
hence for every, choice of base point $x\in X$, the homomorphism
$f_\sharp\co \pi_1(X,x)\to\pi_1(Y,f(x))$ is injective.
\EndNumber

\Number\label{standardaction}
Suppose that $M$ is a manifold or a polyhedron, that $x\in M$ is a
base point and that $(\tilde M,p)$ is a regular covering space of
$M$. Each choice of basepoint $\tilde x\in p^{-1}(x)\subset\tilde M$
determines an action of $\pi_1(M,x)$ on $\tilde M$.  An action of
$\pi_1(M,x)$ on $\tilde M$ will be termed {\it standard} if it arises
in this way from some choice of basepoint $\tilde x\in
p^{-1}(x)$. Note that if $x$ and $y$ are any base points in $X$ and if
$J\co \pi_1(M,x)\to\pi_1(M,y)$ is the isomorphism determined by some path
from $x$ to $y$, then pulling back a standard action of $\pi_1(M,y)$
via $J$ gives a standard action of $\pi_1(M,x)$.
\EndNumber

\Number
A path connected subspace $A$ of a path connected space $X$ will be
termed ``$\pi_1$--injective'' if the inclusion homomorphism from
$\pi_1(A)$ to $\pi_1(X)$ is injective. More generally, if $A$ is a
subspace of a space $X$, to say that $A$ is $\pi_1$--injective in $X$
will mean that each path component of $A$ is $\pi_1$--injective in the
path component of $X$ containing it.
\EndNumber

\Number\label{collarings}
Suppose that $F$ is a properly embedded codimension--$1$ submanifold of
a manifold $M$. By a {\it collaring} of $F$ in $M$ we mean an
embedding neighborhood of $F$, and $h\co F\times [-1,1]\to M$ such that
$h(x,0)=x$ for every $x\in F$ and $h(F\times[-1,1])\cap\partial
M=(\partial F)\times [-1,1]$. If $h$ is a collaring of $F$, we shall
set $V_h=h(F\times[-1,1])$, $V_h^{+1}=h(F\times[0,1])$ and
$V_h^{-1}=h(F\times[-1,0])$. We define a {\it collar neighborhood} of
$F$ to be a set that has the form $V_h$ for some collaring $h$ of $F$.

By definition, the submanifold $F$ is {\it two-sided} if it has a
collaring.  A {\it transverse orientation} of a two-sided submanifold
$F$ is an equivalence class of collarings, where two collarings $h$
and $h'$ are defined to be equivalent if $V_h^{+1}\cap V_{h'}^{-1}=F$.
\EndNumber

\Number
 As usual, we define a {\it homotopy} to be a continuous map
$H\co X\times I\to Y$ where $X$ and $Y$ are spaces, and for each $t\in I$
we denote the map $x\mapsto H(x,t)$ by $H_t$. The {\it inverse} $\bar
H$ of the homotopy $H\co X\times I\to Y$ is defined by $\bar
H(x,t)=H(x,1-t)$. By a {\it reparametrization} of $H$ we mean a map
$H'\co X\times[a,b]\to Y$, where $[a,b]\subset\mathbb{R}$ is a
non-degenerate interval, defined by $H'(x,t)=H(x,\alpha(t))$, where
$\alpha\co [a,b]\to I$ is a homeomorphism with $\alpha(a)=0$.

A homotopy $H\co X\times I\to Y$ is a {\it composition} of
homotopies $H^1,\ldots,H^k\co X\times I\to Y$ if there are real numbers
$t_0,\ldots,t_k$ with $0=t_0<t_1<\cdots<t_k=1$ such that
$H\vert_{[t_{i-1},t_i]}$ is a reparametrization of $H^i$ for $i=1,\ldots,k$.

A {\it path} in a space $Y$, ie, a map from $I$ to $Y$, may be
regarded as a homotopy $\{\star\}\times I\to Y$, where $\{\star\}$ is
a one-point space. By specializing the definitions given above we
obtain definitions of the inverse of a path, a reparametrization of a
path and a composition of paths.
\EndNumber

\Number\label{reducedhomotopydef}
Suppose that $M$ is a compact manifold and that $F\subset M$ is a
properly embedded submanifold of codimension $1$.  By a {\it homotopy
in $(M, F)$} we shall mean a homotopy $H\co K\times I\to M$, where $K$ is
some polyhedron, such that $H(K\times \partial I) \subset F$; we may
regard $H$ as a map of pairs $H\co (K\times I,K\times\partial
I)\to(M,F)$.
 
Now suppose that we are given a transverse orientation of $F\subset
M$, and an element $\omega$ of $\{-1,+1\}$. A homotopy $H$ in $(M,F)$
will be said to {\it start on the $\omega$ side} (or, respectively, to
{\it end on the $\omega$ side} if for some $\delta>0$ we have
$H(K\times[0,\delta])\subset V_h^{\omega}$ (or, respectively,
$H(K\times[1-\delta,1])\subset V_h^{\omega}$), where $h$ is a
collaring of $F$ realizing its transverse orientation; the condition
is independent of the choice of a collaring realizing the given
transverse orientation.
  
A homotopy $H$ in $(M,F)$ is a {\it basic homotopy} if $H^{-1}(F) = K
\times \partial I$. Note that every basic homotopy starts on the
$\omega$ side and ends on the $\omega'$ side for some
$\omega,\omega'\in\{-1,+1\}$.
  
Specializing these definitions to the case in which $K$ is a point, we
obtain the definitions of a path in $(M,F)$, of a basic path in
$(M,F)$, and of a path in $(M,F)$ which starts or ends on the $\omega$
side.

A basic path $\alpha$ in $(M,F)$ will be termed {\it essential} if it
is not fixed-endpoint homotopic to a path in $F$.  A basic homotopy
$H\co (K\times I,K\times\partial I)\to(M,F)$ will be termed {\it
essential} if for every $x\in K$ the basic path $\alpha_x\co t\mapsto
H(x,t)$ in $(M,F)$ is essential. Note that the condition that
$\alpha_x$ be essential depends only on the connected component of $x$
in $K$.
  
Now suppose that $F$ is a properly embedded, codimension--$1$
submanifold $F$ of a compact manifold $M$, and that $k$ is a positive
integer.  A homotopy $H\co  (K \times I, K \times \partial I)\to (M, F)$
will be called a {\it reduced homotopy of length $k$ in $(M,F)$} if we
may write $H$ as a composition of $k$ essential basic homotopies $H^1,
\ldots, H^k$ in such a way that, given a transverse orientation of
$F$, for each $i\in\{1,\ldots,k-1\}$ there is an element $\omega_i$ of
$\{-1,+1\}$ such that $H^i$ ends on the $\omega_i$ side and $H^{i+1}$
starts on the $-\omega_i$ side. Note that this condition is independent
of the choice of transverse orientation. Note also that, for any
choice of transverse orientation, $H$ starts on the same side as $H_1$
and ends on the same side as $H_k$.

We define a {\it reduced homotopy of length $0$} in $(M,F)$ to be a
map $H$ from $K$ to $F$. In this case we set $H_0=H_1=H$. If $H$ is a
reduced homotopy of length $0$ and $H'$ is a reduced homotopy of
length $\ge0$ for which $H'_1$ (or $H'_0$) is equal to $H$, we define
the {\it composition} of $H$ with $H'$ (or of $H'$ with $H$) to be
$H'$.
\EndNumber

\Number\label{minimalintersection} 
By a {\it closed curve} in a topological space $X$ we mean a map
$c\co S^1\to X$. If $c$ is a closed curve in a manifold $M$, and $F$ is a
properly embedded submanifold of codimension $1$ in $M$, we define the
{\it geometric intersection number of $c$ with $F$}\/, denoted
$\Delta(c,F)$ (or $\Delta_M(c,F)$ when we need to be more explicit),
to be the minimum cardinality of $g^{-1}(F)$, where $g$ ranges over
all closed curves homotopic to $c$.
\EndNumber

\Number
A {\it simple closed curve} in a manifold $M$ is a connected closed
$1$--manifold $C\subset M$. With a simple closed curve $C$ we can
associate a closed curve $c$ in $M$, well-defined modulo composition
with self-homeomorphisms of $S^1$, such that $c(S^1)=C$. If $F\subset
M$ is a properly embedded submanifold of codimension $1$, the {\it
geometric intersection number} $\Delta(C,F)=\Delta(c,F)$ is
well-defined, since composing $c$ with a self-homeomorphism of $S^1$
clearly does not change its geometric intersection number with $F$.
  
In particular, for any two simple closed curves $C$ and $C'$ in a
closed 2--manifold, $\Delta(C,C')$ is the geometric intersection number
of $c$ and $C'$ in the familiar sense.
\EndNumber

\Number\label{unislope}
If $T$ is a $2$--dimensional torus, we define a {\it slope} on $T$ to
be an isotopy class of homotopically non-trivial simple closed curves
in $T$. If $s_1$ and $s_2$ are slopes, we shall write
$\Delta(s_1,s_2)=\Delta(C_1,C_2)$ for any simple closed curves $C_i$
realizing the slopes $s_i$.
  
The isotopy classes of homotopically non-trivial {\it oriented} simple
closed curves in $T$ are in natural bijective correspondence with
elements of $H_1(T;\mathbb{Z})$ which are {\it primitive} in the sense of
not being divisible by any integer $>1$. Thus there is a natural
two-to-one map from the set of primitive elements of $H_1(T;\mathbb{Z})$
onto the set of slopes on $T$. We shall denote this map by
$\alpha\mapsto\langle\alpha\rangle$. We have
$\langle\alpha\rangle=\langle\alpha'\rangle$ if and only if
$\alpha'=\pm\alpha$.

If $T$ is a $2$--torus and $\alpha$ and $\beta$ are primitive elements
of $H_1(T;\mathbb{Z})$, then
$\Delta(\langle\alpha\rangle,\langle\beta\rangle)$ is the absolute
value of the homological intersection number of $\alpha$ and $\beta$.
\EndNumber
 
\Number\label{multislope}
If $C$ is a non-empty closed $1$--manifold in a $2$--torus $T$, and $C$
has no homotopically trivial components, then all components of $C$
have the same slope $s$. We call $s$ the {\it slope} of $C$.

Let $C_1,C_2$ be closed $1$--manifolds, with no homotopically trivial
components, in a torus $T$. Let $s_i$ and $m_i$ denote respectively
the slope and the number of components of $C_i$. Then $C_1$ and $C_2$
are isotopic to $1$--manifolds $C_1^0$ and $C_2^0$ such that
$\#(C_1^0\cap C_2^0)=m_1m_2\Delta(s_1,s_2)$. If $\#(C_1\cap
C_2)=m_1m_2\Delta(s_1,s_2)$ we shall say that $C_1$ and $C_2$ {\it
intersect minimally}\/. This implies that no arc in $C_1$ is
fixed-endpoint homotopic to any arc in $C_2$.
\EndNumber

\Number
An {\it essential surface} in an irreducible, orientable $3$--manifold
$M$ is a two-sided properly embedded surface in $M$ which is non-empty
and $\pi_1$--injective, and has no $2$--sphere components and no
boundary-parallel components.
\EndNumber

\Number\label{denominator}
We define a {\it knot manifold} to be a connected, compact, orientable
$3$--manifold $M$ such that $\partial M$ is a torus.

We will say that a knot manifold is {\it hyperbolic}
if it is homeomorphic to the compact core of a complete
hyperbolic manifold with finite volume.  

If $M$ is a knot manifold we will say that an element $\gamma\in
\pi_1(M)$ is {\it peripheral} if it is conjugate to an element of
the subgroup $\image(\pi_1(\partial M)\to\pi_1(M))$.

If $K$ is a (tame) knot in a closed, orientable $3$--manifold $\Sigma$,
the {\it exterior} of $K$, defined to be the complement of an open
tubular neighborhood of $K$, will be denoted by $M(K)$. Note that
$M(K)$ is well-defined up to ambient isotopy in $\Sigma$, and in
particular up to homeomorphism, and that it is a knot manifold. A {\it
  meridian} of $K$ is a non-trivial simple closed curve in the torus
$\partial M(K)$ which bounds a disk in the tubular neighborhood
$\Sigma-\inter M(K)$. Such a curve exists and is unique up to isotopy.
Thus there is a well-defined {\it meridian slope} in the torus
$\partial M(K)$. A {\it meridian class} for $K$ is a primitive element
$\mu$ of $H_1(\partial M(K);\mathbb{Z})$ such that $\langle\mu\rangle$ is
the meridian slope.  According to \ref{unislope}, $K$ has exactly
two meridian classes, and they differ by a sign.
  
We define a {\it framing} for $K$ to be an ordered basis
$(\mu,\lambda)$ for $H_1(\partial M(K);\mathbb{Z})$ such that $\mu$ is a
meridian class.  In the special case where $\Sigma$ is an integer
homology $3$--sphere we define a framing $(\mu,\lambda)$ to be 
{\it standard} if $\lambda$ generates the kernel of the inclusion
homomorphism $H_1(\partial M(K); \mathbb{Z})\to H_1(M(K); \mathbb{Z})$.

If $(\mu,\lambda)$ is an arbitrary framing for $K$, there is a bijective
correspondence between slopes in $\partial M(K)$ and elements of
$\mathbb{Q}\cup\{\infty\}$ defined by
$$\langle\alpha\rangle\mapsto\omega(\alpha,\lambda)/\omega(\alpha,\mu),$$
where $\omega$ denotes homological intersection number. If $C$ is a
non-empty closed $1$--manifold in $\partial M(K)$ whose components are
homotopically non-trivial, we define the {\it numerical slope} of $C$,
with respect to a given framing, to be the element of $\mathbb{Q}\cup\{\infty\}$ corresponding to the slope of $C$ (in the sense of
\ref{unislope}).

Note that if $C$ is a non-empty closed $1$--manifold in $\partial
M(K)$, and if the numerical slope of $C$ in a given framing is written
in the form $p/q$, where $p$ and $q$ are relatively prime integers and
$q>0$, then $q=\Delta(s,\emm)$, where $s$ denotes the slope of $C$ and
$\emm$ denotes the meridian slope of $K$ in particular $q$ is
independent of the choice of framing. For this reason, if $s$ is a
slope on $\partial M(K)$, it is natural to refer to $\Delta(s,\emm)$
as the {\it denominator} of $s$.
\EndNumber

\Number \label{boundaryslope}
Let $M$ be a compact orientable $3$--manifold such that every component
of $\partial M$ is a torus.  Let $T$ be a component of $\partial M$.
If $F$ is an essential surface in $M$ that meets $T$ then $\partial F
\cap T$ is a $1$--manifold in $T$ having no homotopically trivial
components. Thus by \ref{multislope}, $\partial F \cap T$ has a
well-defined slope $s$, which we call the {\it boundary slope of $F$
on $T$}.

If $F$ is a bounded essential surface in a knot manifold $M$ then we
will refer to the boundary slope of $F$ on $\partial M$ simply as the
{\it boundary slope of $F$}.

If $K$ is a knot in a closed orientable $3$--manifold $\Sigma$ and if $F$
is a bounded essential surface in $M(K)$ then we define the {\it
numerical boundary slope} of $F$ with respect to any given framing
$(\mu,\lambda)$ to be the numerical slope of $\partial F$ with respect
to $(\mu,\lambda)$, in the sense of \ref{denominator}.

We define a {\it boundary class} of a bounded essential surface $F$ in
a knot manifold $M$ to be a primitive element $\alpha$ of
$H_1(\partial M;\mathbb{Z})$ such that $\langle\alpha\rangle$ is the
boundary slope of $F$. According to \ref{unislope}, a bounded
essential surface in a knot manifold has exactly two boundary classes,
and they differ by a sign.
\EndNumber

\Number\label{strict}
Suppose that $M$ is a compact, orientable irreducible $3$--manifold
whose boundary components are tori.  A connected essential surface in
$M$ is called a {\it semi-fiber} if either $F$ is a fiber in a
fibration of $M$ over $S^1$, or $F$ is the common frontier of two
$3$--dimensional submanifolds of $M$, each of which is a twisted
$I$--bundle with associated $\partial I$--bundle $F$. An essential
surface $F\subset M$ is termed {\it strict} if no component of $F$ is
a semi-fiber.  A strict essential surface has no disk components,
since an irreducible knot manifold which has an essential disk must be
a solid torus, and the essential disk in a solid torus is a fiber.
\EndNumber

\Number\label{meetsall}
Since a semi-fiber in a bounded 3--manifold $M$ must meet every
component of $\partial M$, any essential surface that is disjoint from
at least one component of $\partial M$ must be strict.
\EndNumber

\Number\label{strictlyvertical} 
Suppose that the orientable $3$--manifold $M$ is either a Seifert fibered
$3$--manifold or an $I$--bundle over a surface.  We will define a
surface in $M$ to be {\it vertical} if it is a union of fibers, and to
be {\it horizontal} if it is everywhere transverse to the fibers.  If
$M$ is an $I$--bundle over a surface $B$, then the {\it
vertical boundary} of $M$ is the inverse image of $\partial B$ under
the projection map.

It is known that if $M$ is a Seifert fibered manifold then an
essential surface in $M$ is either isotopic to a vertical surface or
to a horizontal surface. (A stronger version of this statement, for
essential laminations, is proved in \cite{brittenham}.  See also
\cite[\S II.7]{charsub} and \cite{eineklasse}.)  It is clear that
the manifold obtained by splitting a Seifert fibered manifold along a
horizontal surface has a natural $I$--bundle structure, and hence that
horizontal essential surfaces in Seifert fibered manifolds are never
strict.

An essential vertical annulus in a Seifert fibered manifold is the
inverse image under the Seifert fibration map of an essential properly
embedded arc in the base surface.  An essential vertical annulus in an
$I$--bundle is the inverse image under the fibration map of an
essential simple closed curve in the base.

Suppose that $M$ is a trivial $I$--bundle and that $F$ is a properly
embedded $\pi_1$--injective surface in $M$ such that all components of
$\partial F$ are contained in the same component $B$ of the $\partial
I$--bundle associated to $M$.  It follows from
\cite[Proposition 3.1]{waldhausen} that $F$ is isotopic to a
horizontal surface by an ambient isotopy that preserves the vertical
boundary of $M$, and that each component of $F$ is parallel to a
subsurface of $B$.

As a consequence of this fact we observe that if $M$ is a trivial
$I$--bundle, and $F$ is a properly embedded $\pi_1$--injective
surface in $M$ such that $\partial F$ is contained in the vertical
boundary of $M$, then $F$ is isotopic to a horizontal surface by an
ambient isotopy that preserves the vertical boundary of $M$.

Suppose that $M$ is an $I$--bundle and that $A$ is a disjoint union of
properly embedded annuli in $M$ none of which is parallel to an
annulus contained in the $\partial I$--bundle associated to $M$.  It
follows from \cite[Lemma 3.4]{waldhausen} in the case that $M$ is a
trivial $I$--bundle, and from \cite[Lemma 2]{brittenham} in the
twisted case that $A$ is isotopic to a vertical surface.

Suppose that $F$ is a properly embedded $\pi_1$--injective surface in
an $I$--bundle $M$ such that $\partial F$ is contained in the vertical
boundary of $M$.  Then $F$ is isotopic to a horizontal surface.  This
follows from \cite[Proposition 3.1 and Proposition 4.1]{waldhausen}.
\EndNumber

\Number\label{curveclass}
A closed curve $c$ in a path-connected space $X$ defines a conjugacy
class in $\pi_1(X)$, which we shall denote by $[c]$.
\EndNumber

\Number\label{graphdef}
By a {\it graph} we mean a CW--complex of dimension $\le 1$. Thus a
graph $\Gamma$ has an underlying space, which we shall denote by
$|\Gamma|$, and which need not be connected; $|\Gamma|$ is a disjoint
union of $0$--cells, called {\it vertices}, and open $1$--cells, called
{\it edges}. Each edge has the structure of an affine interval.

The vertices in the frontier of an edge will be called its {\it
endpoints}\/; each edge has either one or two endpoints.  We sometimes
will need to consider {\it oriented edges} in a graph. For every edge
$e$ there are two oriented edges whose underlying edge is $e$; these
will be called {\it orientations} of $e$.  If $\omega$ is an oriented
edge we shall denote by $|\omega|$ the underlying edge of $\omega$ and
by $-\omega$ the opposite orientation to $\omega$.  An oriented edge
$\omega$ has an {\it initial vertex} denoted $\init(\omega)$ and a
{\it terminal vertex} denoted $\term(\omega)$.  The {\it valence} of a
vertex $v$ is the number of oriented edges whose initial vertex is
$v$.

In the last four sections of this paper, the underlying space
$|\Gamma|$ of a graph $\Gamma$ will often arise as a subpolyhedron of a
PL manifold.

A {\it subgraph} of a graph $\Gamma$ is a graph $\Gamma'$ such that
$|\Gamma'|\subset|\Gamma|$ and every vertex or edge of $\Gamma'$ is a
vertex or edge of $\Gamma$. A graph $\Gamma$ is said to be {\it
connected} if $|\Gamma|$ is connected, and a {\it component} of a
graph $\Gamma$ is a subgraph $C$ such that $|C|$ is a connected
component of $|\Gamma|$.

An {\it edge path} of length $n>0$ in a graph $\Gamma$ is a sequence
$(\omega_1, \ldots, \omega_n)$ of oriented edges of $\Gamma$ such that
$\term(\omega_i) = \init(\omega_{i+1})$ for $i = 1, \ldots, n-1$.  If
$\omega_i \not= -\omega_{i+1}$ for $i = 1, \ldots, n-1$ then we will
say that the edge path is {\it reduced}.  The {\it track} of an edge
path $(\omega_1, \ldots, \omega_n)$ is the subgraph of $\Gamma$ whose
edges are $|\omega_1|, \ldots ,|\omega_n|$ and whose vertices are their
endpoints.

An {\it arc} in a graph $\Gamma$ is a subgraph $A$ of $\Gamma$ such
that $|A|$ is homeomorphic to a (possibly degenerate) closed interval
in $\mathbb{R}$. A {\it circuit} in $\Gamma$ is a subgraph $C$ of $\Gamma$ such
that $|C|$ is homeomorphic to $S^1$. 

The {\it length} of a finite graph $\Gamma$ is the number of edges of
$\Gamma$.
\EndNumber

\section{Essential surfaces}
\label{essentialsurfacesection}

In this section we collect several general results about essential
surfaces which will be used in Sections
\ref{actionsfromidealpoints}--\ref{generalsection}.

The next two results are proved in \cite{klaffshalen}.  We restate
them here for completeness.

\Proposition[Proposition 1.1 of \cite{klaffshalen}]
\label{boundaryincompressible}
Suppose that $F$ is a bounded essential surface in an irreducible knot
manifold $M$, and suppose that $\alpha$ is a path in $F$ which has its
endpoints in $\partial F$ and is fixed-endpoint homotopic in $M$ to a
path in $\partial M$. Then $\alpha$ is fixed-endpoint homotopic in $F$
to a path in $\partial F$.
\EndProposition

\Proposition[Proposition 1.3 of \cite{klaffshalen}]
\label{mfldwithtorus}
Let $M$ be a compact orientable irreducible $3$--manifold containing an
essential torus $T$ and let $M'$ be the manifold obtained by splitting
$M$ along $T$ and let $q\co M'\to M$ denote the quotient map.  Let $F$ be
a connected properly embedded surface in $M$ which is not isotopic to
$T$.  Then $F$ is a strict essential surface if and only if it is
isotopic to a surface $S$ transverse to $T$ such that
\Parts
  \Part{(1)} each component of $q^{-1}(S)$ is essential in the component
  of $M'$ containing it; and
  \Part{(2)} some component of $q^{-1}(S)$ is a strict essential surface
  in the component of $M'$ containing it.
\EndParts
\EndProposition

\Proposition\label{whychiisnegative}
For any essential surface $F$ in a hyperbolic knot manifold $M$ we
have $\chi(F)<0$.
\EndProposition

\Proof
Since $M$ is orientable and $F$ is two-sided, $F$ must be orientable.
Since $M$ is irreducible and is not a solid torus, $F$ cannot be a
$2$--sphere or a disk, and the hyperbolicity of $M$ implies that $F$
cannot be a torus.  Now suppose that $F$ is an annulus.  Consider the
submanifold $N$ of $M$ which is the closure of a regular neighborhood
of $F \cup \partial M$.  The frontier of $N$ has either one or two
components, each of which is a torus containing a simple closed curve
that is isotopic to a core curve of $F$.  In particular, no frontier
component of $N$ can be contained in a ball.  Furthermore since $F$ is
essential, no frontier component of $N$ can be boundary-parallel.
Since $M$ has no essential tori, it follows that each frontier
component of $N$ bounds a solid torus, which must be a component of
$\overline{M-N}$.  Since $F$ is a $\pi_1$--injective annulus and each
component of $\overline{M-N}$ is a solid torus it follows that $M$ is
Seifert-fibered.  This contradicts the fact that $M$ is hyperbolic.
\EndProof

\Proposition\label{linemup}
Let $S$ and $F$ be essential surfaces in a compact, irreducible,
orientable $3$--manifold $M$.  Suppose that $F$ is connected, and that
every component of $S$ is isotopic to $F$. Then there exist a
collaring $h$ of $F$ in $M$ and a finite set $Y\subset [-1,1]$ such
that $S$ is isotopic to $h(F\times Y)$.
\EndProposition

\Proof
Let $\mathcal N$ denote the set of all connected submanifolds $N$ of $M$
such that either (a) $N$ is a component of $S$, or (b) $N$ is
$3$--dimensional, each component of $\frontier_M N$ is a component of
$S$, and $N$ can be given the structure of a trivial $I$--bundle in such a way
that $N\cap\partial M$ is the vertical boundary of $N$. Then $\mathcal N$
is finite, and is non-empty since $S\ne\emptyset$. Hence we may choose
$N_0\in{\mathcal N}$ which is maximal with respect to inclusion.

We claim that $S\subset N_0$. Suppose this is false. Then some
component $S_0$ of $S$ is disjoint from $N_0$. Let us choose a
component $S_1$ of $\frontier_M N_0$. Since $S_0$ and $S_1$ are isotopic
to $F$ and hence to each other, it follows from \cite[Lemma 5.3]{waldhausen}
that $S_0\cup S_1$ is the associated $I$--bundle of a
trivial $I$--bundle $H\subset M$. In particular $H\in{\mathcal N}$. Hence
if $H\supset N_0$, we have a contradiction to the maximality of
$N_0$. The other possibility is that $H\cap N_0=S_1$. However, in this
case, since $H,N_0\in{\mathcal N}$, it is clear that $H\cup N_0\in{\mathcal
N}$, and we again have a contradiction to maximality. This proves our
claim.

Now let $M_0$ denote a regular neighborhood of $N_0$ in $M$. Then
$M_0$ may be given the structure of a trivial $I$--bundle in such a way that
$M_0\cap\partial M$ is the vertical boundary of $M_0$. Since
$S\subset N_0$, we may regard $S$ as a $\pi_1$--injective surface in
$M_0$ whose boundary is contained in the vertical boundary of
$M_0$. It now follows from \ref{strictlyvertical} that $S$ is
isotopic to a horizontal surface in $M_0$. This implies the conclusion
of the proposition.
\EndProof

\Proposition\label{thicknessmultiplies}
Let $F$ be a connected essential surface in a compact, irreducible,
orientable $3$--manifold $M$. Let $h$ denote a collaring of $F$ in
$M$. Suppose that $Y$ is a finite subset of $[-1,1]$, and let $S$
denote the essential surface $h(F\times Y)\subset M$. Let $y_0\in Y$
be given, let $K\subset F$ be a compact polyhedron and set $K_0 =
h(K\times {y_0})\subset S$.  Then we have $t_S(K_0)\ge\#(Y)\cdot
t_F(K)$.
\EndProposition

\Proof
Set $\nu=\#(Y)$.  If $\nu=1$ the assertion is trivial.  If
$\nu>1$ we may assume without loss of generality that $\{-1,+1\}\subset
Y$.  According to Definition \ref{thicknessdef}, what we need to
prove is that if for a given positive integer $\theta$ there is a
reduced homotopy $H\co (K\times I,K\times\partial I)\to(M,F)$ of length
$\theta-1$ such that, for some $c\in I$, the map $H_c$ is the
inclusion map $K\hookrightarrow F\subset M$, then there is a reduced
homotopy $H'\co (K\times I,K\times\partial I)\to(M,S)$ of length
$\nu\theta-1$ such that, for some $c'\in I$, the map $H'_{c}$ is the
inclusion map $K_0\hookrightarrow S\subset M$.

Let us give $F$ the transverse orientation determined by the collaring
$h$. We may write $H$ as a composition of $\theta-1$ essential basic
homotopies $H^1, \ldots, H^{\theta-1}$ in such a way that for each
$i\in\{1,\ldots,\theta-2\}$ there is an element $\omega_i$ of $\{-1,+1\}$
such that $H^i$ ends on the $\omega_i$ side and $H^{i+1}$ starts on
the $-\omega_i$ side. Let $\omega_0$ and $\omega_{\theta-1}$ denote
the elements of $\{-1,+1\}$ such that $H^1$ starts on the $-\omega_0$
side and $H^{\theta-1}$ ends on the $\omega_{\theta-1}$ side.

We set $M'=M-h(F\times(-1,1))=\overline{M-V_h}$, and we fix a map
$q\co M'\to M$ such that $q(h(x,j))=x$ for every $(x,j)\in
F\times\{-1,1\}$, and such that $q$ maps $M-V_h$ homeomorphically onto
$M-F$. For $i=1,\ldots\theta-1$, since $H^i$ is a basic homotopy,
there is a homotopy $\tilde H^i\co K\times I\to M'\subset M$ such that
$H^i=q\circ\tilde H^i$. Now for $i=0,\ldots,\theta-1$, define a
homotopy $J^i\co K\times I\to M$ by $J^i(x,t)=h(x,-\omega_i(2t-1))$.
Since $\{-1,+1\}\subset Y$, the $\tilde H^i$ are essential basic
homotopies and the $J^i$ are reduced homotopies of length $\nu-1$. We
may define the required homotopy $H'$ to be a composition of
$J^0,\tilde H^1,J^1,\tilde H^2,\ldots,\tilde
H^{\theta-1},J^{\theta-1}$.  (In particular for $\theta = 1$ we have
$H' = J^0$.)
\EndProof

\Remark
The inequality in Proposition \ref{thicknessmultiplies} can
presumably be shown to be an equality, but we will not need this.
\EndRemark

\section{Dual surfaces}
\label{dualsurfacesection}

The material in this section overlaps with material that has been
presented in \cite{cgls} and \cite{ccgls}, but we have found it convenient to
provide a self-contained account of it.

\Number\label{treedef}
By a {\it tree} we mean a graph $T$ such that $|T|$ is
$1$--connected. Since an edge in a tree is determined by its endpoints
and its endpoints are always distinct, a tree has the structure of a
geometric simplicial complex arising from the affine structure on the
edges.  If $T$ is a tree, $E_T\subset|T|$ will denote the set of all
midpoints of edges of $T$. For any two vertices $s,s'$ of $T$ there is
a unique arc having $s$ and $s'$ as endpoints. The length of this arc
will be denoted by $d_T(s,s')$, or simply by $d(s,s')$ when it is
clear which tree is involved. If we regard $d_T$ as a distance
function, the set of vertices of $T$ becomes an integer metric space.
\EndNumber

\Number\label{gammatree}
Suppose that $\Gamma$ is a group.  By a {\it $\Gamma$--tree} $T$ we
will mean a tree $T$ equipped with a simplicial action of $\Gamma$.
More explicitly, this means an action on the underlying space $|T|$
under which vertices are always carried to vertices, and edges are
carried to edges via affine homeomorphisms.  In general we will leave
the action itself unnamed and implicit in the notation for a
$\Gamma$--tree; the effect of an element $\gamma\in\Gamma$ on a point
$x$ of $|T|$ will ordinarily be denoted $\gamma\cdot x$. We will say
that a $\Gamma$--tree $T$ is {\it trivial} if for some vertex $s$ of
$T$ we have $\Gamma\cdot s=s$.

If $\rho\co  \Gamma \to G$ is a homomorphism of groups and if $T$ is a
$G$--tree then we may define a simplicial action of $\Gamma$ on the
tree $T$ by $\gamma\cdot x = \rho(\gamma)\cdot x$ for any point $x$ in
$|T|$.  The resulting $\Gamma$--tree will be called the {\ it pull-back}
of the $G$--tree $T$ via $\rho$.
\EndNumber

\Definition\label{bipartitedef}
Let $\Gamma$ be a group. A $\Gamma$--tree $T$ will be termed {\it
bipartite} if for every vertex $s$ of $T$ and every $\gamma\in\Gamma$,
the integer $d_T(s,\gamma \cdot s)$ is even.

A $\Gamma$--tree is said to be {\it without inversions} if for every
$\gamma\in\Gamma$ and every edge $e$ of $T$ such that $\gamma\cdot
e=e$, the element $e$ fixes both endpoints of $e$ (and hence fixes $e$
pointwise).  Note that a bipartite $\Gamma$--tree is in particular a
$\Gamma$--tree without inversions.
\EndDefinition

The term ``bipartite'' is motivated by the following result.

\Proposition\label{whybipartite}
Suppose that $\Gamma$ is a
group and that $T$ is a bipartite $\Gamma$--tree. Then the set of
vertices of $T$ is a disjoint union of two $\Gamma$--invariant subsets
$X_0$ and $X_1$ such that each edge of $T$ has one endpoint in $X_0$
and one endpoint in $X_1$.  In particular, the quotient graph
$T/\Gamma$ is bipartite.
\EndProposition
  
\Proof
Fix a vertex $s_0\in T$. For $i=0,1$, define $X_i$ to be
the set of all vertices $s$ of $T$ such that $d(s_0,s)\equiv
i\pmod2$. If $s$ and $s'$ are any two vertices of $T$, we have
$d(s,s')=d(s_0,s)+d(s_0,s')-2l$, where $l$ is the length of the
intersection of the  arcs joining $s_0$ to $s$ and to $s'$. In
particular, $d(s,s')\equiv d(s_0,s)+d(s_0,s')\pmod2$. It follows that
the distance between two vertices of $X_0$ or between two vertices of
$X_1$ is even, while the distance between a vertex of $X_0$ and a
vertex of $X_1$ is odd. The definition of a bipartite $\Gamma$--tree
therefore implies that $X_0$ and $X_1$ are $\Gamma$--invariant.
Furthermore, if two vertices $s$ and $s'$ are joined by an edge of $T$
then $d(s,s')$ is the odd number $1$, and hence one of the vertices
$s,s'$ must be in $X_0$ and the other in $X_1$.
\EndProof

\Definition
Let $\Gamma$ be a group and let $T$ be a
$\Gamma$--tree. We shall define the {\it length} of an element
$\gamma\in\Gamma$ relative to the $\Gamma$--tree $T$, denoted
$\lambda_T(\gamma)$, by
$$\lambda_T(\gamma)=\min_sd_T(s,\gamma\cdot s),$$
where $s$ ranges over the vertices of $T$.  It is clear that conjugate
elements of $\Gamma$ have the same length.  If $\Gamma = \pi_1(X)$ for
some path-connected space $X$, and if $c$ is a closed curve in $X$,
then we will set
$$\lambda_T(c) = \lambda_T(\gamma) $$
where $\gamma$ is an arbitrary
element of the conjugacy class $[c]$ (see \ref{curveclass}).
\EndDefinition
  
\Definition\label{sigmamapdef} 
Let $M$ be a $3$--manifold, let $(\tilde M,p)$ denote its universal
covering space, and let $T$ be a $\pi_1(M)$--tree without
inversions. We shall say that a map $f\co \tilde M\to |T|$ is {\it
equivariant} if it is $\pi_1(M)$--equivariant with respect to the
action of {$\pi_1(M)$} on $T$ and some standard action (see
\ref{standardaction}) of $\pi_1(M)$ on $\tilde M$.  We shall say that
$f$ is {\it transverse} if it is transverse to $E_T$.  If $f\co \tilde
M\to {|T|}$ is a transverse equivariant map, we have
$f^{-1}(E_T)=p^{-1}(S)$ for a unique properly embedded surface
$S\subset M$. The surface $S$ will be denoted by $S_f$.
  
When we are given a $3$--manifold $M$ and a $\pi_1(M)$--tree $T$ without
inversions, we define a {\it \sigmasurface} in $M$ to be a surface
that has the form $S_f$ for some transverse equivariant map $f\co \tilde
M\to |T|$, where $\tilde M$ denotes the universal covering space of
$M$.
\EndDefinition

\Remark
In Definition \ref{sigmamapdef}, $\pi_1(M)$ is understood to be
defined in terms of an unspecified base point. It follows from the
remark on change of base point in \ref{standardaction} that if $x$
and $y$ are points of $M$, if $T$ is a $\pi_1(M,x)$--tree, and if we
give $T$ the structure of a $\pi_1(M,y)$--tree by pulling back the
action of $\pi_1(M,x)$ via the isomorphism
$J\co \pi_1(M,y)\to\pi_1(M,x)$ determined by some path from $y$ to $x$,
then a map $f\co \tilde M\to {|T|}$ is equivariant (in the sense of
\ref{sigmamapdef}) when we regard $T$ as a $\pi_1(M,x)$--tree if and
only if it is equivariant when we regard $T$ as a
$\pi_1(M,y)$--tree. From this it follows that the statements made in
this section are independent of the choice of a base point, and in
accordance with the convention described in \ref{basepoints}, base
points will be suppressed.
\EndRemark
  
\Proposition\label{surfacenonempty}
If $M$ is an orientable $3$--manifold and $T$ is a $\pi_1(M)$--tree
without inversions which is non-trivial (see \ref{gammatree}), then
any \sigmasurface\ in $M$ is non-empty.
\EndProposition

\Proof
Let $\tilde M$ denote the universal cover of $M$. Suppose
that $f\co \tilde M\to |T|$ is a transverse equivariant map such that
$S_f=\emptyset$. Then $f$ maps $\tilde M$ into a component $C$ of
$|T|-E_T$. Such a component contains only one vertex, say $s$. Since
$f$ is equivariant, $\Gamma$ must leave $C$ invariant, and since
$\Gamma$ acts simplicially on $T$ it must fix $s$ and therefore be a
trivial action.
\EndProof

\Remark\label{coho}
A non-trivial homomorphism $f\co \pi_1(M)\to\mathbb{Z}$
determines a non-trivial $\pi_1(M)$--tree $T$, where $|T|$ is the real
line, and the vertices of $T$ are the integers.  If $\phi \in
H^1(M;\mathbb{Z})$ corresponds to $f$ under the natural isomorphism
between $\mathop{\rm Hom}(\pi_1(M), \mathbb{Z})$ and $H^1(M;\mathbb{Z})$,
and if $S$ is a $T$--surface, then $S$ is an essential surface which
represents the class that is the Poincar\'e--Lefschetz dual of $\phi$ in
$H_2(M, \partial M; \mathbb{Z})$.
\EndRemark

\Proposition\label{bipartitemanifold}
If $M$ is an orientable $3$--manifold and $T$ is a bipartite
$\pi_1(M)$--tree, then for any \sigmasurface\ $S\subset M$ there are
closed subsets $A_0$ and $A_1$ of $M$ which are $3$--dimensional
submanifolds, such that $A_0\cap A_1=\frontier A_0=\frontier A_1=S$.
\EndProposition

\Proof
According to Proposition \ref{whybipartite}, the set of
vertices of $T$ is a disjoint union of two $\Gamma$--invariant subsets
$X_0$ and $X_1$, such that each edge of $T$ has one endpoint in $X_0$
and one endpoint in $X_1$. For $i=0,1$, let $Y_i$ denote the union of
the closures of all components of $|T|-E_T$ which contain vertices in
$X_i$. Then the $Y_i$ are $\Gamma$--invariant, and $Y_0\cap
Y_1=\frontier Y_0=\frontier Y_1 =E_T$. Now suppose that $S\subset M$
is a \sigmasurface, so that $S=S_f$ for some transverse equivariant
map $f\co \tilde M\to {|T|}$, where $(\tilde M,p)$ denotes the universal covering
space of $M$. Since $f$ is $\pi_1(M)$--equivariant and transverse to
$E_T$, the closed set $\tilde A_i=f^{-1}(Y_i)\subset\tilde M$ is a
$\pi_1(M)$--invariant $3$--dimensional submanifold, and $\tilde A_0\cap
\tilde A_1=\frontier \tilde A_0=\frontier \tilde
A_1=f^{-1}(E_T)=p^{-1}(S)$. Hence $\tilde A_i=p^{-1}(A_i)$ for some
closed set $A_i\subset M$ which is a $3$--dimensional submanifold, and
$A_0\cap A_1=\frontier A_0=\frontier A_1=S$.
\EndProof

\Proposition\label{Deltaandtranslation}
Suppose that $M$ is a compact, orientable $3$--manifold, that $T$ is a
$\pi_1(M)$--tree without inversions, and that $S\subset M$ is a
\sigmasurface. Then for any closed curve $c$ in $M$ we have
$\lambda_T(c)\le\Delta_M(c,S)$.
\EndProposition

\Proof
Set $\Delta=\Delta_M(c,S)$ and $E=E_T$. We may assume $c$
to be chosen within its homotopy class so that $\#(c^{-1}(S))=\Delta$.
Let $(\tilde M,p)$ denote the universal covering of $M$.  According to
the definition of a \sigmasurface, we have $S=S_f$ for some transverse
equivariant map $f\co \tilde M\to |T|$; in particular,
$f^{-1}(E)=p^{-1}(S)$. We first consider the degenerate case in which
$f(p^{-1}(c(S^1)))$ contains no vertex of $T$. In this case, $f$ maps
each component of $p^{-1}(c(S^1))$ into a single edge $e$ of $T$.  It
then follows from equivariance that some element of $[c]$ leaves $e$
invariant, and hence fixes the endpoints of $e$ since {$T$ is a
$\pi_1(M)$--tree without inversions}. Hence we have $\lambda_T(c)=0$ in
this case, and the conclusion follows.

We may therefore assume that $f$ maps some point $\tilde x\in
p^{-1}(c(S^1))$ to a vertex $s$ of $T$. After reparametrizing $c$ if
necessary we may assume that $p(\tilde x)=c(1)$, where $1$ is the
standard base point of $S^1$. Let $q\co I\to S^1$ be a path representing
a generator of $\pi_1(S^1,1)$, set $\alpha=c\circ q$, choose a lift
$\tilde\alpha\co I\to\tilde M$ of $\alpha$, and set
$\beta=f\circ\alpha\co I\to |T|$.  Then $\beta(0)=s$ and
$\beta(1)=\gamma\cdot s$ for some element $\gamma$ of $[c]$ in
$\pi_1(M)$.

Since $p^{-1}(S)=f^{-1}(E)$, we have
$$\#(\beta^{-1}(E))=\#(\alpha^{-1}(S))=\Delta.$$ 
Hence if $\mathcal E$ denotes the set of edges of $T$ whose midpoints lie
in $\beta(I)$, we have $\#({\mathcal E})\le\Delta$. If $X$ denotes the
subgraph of $T$ consisting of all vertices of $T$ and of those edges
that belong to $\mathcal E$, then $\beta$ can clearly be deformed to a
path in $X$, and hence to an arc in $X$. This arc has length at most
$\Delta$ since $X$ has at most $\Delta$ edges. Hence
$$d_T(s,\gamma\cdot s)=d_T(\beta(0),\beta(1))\le\Delta,$$
and by the definition of translation length we have
$\lambda_T(c)\le\Delta$.
\EndProof

\Definition\label{thicknessdef}
Let $S$ be an essential surface in a compact, orientable, irreducible
$3$--manifold $M$. Let $K\subset S$ be a compact polyhedron which is
$\pi_1$--injective in $S$. We define the {\it thickness} of $K$
(relative to $S$) to be the supremum of all integers $\theta>0$ for
which there is a reduced homotopy $H\co (K\times I,K\times\partial
I)\to(M,S)$ of length $\theta-1$ such that, for some $t\in I$, the map
$H_t$ is the inclusion map $K\hookrightarrow S\subset M$.
The thickness of $K$ will be denoted by $t_S(K)$, or by $t(K)$ when
there is no danger of confusion.  Note that $t_S(K)$ is either a
strictly positive integer or $+\infty$.  Moreover, if $S$ is a
semi-fiber then $t_S(S) = +\infty$, and hence $t_S(K) = +\infty$ for any
compact $\pi_1$--injective polyhedron $K\subset S$.
\EndDefinition

\Theorem\label{dualsurface}
Suppose that $M$ is an irreducible knot manifold and that $T$ is a
non-trivial (see \ref{gammatree}) bipartite $\pi_1(M)$--tree. Then
there is an essential \sigmasurface\ $S\subset M$ which has the
following properties.
\Parts
  \Part{(1)} For any closed curve $c$ in $\partial M$ we have
  $\lambda_T(c)=\Delta_{\partial M}(c,\partial S)$.
  \Part{(2)} If $K\subset S$ is any $\pi_1$--injective, connected,
  compact polyhedron such that $\chi(K)<0$, if $t\le t_S(K)$ is a
  positive integer, and if $\Theta\le\pi_1(M)$ is the subgroup defined
  up to conjugacy by $\Theta=\image(\pi_1(K)\to\pi_1(M))$, then
  $\Theta$ fixes an arc of length $t$ in $T$.
\EndParts
\EndTheorem

The next \howmany lemmas are needed for the proof of Theorem
\ref{dualsurface}.

\Lemma\label{mapextends} 
Suppose that $M$ is a compact orientable $3$--manifold and that $T$ is
a $\pi_1(M)$--tree without inversions. Let $(\tilde M,p)$ denote the
universal covering space of $M$ and fix a standard action of
$\pi_1(M)$ on $M$.  Suppose that $L\subset K\subset M$ are compact
polyhedra and that $\tilde K$ is a union of components of
$p^{-1}(K)\subset\tilde M$.  Set $\tilde L=\tilde K\cap p^{-1}(L)$.
Suppose that $H$ is a subgroup of $\pi_1(M)$ which stabilizes $\tilde
K$ and that $V$ is a connected $H$--invariant subset of $|T|$ such that
$\frontier_{|T|} V$ contains no vertices of $T$.  Suppose that
$g_L\co \tilde L\to \bar V$ is a PL map such that $g_L(h\cdot x) = h\cdot
g_L(x)$ for all $x \in \tilde L$ and $h\in H$.  Then $g_L$ may be
extended to a PL map $g_K\co \tilde K\to \bar V$ such that $g_K(\tilde
K-\tilde L)\subset V$, and $g_K(h\cdot x) = h\cdot g_K(x)$ for all $x
\in \tilde K$ and $h\in H$.
\EndLemma

\Proof
Fix a triangulation of $K$ in which $L$ is a subcomplex,
and give $\tilde K$ the triangulation inherited from that of $K$.  For
$i=-1,0,1,2,3$, let $\tilde K^{(i)}$ denote the $i$--skeleton of
$\tilde K$ (so that $\tilde K^{(-1)}=\emptyset$), and set $\tilde
L^{(i)}=\tilde K^{(i)}\cap\tilde L$. We shall recursively construct,
for $i=-1,0,1,2,3$, a piecewise-linear map $g^{(i)}\co \tilde K^{(i)}\cup
\tilde L\to \bar V$ which extends $g_L$, maps $\tilde K^{(i)}-\tilde
L^{(i)}$ into $V$, and is $H$--equivariant in the sense that
$g_i(h\cdot x) = h\cdot g_i(x)$ for all $x \in \tilde K^{(i)}\cup
\tilde L$ and $h\in H$. We take $g^{(-1)}=g_L$.  Suppose that
$g^{(i)}$ has been constructed for a given $i\le2$. Let $\mathcal D$ be a
complete set of orbit representatives for the action of $H$ on the set
of $(i+1)$--simplices of $\tilde K$ that are not contained in $\tilde
L$. For each $\delta\in{\mathcal D}$ we extend $g^{(i)}\vert_{\partial\delta}$ to
a PL map $h_\delta\co \bar\delta\to \bar V$; the extension exists because
$V$, being a connected subset of the underlying space of the tree $T$,
is contractible. Furthermore, since $\frontier_{|T|} V$ contains no
vertices of $T$, there is a neighborhood $N$ of $\frontier_{|T|} V$
relative to $\bar V$ such that $N$ is a $1$--manifold with boundary and
$\frontier_{|T|} V\subset \partial N$.  Hence by general position we
may choose the extension $h_\delta$ so that it maps the open simplex
$\delta$ into $V$.

For each point $x\in\tilde K^{(i+1)}-(\tilde K^{(i)}\cup
\tilde L^{(i+1)})$ there exist a unique $\gamma\in\pi_1(M)$ and a unique
$\delta\in{\mathcal D}$ such that $\gamma\cdot x\in\delta$. We set
$g^{(i+1)}(x)= \gamma^{-1}\cdot h_\delta(\gamma\cdot x)$. For
$x\in \tilde K^{(i)}\cup \tilde L$ we set $g^{(i+1)}(x)=g^{(i)}(x)$. The
extension $g^{(i+1)}$ of $g^{(i)}$ defined in this way is clearly
piecewise-linear and $H$--equivariant. Since $V$ is $H$--invariant and
since $g^{(i)}(\tilde K^{(i)}-\tilde L^{(i)})\subset V$, we have
$g^{(i+1)}(\tilde K^{(i+1)}-\tilde L^{(i+1)})\subset V$.
\EndProof

\Lemma\label{equivariantmapexists}
Suppose that $M$ is a compact orientable $3$--manifold and that $T$ is
a $\pi_1(M)$--tree without inversions.  Then there exists a transverse
equivariant map $f$ from the universal cover $\tilde M$ of $M$ to
$|T|$.
\EndLemma

\Proof
First we fix a standard action of $\pi_1(M)$ on $\tilde M$
and apply Lemma \ref{mapextends}, taking $K=M$,
$L=\emptyset$, $\tilde K=\tilde M$, $H=\pi_1(M)$, $V=|T|$, and
taking $g_L=g_\emptyset$ to be the empty map.  This gives
a $\pi_1(M)$--equivariant PL map $g=g_M\co \tilde M\to |T|$.

If we subdivide the triangulation of $\tilde M$, and subdivide the
simplicial complex $T$, so that $g$ is simplicial, then $g$ is
transverse to every non-vertex point in the subdivision of $T$. In
particular, every (open) edge $e$ of $T$ contains a point to which $g$
is transverse. Let $\mathcal E$ denote a complete set of orbit
representatives for the action of $\pi_1(M)$ on the set of edges of
$T$. For each $e\in {\mathcal E}$ choose a point $z_e\in e$ such that $g$
is transverse to $z_e$, and set $E_0=\{\gamma\cdot
z_e\co \gamma\in\pi_1(M),e\in{\mathcal E}\}$. Since $T$ is a $\pi_1(T)$--tree
without inversions, $E_0$ contains exactly one point in each edge of
$T$, and there is a $\pi_1(M)$--equivariant self-homeomorphism $\eta$
of $|T|$ such that $\eta(E_0)=E_T$.  Then $f=\eta\circ g$
is a transverse equivariant map.
\EndProof

\Lemma\label{fixupthesurface} 
Suppose that $M$ is a compact orientable $3$--manifold and that $T$ is
a $\pi_1(M)$--tree without inversions. Let $(\tilde M,p)$ denote the
universal cover of $M$, and suppose that $f\co \tilde M\to |T|$ is a
transverse equivariant map.  Suppose that $S_f$ is the frontier of a
compact $3$--dimensional submanifold $A$ of $M$. Suppose that $X\subset
A$ is a compact connected $3$--manifold with the following properties:
\Parts
  \Part{(i)} every component of $\frontier_AX$ is a properly embedded
  $2$--manifold $C\subset A$ with $\partial C\subset \inter S_f$; and
  \Part{(ii)} for some component $\tilde X$ of $p^{-1}(X)$, $f(\tilde
  X\cap p^{-1}( S_f))$ is a single point.
\EndParts
Then $\frontier_M\overline{X-A}$ is a \sigmasurface.
\EndLemma

\Remarks\label{degeneratecase}
(1)\qua Of course condition (i) in the hypothesis of Lemma
\ref{fixupthesurface} holds vacuously in the special case where $X$
is a component of $A$, since then $\frontier_AX=\emptyset$.

(2)\qua If $f$ is a transverse equivariant map then it follows from the
definition of $S_f$ that $f( p^{-1}( S_f))\subset E_T$. Thus condition
(ii) in the hypothesis of Lemma \ref{fixupthesurface} may be
paraphrased by saying that $f(\tilde X\cap p^{-1}( S_f))$ is a single
point of $E_T$.

(3)\qua Condition (ii) in the hypothesis of Lemma \ref{fixupthesurface},
together with the equivariance of $f$, implies that for {\it every}
component $\tilde X$ of $p^{-1}(X)$, $f(\tilde X\cap p^{-1}( S_f))$ is
a single point of $E_T$.
\EndRemarks

\Proof[Proof of Lemma \ref{fixupthesurface}]
We fix a standard action of $\pi_1(M)$ on $\tilde M$ that makes $f$ a
$\pi_1(M)$--equivariant map. According to the hypotheses, we may choose
a component $\tilde X_0$ of $p^{-1}(X)$ and a point $\mu\in E_T$
(cf Remark \ref{degeneratecase}(2)) such that $f(\tilde X_0\cap p^{-1}(
S_f))=\{\mu\}$. We denote by $e$ the edge of $T$ whose midpoint is
$\mu$, and by $H$ the stabilizer of $\tilde X_0$ in $\pi_1(M)$. Then
$H$ stabilizes $\tilde X_0\cap p^{-1}( S_f)$, and by the
$\pi_1(M)$--equivariance of $f$ it follows that $H$ fixes $\mu$.  Since
$T$ is a $\pi_1(M)$--tree without inversions it follows that $H$ fixes
$e$.
  
We denote by $B$ the closure of $M-A$ in $M$, and by $C$ the closure
of $A-X$ in $A$. We set $F=\frontier_AX\subset\partial C$ and $J=X\cap
S_f\subset\partial B$. We denote by $Z$ a submanifold of $C$ such that
$Z\cap\partial C=F$ and such that the pair $(Z,F)$ is homeomorphic to
$(F\times I,F\times\{0\})$; we set $F^\sharp=(\partial Z)-F$. (Note
that if $X$ is a component of $A$ then $F^\sharp=F=\emptyset$,
cf Remark \ref{degeneratecase}(1).)  Likewise, we denote by $Y$ a
submanifold of $B$ such that $Y\cap\partial B=J$ and such that the
pair $(Y,J)$ is homeomorphic to $(J\times I,J\times\{0\})$; we set
$J^\sharp=(\partial Y)-J$. Then the $3$--manifold $X^\sharp=X\cup Y\cup
Z\subset M$ deform-retracts to $X$. Hence the component $\tilde
X_0^\sharp$ of $p^{-1}(X^\sharp)$ containing $\tilde X_0$ is precisely
invariant under $H$, in the sense that $\gamma\cdot\tilde
X_0^\sharp=\tilde X_0^\sharp$ for any $\gamma\in H$, while
$(\gamma\cdot\tilde X_0^\sharp)\cap\tilde X_0^\sharp=\emptyset$ for
any $\gamma\in \pi_1(M)-H$. Note that $\frontier_M
X^\sharp=J^\sharp\cup F^\sharp$, and hence that $\frontier_{\tilde M}
\tilde X_0^\sharp=\tilde J_0^\sharp\cup \tilde F_0^\sharp$, where
$\tilde J_0^\sharp=p^{-1}(J^\sharp)\cap\tilde X_0^\sharp$ and $\tilde
F_0^\sharp=p^{-1}(F^\sharp)\cap\tilde X_0^\sharp$.

Let $\tilde A_0$ denote the component of $p^{-1}(A)$ containing
$\tilde X_0$, and let $V$ denote the component of $|T|-E_T$ containing
$f(\inter\tilde A_0)$. Then 
$$f(\tilde F_0^\sharp)\subset f(\tilde
  X_0^\sharp\cap p^{-1}(A))\subset f(\tilde A_0)\subset 
\bar V.$$
Note that $V$ is one of the two components of $|T|-E_T$ whose closures
contain $\mu$; we shall denote the other one by $W$. Since $H$ fixes
$e$, it leaves $V$ and $W$ invariant. Since $f$ is transverse to $E_T$
and maps $\tilde X_0\cap p^{-1}( S_f)$ to $\mu$, every component of
$p^{-1}(B)$ which meets $\tilde A_0$ must be mapped into $\bar W$ by
$f$. In particular we have
$$f(\tilde J_0^\sharp)\subset f(\tilde X_0^\sharp\cap
  p^{-1}(B))\subset\bar W.$$
We set $\tilde F_0= p^{-1}(F)\cap\tilde X_0$, and we define a map
$g_0\co \tilde J_0^\sharp\cup \tilde F_0^\sharp\cup\tilde F_0\to |T|$ to
agree with $f$ on $\tilde J_0^\sharp\cup \tilde F_0^\sharp$ and to map
$\tilde F_0$ to $\mu$. Then $g_0$ is well-defined since $(\tilde
J_0^\sharp\cup \tilde F_0^\sharp)\cap\tilde F_0\subset
p^{-1}(S_f)\cap\tilde X_0\subset f^{-1}(\{\mu\})$, and it is
$H$--equivariant because $H$ fixes $\mu$.  Now set $P=X\cup Y\subset
X^\sharp$, and note that $\tilde X_0^\sharp$ is the union of the two
$H$--invariant sets $\tilde Z_0=\tilde X_0^\sharp\cap p^{-1}(Z)$ and
$\tilde P_0=\tilde X_0^\sharp\cap p^{-1}(P)$, and that $\tilde
Z_0\cap\tilde P_0=\tilde F_0$. It follows from Lemma \ref{mapextends}
that $g_0\vert_{\tilde F_0^\sharp\cup\tilde F_0}$ may be extended to a PL
$H$--equivariant map $g_Z\co \tilde Z_0\to \bar V$ such that $g_Z(\tilde
Z_0-(\tilde F_0^\sharp\cup\tilde F_0))\subset V$, and that
$g_0\vert_{\tilde J_0^\sharp\cup\tilde F_0}$ may be extended to a PL $H$--equivariant map
$g_P\co \tilde P_0\to \bar V$ such that $g_P(\tilde P_0-(\tilde
J_0^\sharp\cup\tilde F_0))\subset W$. Now define a map
$g_{X^\sharp}\co \tilde X_0^\sharp\to \bar V$ to agree with $g_Z$ on
$\tilde Z_0$ and with $g_P$ on $\tilde P_0$. Since $g_{X^\sharp}$ is
$H$--equivariant and agrees with $f$ on $\frontier_{\tilde M}\tilde
X_0^\sharp$, and since $\tilde X_0^\sharp$ is precisely invariant
under $H$, there is a unique $\pi_1(M)$--equivariant map $f'\co \tilde
M\to |T|$ which agrees with $g_{X^\sharp}$ on $\tilde X^\sharp_0$ and
with $f$ on $\tilde M-\pi_1(M)\cdot\tilde X^\sharp_0$. (For any
$x\in\tilde X^\sharp_0$ and any $\gamma\in\pi_1(M)$ we set
$f'(\gamma\cdot x)=\gamma\cdot g_{X^\sharp}(x)$; the precise
invariance of $\tilde X_0^\sharp$ and the $H$--equivariance of
$g_{X^\sharp}$ guarantee that $f'$ is well-defined.)

If we set $S'=\frontier_M\overline{X-A}=(S_f-(S_f\cap A))\cup F$, it
follows from the construction of $f'$ that
$(f')^{-1}(E_T)=p^{-1}(S')$. The construction also shows that the
restriction of $f'$ to a small neighborhood of $\tilde X_0^\sharp$ is
transverse to $E_T$. Since $f'$ is $\pi_1(M)$--equivariant and agrees
with $f$ outside $\pi_1(M)\cdot\tilde X_0^\sharp$, it is everywhere
transverse to $E_T$. Hence $f'$ is a transverse equivariant map and $S_{f'}=S'$. In
particular, $S'$ is a $T$--surface.
\EndProof

The following slight variant of Lemma \ref{fixupthesurface} will also
be useful. The proof will show that it is essentially a special case
of \ref{fixupthesurface}.

\Lemma\label{fixupthesurfacetoo}
Suppose that $M$ is a compact orientable $3$--manifold and that $T$ is
a $\pi_1(M)$--tree without inversions. Let $(\tilde M,p)$ denote the
universal cover of $M$, and suppose that $f\co \tilde M\to |T|$ is a
transverse equivariant map.  Suppose that $S_f$ is the frontier of a
compact $3$--dimensional submanifold $A$ of $M$. Suppose that $X\subset
A$ is a compact connected $3$--manifold with the following properties:
\Parts
  \Part{(i)} every component of $\frontier_AX$ is a properly embedded
  $2$--manifold $C\subset A$ with $\partial C\subset \inter S_f$;
  \Part{(ii)} $X\cap S_f$ is connected;
  \Part{(iii)} $\pi_1(X\cap S_f)\to\pi_1(X)$ is surjective. 
\EndParts
Then $\frontier_M\overline{X-A}$ is a \sigmasurface.
\EndLemma

\Proof
We will prove this by showing that the hypotheses of Lemma
\ref{fixupthesurfacetoo} imply those of Lemma
\ref{fixupthesurface}. The only point to check is that condition 
(ii) of \ref{fixupthesurface} follows from the hypotheses of Lemma
\ref{fixupthesurfacetoo}. If $\tilde X$ is any component of the
covering space $p^{-1}(X)$ of $X$, the surjectivity of $\pi_1(X\cap
S_f)\to\pi_1(X)$ implies that the induced covering space $\tilde X\cap
p^{-1}( S_f)$ of $X\cap S_f$ is connected. Since $f$ maps $p^{-1}(
S_f)$ into the discrete set $E_T$, it must map the connected subset
$\tilde X\cap p^{-1}( S_f)$ to a single point.
\EndProof

\Definition
Let $M$ be an orientable $3$--manifold and let $(\tilde M,p)$ denote
its universal covering. Let $T$ be a bipartite $\pi_1(M)$--tree, and
let $f\co \tilde M\to |T|$ be a transverse equivariant map. We shall say
that $f$ has a {\it folded boundary-annulus} if there is an annulus
$R\subset\partial M$ such that $\inter R$ is a component of $\partial
M-\partial S_f$, and for some component $\tilde R$ of $p^{-1}(R)$, the
components of $\partial\tilde R=\tilde R\cap p^{-1}(E_T)$ are mapped
by $f$ to the same point of $E_T$. We shall say that $f$ has a {\it
big folded $I$--bundle} if there is a submanifold $X$ of $M$ which is
an $I$--bundle over a compact, connected surface of negative Euler
characteristic, such that (i) $Y=X\cap S_f$ is the associated
$\partial I$--bundle of $X$, (ii) $Y$ is $\pi_1$--injective in $S_f$,
and (iii) for some component $\tilde X$ of $p^{-1}(X)$, the set
$\tilde X\cap p^{-1}(Y)=\tilde X\cap f^{-1}(E_T)$ is mapped by $f$ to
a single point of $E_T$.
\EndDefinition

\Remark\label{whatnobigfoldedIbundlesmeans}
Suppose that $M$ is an orientable $3$--manifold, whose universal cover
we denote by $(\tilde M,p)$. Suppose that $T$ is a bipartite
$\pi_1(M)$--tree $T$, and that $f\co \tilde M\to |T|$ is a transverse
equivariant map. Suppose that a submanifold $X$ of $M$ is an
$I$--bundle over a compact, connected surface of negative Euler
characteristic, that $Y=X\cap S_f$ is the associated $\partial
I$--bundle of $X$, and that $Y$ is $\pi_1$--injective in $S_f$. If
$\tilde X$ is any component of $p^{-1}(X)$, then $\tilde X$ is a
covering space of $X$, and is therefore a connected $I$--bundle whose
associated $\partial I$--bundle is $\tilde Y=\tilde X\cap p^{-1}(Y)$.
Thus $\tilde Y$ has at most two components, and since $\tilde Y\subset
p^{-1}(S_f)$, each component of $\tilde Y$ must be mapped by $f$ to a
point of $E_T$. Furthermore, $\tilde Y$ has exactly two components if
and only if $\tilde X$ is a trivial $I$--bundle. If we assume that the
transverse equivariant map $f$ has no big folded $I$--bundles, then
$f(\tilde Y)$ cannot be a single point; hence in this case the
$I$--bundle $\tilde X$ must be trivial, and $f$ must map the two
components of $\tilde Y$ to distinct points of $E_T$.
\EndRemark

\Lemma\label{controlledfolding}
Suppose that $M$ is an irreducible knot manifold and that $T$ is a
non-trivial bipartite $\pi_1(M)$--tree. Then there is a transverse
equivariant map $f\co \tilde M\to |T|$ such that (i) $S_f$ is essential
and (ii) $f$ has no folded boundary-annuli or big folded $I$--bundles.
\EndLemma

\Proof
We denote by $(\tilde M,p)$ the universal cover of $M$ and fix a
standard action of $\pi_1(M)$ on $\tilde M$.  For any compact,
orientable surface $F$, we set
$$\chiminus(F)=\sum_C\max(0,-\chi(C)),$$ where $C$ ranges over the
components of $F$. For any \sigmasurface\ $F$, we define the {\it
complexity} $c(F)\in\mathbb{N}^4$ to be $(b(F),\chiminus(F),t(F),s(F))$,
where $b(F)$ is the number of components of $\partial F$, $t(F)$ is
the number of components of $F$ that are tori or annuli, and $s(F)$ is
the number of closed components of $F$. We endow the set $\mathbb{N}^4$
with the lexicographical order. It follows from Lemma
\ref{equivariantmapexists} that the set of all \sigmasurfaces\ in $M$
is non-empty. Hence there is a \sigmasurface\ $S$ which has minimal
complexity among all \sigmasurfaces. By the definition of a
\sigmasurface, we have $S=S_f$ for some transverse equivariant map $f\co \tilde M\to
|T|$.  We shall prove Lemma \ref{controlledfolding} by showing that
$S$ is essential and that $f$ has no folded boundary-annuli or big
folded $I$--bundles.
 
Since $T$ is a bipartite $\pi_1(M)$--tree, Proposition
\ref{bipartitemanifold} asserts that there are compact
$3$--dimensional submanifolds $A_0$ and $A_1$ of $M$ such that
$A_0\cap A_1=\frontier A_0=\frontier A_1=S$.

To show that $S$ is essential, we first observe that since the
$\pi_1(M)$--tree $T$ is by hypothesis non-trivial, we have
$S\ne\emptyset$ according to Proposition \ref{surfacenonempty}. Next,
we shall show that $S$ is $\pi_1$--injective. Assume it is not. Then by
a standard consequence of the loop theorem, there is a disk $D\subset
M$ such that $D\cap S=\partial D$, and $\partial D$ does not bound a
disk in $S$. Let $X\subset M$ be a ball such that $X\cap
S\subset\partial X$, and $R=X\cap S$ is a regular neighborhood of
$\partial D$ in $S$. Then $\partial X - \inter R$ is a disjoint union
of two disks $D_1$ and $D_2$. We must have $X\subset A_j$ for some
$j\in\{0,1\}$. The hypotheses of Lemma \ref{fixupthesurfacetoo}
clearly hold with this choice of $X$, and with $A=A_j$. Hence
\ref{fixupthesurfacetoo} implies that the surgered surface
$S'=(S-R)\cup D_1\cup D_2$ is a
\sigmasurface. We shall reach a contradiction by showing that
$c(S')<c(S)$. Note that $b(S')=b(S)$.

Let $S_0$ denote the component of $S$ containing $R$. Then
$S_0'=(S-R)\cup D_1\cup D_2$ has either one or two components, and
$\chi(S_0')=\chi(S_0)+2$. We first consider the case in which
$\chi(S_0)<0$. In this case, at most one component of $S_0'$ can be a
disk; and since the core curve $\partial D$ of $R$ does not bound a
disk in $S_0$, no component of $S_0'$ can be a sphere. It follows that
$$\chiminus(S_0')\le-\chi(S_0')+1=-\chi(S_0)-1=\chiminus(S_0)-1.$$
Since $\chiminus(S_0')<\chiminus(S_0)$, it is clear that
$\chiminus(S')<\chiminus(S)$. Hence $c(S)<c(S')$ in this case. There
remains the case in which $\chi(S_0)\ge0$. Since the core curve of $R$
is homotopically non-trivial in $S$, the only possibilities are that
$S_0$ is a torus and $S_0'$ is a sphere, or that $S_0$ is an annulus
and $S_0'$ consists of two disks. In both subcases we have
$\chiminus(S_0')=\chiminus(S_0)=0$, so that
$\chiminus(S')=\chiminus(S)$, whereas $t(S')<t(S)$. Hence
$c(S')<c(S)$, and the proof of $\pi_1$--injectivity is complete.

Next we show that no component of $S$ is a $2$--sphere.  If $S$ does
have a $2$--sphere component then by irreducibility, any $2$--sphere
component of $S$ must bound a ball $X$, and since we have shown that
$S$ is $\pi_1$--injective, any component of $S$ contained in $X$ must
itself be a sphere. Hence if we take $X$ to be minimal with respect to
inclusion among all balls in $M$ bounded by components of $X$, then
$X\cap S=\partial X$.  We must have $X\subset A_j$ for some
$j\in\{0,1\}$. The hypotheses of Lemma \ref{fixupthesurfacetoo}
clearly hold with this choice of $X$, and with $A=A_j$. (See Remark
\ref{degeneratecase}(1).)  Hence if $S_1$ denotes the sphere
$\partial X$, \ref{fixupthesurfacetoo} implies that $S'=S-S_1$ is a
\sigmasurface. But we obviously have $b(S')=b(S)$,
$\chiminus(S')=\chiminus(S)$, $t(S')=t(S)$ and $s(S')=s(S)-1$. Thus
$c(S')<c(S)$, a contradiction. This shows that no component of $S$ is
a $2$--sphere.

To prove that $S$ is essential, it remains to show that $S$ has no
boundary-parallel component. If $S$ does have a boundary-parallel
component $S_2$, then $S_2$ has a {\it region of
boundary-parallelism}, ie, a submanifold $X$ of $M$ such that
$\frontier X = S_2$ and $(X, S_2)$ is homeomorphic to $(S_2\times I,
S_2\times\{1\})$.  Since we have shown that $S$ is $\pi_1$--injective
and has no sphere components, it follows from \cite{waldhausen} that
any component of $S$ contained in $X$ must itself
have a region of boundary-parallelism which is contained in $X$. Hence
if we choose $S_2$ so that $X$ to be minimal with respect to inclusion
among all regions of boundary-parallelism for components of $S$, then
$X\cap S=S_2$.  We must have $X\subset A_j$ for some $j\in\{0,1\}$.
The hypotheses of Lemma \ref{fixupthesurfacetoo} clearly hold with
this choice of $X$, and with $A=A_j$. (See Remark
\ref{degeneratecase}(1).)  Hence \ref{fixupthesurfacetoo} implies
that $S'=S-S_2$ is a \sigmasurface.  If $\partial S_2\ne\emptyset$
then $b(S')<b(S)$.  If $\partial S_2=\emptyset$ then $b(S')=b(S)$,
$\chiminus(S')\le\chiminus(S)$, $t(S')\le t(S)$ and $s(S')=s(S)-1$.
In either case we conclude that $c(S')<c(S)$, a contradiction. This
completes the proof that $S$ is essential.

We now turn to the proof that $f$ has no folded boundary-annuli or big
folded $I$--bundles. First suppose that $f$ has a folded
boundary-annulus; that is, there is an annulus $R\subset\partial M$
such that $\inter R$ is a component of $\partial M - \partial S$, and
for some component $\tilde R$ of $p^{-1}(R)$, the components of
$\partial\tilde R=\tilde R\cap p^{-1}(E_T)$ are mapped by $f$ to the
same point of $E_T$. We must have $R\subset A_j$ for some
$j\in\{0,1\}$. Let $X$ be a regular neighborhood of $R$ in $A_j$ such
that $X\cap S$ is a regular neighborhood of $\partial R$ in $S$. Thus
$X\cap S$ consists of two disjoint annuli $R_1$ and $R_2$, while
$\frontier_AX$ is an annulus $R_3$. The hypotheses of Lemma
\ref{fixupthesurface} clearly hold with this choice of $X$, and
with $A=A_j$. Hence \ref{fixupthesurface} implies that the ``tubed''
surface $S'=(S-(R_1\cup R_2))\cup R_3$ is a \sigmasurface. But we
have $b(S')<b(S)$ and hence $c(S')<c(S)$, a contradiction.

Finally, suppose that $f$ has a big folded $I$--bundle; that is, there
is a submanifold $X$ of $M$ which is an $I$--bundle over a compact,
connected surface of negative Euler characteristic, such that (i)
$Y=X\cap S$ is the associated $\partial I$--bundle of $X$, (ii) $Y$ is
$\pi_1$--injective in $S$, and (iii) for some component $\tilde X$ of
$p^{-1}(X)$, the components of $\tilde X\cap p^{-1}(Y)=\tilde X\cap
p^{-1}(E_T)$ are mapped by $f$ to the same point of $E_T$. We may
choose $X$ so that $Y\subset\inter S$.  The surface ${\mathcal R}=\partial
X-\inter Y$ is a (possibly empty) disjoint union of annuli. We must
have $X\subset A_j$ for some $j\in\{0,1\}$. The hypotheses of Lemma
\ref{fixupthesurface} clearly hold with this choice of $X$, and with
$A=A_j$. Hence \ref{fixupthesurface} implies that
$S'=(S-Y)\cup {\mathcal R}$ is a \sigmasurface.

Since the base of the $I$--bundle $X$ has negative Euler
characteristic, so does every component of $Y$. Let $S_3$ denote the
union of all components of $S$ that meet $Y$. (There are at most two
such components.) The $\pi_1$--injectivity of $Y$ in $S$ implies that
each component of $S_3$ has negative Euler characteristic, so that
$\chiminus(S_3)=-\chi(S_3)$. Set $S_3'=(S_3-Y)\cup {\mathcal R}$. Then
each component of $S_3'$ contains a component of $\partial Y$. (This
is vacuously true if $Y$ is closed, since $S_3'=\emptyset$ in that
case.)  Since $Y$ is $\pi_1$--injective in the essential surface $S$,
each component of $\partial Y$ is homotopically non-trivial in $M$ and
hence in $S_3'$. This shows that each component of $S_3'$ is
non-simply-connected and hence has non-positive Euler
characteristic. Thus
$$\chiminus(S_3')=-\chi(S_3')=-(\chi(S_3)-\chi(Y)+\chi({\mathcal
R}))<-\chi(S_3),$$ since $\chi({\mathcal R})=0$ and $\chi(Y)<0$. Hence
$\chiminus(S_3')<\chiminus(S_3)$, which implies
$\chiminus(S')<\chiminus (S)$. As it is clear that $b(S')=b(S)$, it
follows that $c(S')<c(S)$, and again we have a contradiction.
\EndProof

\Lemma\label{whydeltaequalslambda}
Suppose that $M$ is an irreducible knot manifold, and that $T$ is a
bipartite $\pi_1(M)$--tree.  Let $(\tilde M,p)$ denote the universal
covering of $M$ and suppose that $f\co \tilde M \to T$ is a transverse
equivariant map which has no folded boundary-annuli. Then for any
closed curve $c$ in $\partial M$ we have
$\lambda_T(c)=\Delta_{\partial M}(c,\partial S_f)$.
\EndLemma

\Proof
We fix a standard action of $\pi_1(M)$ on $\tilde M$ that makes $f$ a
$\pi_1(M)$--equivariant map.

Since $T$ is a bipartite $\pi_1(M)$--tree, Proposition
\ref{bipartitemanifold} asserts that there are compact
$3$--dimensional submanifolds $A_0$ and $A_1$ of $M$ such that
$A_0\cap A_1=\frontier A_0=\frontier A_1=S_f$.

Set $\Delta=\Delta_{\partial M}(c,\partial S_f)$.  According to
Proposition \ref{Deltaandtranslation} we have
$$\lambda_T(c)\le\Delta_M(c,S_f)\le\Delta .$$
It remains to show
$\lambda_T(c) \ge \Delta$.  We may therefore assume that
$\Delta>0$. In particular, $\partial S_f\ne\emptyset$ and $c$ is
homotopically non-trivial. Since $\partial S_f\ne\emptyset$, all the
components of $A_0\cap\partial M$ and $A_1\cap\partial M$ are
annuli. We shall refer to these as {\it complementary annuli.} We may
suppose $c$ to be chosen within its homotopy class so that it is
transverse to $\partial S_f$ and $\#(c^{-1}(\partial
S_f))=\Delta$. This implies that for each component $\alpha$ of
$c^{-1}(A_0)$ or $c^{-1}(A_1)$, $c\vert_\alpha$ is a map of the arc
$\alpha$ into a complementary annulus $R$ which takes the endpoints of
$\alpha$ to different components of $\partial R$.

After reparametrization we may assume that
$$c^{-1}(\partial S_f)=\{\exp(2k\pi\sqrt{-1}/\Delta):0\le
  k<\Delta\}\subset S^1 .$$ 
Let $q\co \mathbb{R}\to S^1$ be the covering map defined by $q(t)=\exp\left(2\pi
\left(t-\frac{1}{2}\right)\sqrt{-1}/\Delta\right)$, set $\ell=c\circ
q\co \mathbb{R}\to\partial
M$ and let $\tilde \ell\co \mathbb{R}\to\tilde M$ denote a lift of $\ell$.
Then for every $t\in\mathbb{R}$ we have
$\tilde\ell(t+\Delta)=\gamma\cdot\tilde\ell(t))$, where $\gamma$ is
an element of $[c]$.

Our parametrization of $c$ guarantees that
$\tilde\ell^{-1}(p^{-1}(\partial S_f))=\mathbb{Z}+\frac{1}{2}$. For each
$n\in\mathbb{Z}$ we let $\tilde S_n$ denote the component of $p^{-1}(S_f)$
containing $\tilde\ell\left(n+\frac{1}{2}\right)$,
and we let $\tilde A_n$ denote the component of $p^{-1}(A_0)$ or $p^{-1}(A_1)$
containing $\tilde\ell\left(\left[n-\frac{1}{2},n+\frac{1}{2}\right]\right)$.  We denote by
$\mu_n\in E_T$ the point $f(\tilde S_n)$, by $e_n$ the edge of $T$
whose midpoint is $\mu_n$, and by $s_n$ the unique vertex of $T$ in
the component of $|T|-E_T$ containing $f(\tilde A_n)$. The
transversality of $f$ to $E_T$ implies that the endpoints of $e_n$ are
$s_n$ and $s_{n+1}$ for every $n\in\mathbb{Z}$.

For every $n\in\mathbb{Z}$ we have
$$\mu_{n+\Delta}\!=\!
  f\left(\tilde\ell\left(n+\frac{1}{2}+\Delta\right)\right)\!=\!
  f\left(\gamma\cdot\tilde\ell\left(n+\frac{1}{2}\right)\right)\!=\!
  \gamma\cdot f\left(\tilde\ell\left(n+\frac{1}{2}\right)\right)\!=\!
  \gamma\cdot \mu_n.$$
Hence $e_{n+\Delta}=\gamma\cdot e_n$ for every $n$.

For each $n$ the interval $\left[n-\frac{1}{2},n+\frac{1}{2}\right]$ is mapped by
$\ell$ into a complementary annulus $R_n$, and
$\ell\left(n-\frac{1}{2}\right)$ and
$\ell\left(n+\frac{1}{2}\right)$ lie in different components of $\partial R_n$.
Hence $\tilde\ell$ maps $\left[n-\frac{1}{2},n+\frac{1}{2}\right]$
into a component
$\tilde R_n\subset \tilde A_n$ of $p^{-1}(R_n)$, and the points
$\ell\left(n-\frac{1}{2}\right)$ and $\ell\left(n+\frac{1}{2}\right)$
lie in different
components $\tilde C_{n-1}\subset\partial \tilde S_{n-1}$ and $\tilde
C_{n}\subset\partial \tilde S_{n}$ of $\partial \tilde R_n$.  But
$\partial \tilde R_n$ must have exactly two components, and the
hypothesis that $f$ has no folded boundary-annuli implies that $f$
maps these components to different points of $E_T$. Hence $\mu_{n-1}
\ne \mu_n$, which implies that $e_{n-1}\ne e_n$ for every $n\in\mathbb{Z}$.  In particular $s_n$ is the unique common vertex of $e_{n-1}$ and
$e_n$. Since $e_{n+\Delta}=\gamma\cdot e_n$ for every $n$ it now
follows that $s_{n+\Delta}=\gamma\cdot s_n$ for every $n$.

Since $T$ is a tree, and since $e_{n-1}\ne e_n$ for every $n$, the
$e_n$ and $s_n$ make up a subgraph $\mathcal A$ of $T$ isomorphic to the
real line, triangulated with a vertex at every integer point. In
particular, for all $m,n\in\mathbb{Z}$, we have $d_T(s_m,s_n)=|m-n|$, and
the arc joining $s_m$ and $s_n$ is contained in $\mathcal A$. Hence for
any $n$ we have $d_T(s_n,\gamma\cdot
s_n)=d_T(s_n,s_{n+\Delta})=\Delta$. Now consider an arbitrary vertex
$s$ of $T$, let $s_n$ be a vertex of $\mathcal A$ for which $d_T(s,s_n)$
is as small as possible, let $\beta$ denote the arc with endpoints $s$
and $s_n$, and let $\alpha\subset{\mathcal A}$ denote the arc with
endpoints $s_n$ and $\gamma\cdot s_n=s_{n+\Delta}$. Then
$\beta\cap{\mathcal A}=\{s_n\}$ and hence $(\gamma\cdot \beta)\cap{\mathcal
A}=\gamma\cdot(\beta\cap{\mathcal A})=\{s_{n+\Delta}\}$.
In particular,
$\beta\cap{\alpha}= \{s_{n}\}$ and $(\gamma\cdot\beta)\cap{\alpha}=
\{s_{n+\Delta}\}$. Since $T$ is a tree it follows that
$d(s,\gamma\cdot s)=\Delta+2d_T(s,s_n)\ge\Delta$. This proves that
$\lambda_T(c)\ge\Delta$, as required.
\EndProof

\Lemma\label{whylongarcisfixed}
Suppose that $M$ is an irreducible knot manifold, and that $T$ is a
bipartite $\pi_1(M)$--tree.  Let $(\tilde M,p)$ denote the universal
covering of $M$, and suppose that $f\co \tilde M \to T$ is a transverse
equivariant map such that $S_f$ is essential and $f$ has no big folded
$I$--bundle. If $K\subset S_f$ is any $\pi_1$--injective, connected,
compact polyhedron such that $\chi(K)<0$, if $\theta\le t_{S_f}(K)$ is
a positive integer (cf \ref{thicknessdef}), and if
$\Theta\le\pi_1(M)$ is the subgroup defined up to conjugacy by
$\Theta=\image(\pi_1(K)\to\pi_1(M))$, then $\Theta$ fixes an arc of
length $\theta$ in $T$.
\EndLemma

\Proof
We choose a base point $\star$ in $K\subset M$ and a standard action of
$\pi_1(M) = \pi_1(M,\star)$ on $\tilde M$ that makes $f$ a
$\pi_1(M)$--equivariant map.

According to \ref{bipartitemanifold}, there are closed subsets $A_0$
and $A_1$ of $M$ which are $3$--dimensional submanifolds, such that
$A_0\cap A_1=\frontier A_0=\frontier A_1=S_f$. Since $\theta\le
t_{S_f}(K)$, there is a reduced homotopy $H\co (K\times I,K\times\partial
I)\to(M,S_f)$ of length $\theta-1$ such that, for some $t\in I$, the map
$H_t$ is the inclusion map $K\hookrightarrow S_f\subset M$. By
definition $H$ is a composition of essential basic homotopies
$H^1,\ldots,H^{\theta-1}$, and by symmetry we may assume that $H^i(K\times
I)\subset A_{[i]}$, where $[i]$ denotes the least residue of $i$
modulo $2$. Thus there are points $0=t_0<t_1<\cdots<t_{\theta-1}=1$ of
$I$ such that $H\vert_{K\times[t_{i-1},t_i]}$ is a reparametrization of $H^i$
for $i=1,\ldots,\theta-1$, and there is some
$m\in\{0,\ldots,\theta-1\}$ such that $H_{t_m}$ is the inclusion
$K\hookrightarrow M$. We let $\xi\co \pi_1(K,\star)\to\pi_1(M,\star)$ denote the
inclusion homomorphism.

Let $(\tilde K, q)$ denote the universal covering of $K$, set
$h=H\circ(q\times\id)\co \tilde K\times I\to M$, and choose a lift
$\tilde h\co \tilde K\times I\to \tilde M$ of $h$. Note that with respect to
our chosen standard action on $\tilde M$ and some standard
action of $\pi_1(K)$ on $\tilde K$, the map $\tilde h$ is
$\xi$--equivariant in the sense that for every $(\tilde
z,t)\in\tilde K\times I$ and every $\gamma\in\pi_1(K)$ we have $\tilde
h(\gamma\cdot \tilde z,t)=\xi(\gamma)\cdot h(\tilde z,t)$.

For $i=0,\ldots,\theta-1$, let $\tilde S_i$ denote the component of
$p^{-1}(S_f)$ containing $\tilde h(\tilde K\times\{t_i\})$. Then $f(\tilde
S_i)$ is a point of $E_T$, which means that it is the midpoint of a
well-defined edge $e_{i+1}$ of $T$.

For $i=0,\ldots,\theta-1$, since $\tilde h$ maps $\tilde
K\times\{t_i\}$ into the component $\tilde S_i$ of the
$\pi_1(M)$--invariant set $p^{-1}(S_f)\subset\tilde M$, the
$\xi$--equivariance of $\tilde h$ implies that $\tilde S_i$ is
invariant under the subgroup $\Theta=\xi(\pi_1(K))$ of $\pi_1(M)$.
The $\pi_1(M)$--equivariance of $f$ then implies that the midpoint of
the edge $e_{i+1}$ is fixed by $\Theta$ for
$i=0,\ldots,\theta-1$. Since the bipartite $\pi_1(M)$--tree $T$ is in
particular a $\pi_1(M)$--tree without inversions (see
\ref{bipartitedef}), the edge $e_{i+1}$ is itself fixed by $\Theta$
for $i=0,\ldots,\theta-1$.

For $i=1,\ldots,\theta-1$, let $\tilde A_i$ denote the component of
$p^{-1}(A_{[i]})$ containing $\tilde h(\tilde
K\times[t_{i-1},t_i])$. Then $f(\tilde A_i)$ is contained in the
closure of a unique component of $T - E_T$.  This component contains a
unique vertex of $T$ which will be denoted $s_i$. It is clear that
$s_i$ is a common endpoint of $e_{i-1}$ and $e_i$ for
$i=1,\ldots,\theta-1$. We denote by $s_0$ the vertex of $e_1$ which is
distinct from $s_1$, and by $s_\theta$ the vertex of $e_{\theta}$
which is distinct from $s_{\theta-1}$.  We denote by $\omega_i$ the
orientation of $e_i$ such that $\init(\omega_i) = s_i$ and
$\term(\omega_i) = s_{i+1}$.  Then
$(\omega_1,\ldots,\omega_{\theta-1})$ is an edge path in $T$.  We
claim that this edge path is reduced.  This amounts to showing that
$e_i$ and $e_{i-1}$ are distinct for any $i$ with $1 \le i \le
\theta$.

To prove this, we note that since $H^i$ is an essential homotopy in
$A_{[i]}$, it follows from
\cite[``Essential Homotopy Theorem'' Chapter III \S 2]{charsub}
that $H^i\co (K\times I,K\times\partial I)\to (A_{[i]},S_f)$ is homotopic
as a map of pairs to a map $J^i$ such that $J^i(K\times I)\subset X_i$
and $J^i(K\times\partial I)\subset Y_i$, where $X_i$ is a submanifold
of $A_{[i]}$, $Y_i$ is a submanifold of $X_i\cap S_f$ which is
$\pi_1$--injective in $S_f$, and either (i) $X_i\subset A_{[i]}$ is an
$I$--bundle over a surface and $Y_i\subset S_f$ is the associated
$\partial I$--bundle, or (ii) $X_i$ is a Seifert fibered space and
$Y_i$ is a saturated subsurface of $\partial X_i$. On the other hand,
since $H_0\co K\to M$ is homotopic in $M$ to the inclusion
$K\hookrightarrow S_f\subset M$, and since $K$ is $\pi_1$--injective in
the essential surface $S_f$ and $\chi(K)<0$, the subgroup
$(J^i_0)_\sharp(\pi_1(K))=(H^i_0)_\sharp(\pi_1(K))$, defined up to
conjugacy in $\pi_1(M)$, is non-abelian. Hence $Y_i$ has a component
with non-abelian fundamental group. This rules out (ii), and shows
that the base of the $I$--bundle given by (i) must have negative Euler
characteristic. Note that since $H^i$ is an essential homotopy, the
homotopy $J^i\co (K\times I,K\times\partial I)\to(X_i,Y_i)$ is also
essential.

Since $H\vert_{K\times[t_{i-1},t_i]}$ is a reparametrization of the
homotopy $H_i$, it follows that the map
$h\vert_{K\times[t_{i-1},t_i]}$ is a reparametrization of
$h^i=H^i\circ(q\times\id)\co \tilde K\times I\to M$.  Hence $\tilde
h\vert_{\tilde K\times[t_{i-1},t_i]}$ is a reparametrization of a lift
$\tilde h^i\co \tilde K\times I\to\tilde M$ of $h^i$.  Thus $\tilde
h^i(\tilde K\times I)\subset \tilde A_{i}$, $\tilde h^i(\tilde
K\times\{0\})\subset\tilde S_{i-1}$, and $\tilde h^i(\tilde
K\times\{1\})\subset\tilde S_{i}$.  Since $H^i\co (K\times
I,K\times\partial I)\to (A_{[i]},S_f)$ is homotopic to $J^i$ as a map
of pairs, the covering homotopy property of covering spaces implies
that $\tilde h^i\co (\tilde K\times I,\tilde K\times\partial I)\to
(\tilde A_i,\tilde A_i\cap p^{-1}(S_f))$ is homotopic as a map of
pairs to some lift $\tilde j^i$ of $j^i=J^i\circ(q\times\id)$. In
particular it follows that $\tilde j^i(\tilde
K\times\{0\})\subset\tilde S_{i-1}$ and that $\tilde j^i(\tilde
K\times\{1\})\subset\tilde S_{i}$.

Let $\tilde X_i$ denote the component of $p^{-1}(X_i)$ containing
$\tilde j^i(\tilde K\times I)$, and set $\tilde Y_i=\tilde X_i\cap
p^{-1}(Y_i)$, so that $\tilde j^i(\tilde K\times I)\subset\tilde Y_i$. 
Since $f$ has no big folded $I$--bundles, it follows from Remark
\ref{whatnobigfoldedIbundlesmeans} that $\tilde X_i$ is a trivial
$I$--bundle with associated $\partial I$--bundle $\tilde Y_i$, and that
$f$ maps the two components of $\tilde Y_i$ to distinct points of $E_T$.

Now consider any point $\tilde z\in \tilde K$, and set $z=q(\tilde
z)\in K$.  Let $w^i_0=\tilde j^i(\tilde z,0)$ and $w^i_1=\tilde
j^i(\tilde z,1)$.  If $w^i_0$ and $w^i_1$ lie in the same component of
$\tilde Y_i$, then the path $t\mapsto\tilde j^i(\tilde z,t)$ in the
trivial $I$--bundle $\tilde X_i$ is fixed-endpoint homotopic to a path
in $\tilde Y_i$. This implies that the path $t\mapsto J^i( z,t)$ is
fixed-endpoint homotopic in $X_i$ to a path in $ Y_i$, which is
impossible since $J^i\co (K\times I,K\times\partial I)\to(X_i,Y_i)$ is an
essential homotopy. Hence $w^i_0$ and $w^i_1$ lie in different
components of $\tilde Y_i$, and hence $f(w^i_0)$ and $f(w^i_1)$ are
distinct points of $E_T$. On the other hand, we have $w^i_0\in\tilde
j^i(\tilde K\times\{0\})\subset\tilde S_{i-1}$ and $w^i_1\in\tilde
j^i(\tilde K\times\{1\})\subset\tilde S_{i}$, which implies that
$f(w^i_0)$ and $f(w^i_1)$ are the midpoints of $e_{i-1}$ and $e_i$
respectively. This shows that $e_{i-1}$ and $e_i$ are distinct edges
of $T$, and establishes the claim that the edge path
$(\omega_1,\ldots,\omega_{\theta-1})$ is reduced.

Since $T$ is a tree, it now follows that the edges
$e_1,\ldots,e_\theta$ form an arc of length $\theta$ in $T$. As we
have seen that $\Theta$ fixes $e_1,\ldots,e_\theta$, we have now
produced the required arc of length $\theta$ fixed by $\Theta$.
\EndProof

\Proof[Proof of Theorem \ref{dualsurface}]
The theorem is an immediate consequence of Lemmas
\ref{controlledfolding}, \ref{whydeltaequalslambda}, and
\ref{whylongarcisfixed}.
\EndProof

\section{The tree for $\gltwo$}

In this section we record a few facts about the tree for $\gltwo$ over
a discretely valued field. Our point of view is close to that of Serre
\cite{serre}, except that to be consistent with the conventions of
\ref{graphdef} and \ref{treedef} we take the tree to be a
$1$--connected geometric simplicial $1$--complex which realizes the
abstract combinatorial structure considered by Serre. We begin by
summarizing some results from \cite{serre}, translated into our
geometric setting.

\Number
Suppose that $F$ is a field with a discrete rank--$1$ valuation $v$.
We always denote the valuation ring associated to $v$ by ${\mathcal O}_v$;
it consists of all elements $x\in F$ with $v(x)\ge0$, where by
convention $v(0)=+\infty$. A {\it lattice} in the
$2$--dimensional vector space $F^2$ is a rank--$2$ ${\mathcal O}$--submodule
of $F^2$.

There is a $\gltwo(F)$--tree, in the sense of Section
\ref{dualsurfacesection}, canonically associated to the valued
field $F$.  We shall always denote this tree by $T_F$, leaving the
valuation $v$ implicit in the notation.  The vertices of $T_F$ are in
bijective correspondence with homothety classes of lattices in
$F^2$. If $L$ is a lattice representing a vertex $s$ of $T_F$, and if
$\pi\in{\mathcal O}$ is a uniformizer (ie, an element such that
$v(\pi)=1$), then any vertex $s'$ can be represented by a lattice
$L'\subset L$ which is generated by $e$ and $\pi^df$ for some integer
$d\ge0$ and some basis $\{e,f\}$ of $L$.  The integer $d$, which is
uniquely determined by the vertices $s$ and $s'$, is equal to the
distance $d_{T_F}(s,s')$; this fact completely characterizes the tree
$T_F$, since two vertices $s,s'$ are joined by an edge if and only if
$d_{T_F}(s,s')=1$. The action of $\gltwo(F)$ on $T_F$, which is
transitive on the vertices, is characterized by the fact that an
element $A\in\gltwo(F)$ carries the vertex represented by a lattice
$L$ to the vertex represented by $A(L)$.
\EndNumber

\Proposition\label{sltwobipartite}
If $F$ is a field with a discrete rank--$1$ valuation $v$, the
$\sltwo(F)$--tree $T_F$ is bipartite. Furthermore, if an element $A$ of
$\sltwo(F)$ fixes a vertex of $T_F$ represented by a lattice $L$, then
$A(L)=L$.
\EndProposition

\Proof
Let $s$ be any vertex of ${T_F}$, and $L$ a lattice representing
$s$. Let $A\in\sltwo(F)$ be given, and set $s'=A\cdot s$ and
$d=d_{T_F}(s,s')$. Then $s'$ is represented by a lattice $L'\subset L$
which is generated by $e$ and $\pi^df$ for some basis $\{e,f\}$ of
$L$. As $A(L)$ also represents $s'$, we must have $L'=\pi^kA(L)$ for
some $k\in\mathbb{Z}$. Hence if $B$ is the element of $\gltwo(F)$ defined
by $B(e)=\pi^{-k}e$ and $B(f)=\pi^{d-k}(f)$, we have $B(L)=A(L)$. Thus
$A^{-1}B$ leaves $L$ invariant, and $\det(A^{-1}B)=\det B$ must be a
unit in $\mathcal O$, ie, $v(\det B)=0$. But $B$ is conjugate in
$\gltwo(F)$ to
$$\begin{pmatrix}
\pi^{-k}&0\\0&\pi^{d-k}
\end{pmatrix}
,$$ 
so that $v(\det B)=d-2k$. Hence $d=2k$. In particular $d$ is always
even, so that $T_F$ is a bipartite $\sltwo(F)$--tree.

Now if $A\cdot s=s$, so that $d=0$, then $k=0$ and hence
$A(L)=L'=L$.
\EndProof

\Proposition\label{valuationandtranslation}
Suppose that $F$ is a field with a discrete rank--$1$ valuation
$v$. Then for every $A\in\sltwo(F)$ we have
$$\lambda_{T_F}(A)=2\max(0,-v(\trace(A))).$$
\EndProposition
 
\Proof
We set ${\mathcal O}={\mathcal O}_v$, and we denote by $s_0$ the
homothety class of the standard lattice ${\mathcal O}^2\subset F^2$. We
set $v_0=\max(0,-v(\trace(A)))$.  The proposition asserts that
$\lambda_{T_F}(A)=2v_0$.  We first show that
$\lambda_{T_F}(A)\ge2v_0$; for this step we may assume without loss of
generality that $v_0>0$. We need to prove that for any vertex $s$ of
$T_F$ we have $d_{T_F}(s,A\cdot s)\ge2v_0$.  Since $\gltwo(F)$ acts
transitively on the vertices of $T_F$, and since the length function
$\lambda_{T_F}$ is constant on $\gltwo(F)$--conjugacy classes, it
suffices to show that for any conjugate $B$ of $A$ in $\gltwo(F)$ we
have $d_{T_F}(s_0,B\cdot s_0)\ge2v_0$.

Since a rank--$1$ discrete valuation ring is Euclidean, we may reduce
the matrix $B$ to a diagonal matrix $D$ using row and column operations over
$\mathcal O$: a column operation  over $\mathcal O$ has the form
$$\begin{pmatrix}a&b\\c&d\end{pmatrix}
  \to 
  \begin{pmatrix}a&b+\alpha a\\c&d+\alpha c\end{pmatrix}
  \qquad{\rm or}\qquad 
  \begin{pmatrix}a&b\\ c&d\end{pmatrix}
  \to 
  \begin{pmatrix}a+\alpha b&b\\ c+\alpha d&d\end{pmatrix}$$
for some $\alpha\in{\mathcal O}$, and row operations over $\mathcal O$ are
defined similarly. In particular we have $D=XBY$ for some
$X,Y\in\sltwo({\mathcal O})$, and hence $D\in\sltwo(F)$. Thus
$$D=\begin{pmatrix}\delta&0\\0&\delta^{-1}\end{pmatrix}$$
for some $\delta\in{F}$.

Let us define the {\it height} of an arbitrary matrix $M\in\sltwo(F)$
to be the minimum of $v(\alpha)$, where $\alpha$ ranges over the
entries of $M$.  Note that $v(\trace B)=v(\trace A)=-v_0$ (since
$v_0>0$), and hence $v(\beta)\le-v_0$ for at least one diagonal entry
$\beta$ of $B$.  Hence $B$ has height at most $-v_0$. It is clear that
a row or column operation defined over $\mathcal O$ does not affect the
height of a matrix, and hence $\height D=\height B\le-v_0$. This means
that either $v(\delta)\le- v_0$ or that $v(\delta^{-1})=-v(\delta)\le
-v_0$, so in either case $|v(\delta)|\ge v_0>0$. Now $B({\mathcal O}^2)
=X^{-1}DY^{-1}({\mathcal O}^2)
=X^{-1}D({\mathcal O}^2)
=X^{-1}(\delta{\mathcal O}\oplus\delta^{-1}{\mathcal O})
= \delta{\mathcal O}e+\delta^{-1}{\mathcal O}f$,
where $\{e,f\}$ is the image of the standard basis for ${\mathcal O}^2$
under $X^{-1}$, and is itself a basis of ${\mathcal O}^2$ since
$X^{-1}\in\sltwo({\mathcal O})$. We define a lattice $L\subset{\mathcal O}^2$
by $L=\delta^2{\mathcal O} e+{\mathcal O}f$ if $v(\delta)>0$ and by $L=\{{\mathcal
O} e+\delta^{-2}{\mathcal O}f\}$ if $v(\delta)<0$. In either case $L$ is
homothetic to $B({\mathcal O}^2)$ and hence represents the vertex $B\cdot
s_0$, and the definition of distance in $T_F$ implies that
$d_{T_F}(s_0,B\cdot s_0)=|v(\delta^2)|=2|v(\delta)|\ge2v_0$, as
required.

It remains to show that $\lambda_{T_F}(A)\le2v_0$. This is
trivial if $A=\pm Id$. If $A\ne\pm Id$, we may assume after a
conjugation in $\gltwo(F)$ that
$$A=\begin{pmatrix}0&-1\\1&\tau\end{pmatrix}$$ for some $\tau\in F$, and it is apparent
that $\tau=\trace A$. If $v(\tau)\ge0$ then $A\in\sltwo(F)$; hence $A$
fixes $s_0$, so that $\lambda_{T_F}(A)=0\le2v_0$. Now suppose that
$v(\tau)<0$.  If $s_1$ denotes the homothety class of the lattice
$L_1$ generated by $(1,0)$ and $(0,\tau^{-1})$, then $A\cdot s_1$ is
represented by the lattice $\tau^{-1} A(L_1)\subset L_1$, which is
generated by $(\tau^{-2},0)$ and $(0,\tau^{-1})$. The definition of
distance in $T_F$ implies that $d_{T_F}(s_1,A\cdot
s_1)=v(\tau^{-2})=-2v(\tau)=2v_0$. Hence $\lambda_{T_F}(A)\le2v_0$.
\EndProof

\Proposition\label{valuationandfixedarc}
Suppose that $F$ is a field with a discrete rank--$1$ valuation
$v$. Suppose that $J$ is a subgroup of $\sltwo(F)$ which fixes an arc
of length $t$ in $T_F$.  Then for every element $A$ of the commutator
subgroup $[J,J]\le J\le\sltwo(F)$, we have
  $$v((\trace(A))-2)\ge t.$$
\EndProposition

\Proof
We set ${\mathcal O}={\mathcal O}_v$. The hypothesis implies that
there are vertices $s$ and $s'$ of $T_F$ such that $d_{T_F}(s,s')=t$,
$H\cdot s=s$ and $H\cdot s'=s'$. Since $d_{T_F}(s,s')=t$, the vertices
$s$ and $s'$ are represented by lattices $L$ and $L'$ which
respectively have bases of the forms $\{e,f\}$ and $\{e,\pi^tf\}$.
After conjugating by an element of $\gltwo(F)$ we may assume that
$\{e,f\}$ is the standard basis for $F^2$. This implies that
$L={\mathcal O}^2$ and that $L'=C({\mathcal O}^2)$, where
$$C=\begin{pmatrix}1&0\\0&\pi^t\end{pmatrix}\in\gltwo(F).$$
Hence $s'=C\cdot s$. It follows that the subgroups $J$ and $C^{-1}JC$
of $\sltwo(F)$ both fix the vertex $s$. Hence by the second assertion
of Proposition \ref{sltwobipartite}, they are both contained in
$\sltwo({\mathcal O})$.
$$X=\begin{pmatrix}a&b\\c&d\end{pmatrix}\in J,\leqno{\hbox{For any}}$$ 
$$C^{-1}XC=\begin{pmatrix}a&\pi^tb\\\pi^{-t}c&d\end{pmatrix}.\leqno{\hbox{we have}}$$
As $X$ and $C^{-1}XC$ both belong to $\sltwo({\mathcal O})$, we have
$a,b,d\in\mathcal O$ and $c\in\pi^{t}{\mathcal O}$. It follows that
the natural homomorphism $\eta\co \sltwo({\mathcal O})$ to
$\sltwo({\mathcal O}/\pi^t{\mathcal O})$ maps $J$ onto a group $\bar
J$ of upper triangular matrices. Hence $\eta([J,J])\le[\bar J,\bar J]$
is a subgroup of $\sltwo({\mathcal O}/\pi^t{\mathcal O})$ consisting
of upper triangular matrices whose diagonal entries are equal to
$1$. In particular the trace of any element of $\eta([J,J])$ is equal
to $2$. This means that for any $A\in [J,J]$ we have $\trace
A\in2+\pi^t{\mathcal O}$, so that $v((\trace A)-2)\ge t$, as asserted.
\EndProof

\section{Curves, norms and actions associated to ideal\nl points}
\label{actionsfromidealpoints}

\vspace{-12pt}

\Number\label{sigmadef}
In this subsection we review notation for character varieties as used
in \cite{cgls}, and introduce some additional notation that will be
needed in this paper. 

We begin with some algebraic geometric conventions.  Suppose that $K$
denotes the field of rational functions on an irreducible complex
projective algebraic curve $C$ and that $x$ is a smooth point of $C$.
For a non-zero element $f$ of $K$ we will write $Z_x(f)$ to denote the
order of zero of $f$ at $x$, or $0$ if $f$ does not have a zero at
$x$.  Similarly we we will let $\Pi_x(f)$ denote the order of pole of
$f$ at $x$, or $0$ if $f$ does not have a pole at $x$.  The function
$v_x\co F^*\to \mathbb{Z}$ defined by $v_x(f) = Z_x(f) - \Pi_x(f)$ is a
discrete rank--1 valuation on the field $K$.  We will denote the
valued field $(K,v_x)$ by $K_x$.

Given an irreducible complex affine algebraic curve $A$ we will denote
by $\tilde A$ the unique smooth projective curve that admits a
birational correspondence $\phi\co \tilde A\to A$.  (The curve $\tilde A$
can be constructed by desingularizing a projective completion of $A$.)
We will say that a point $x\in \tilde A$ is an {\it ideal point} if it
does not correspond to any point of $A$ under $\phi$.

Now let $\Gamma$ be a finitely generated group.  We
will denote by $R(\Gamma)$ the complex affine algebraic set of
representations of $\Gamma$ in $SL_2(\mathbb{C})$.  If $R_0$ is an
irreducible subvariety of $R(\Gamma)$ and if $F$ denotes the field of
rational functions on $R(\Gamma)$ then the {\it tautological
representation} $P\co \Gamma \to \sltwo(F)$ associated to $R_0$ is defined
by
$$P(g) = \begin{pmatrix}a&b\\c& d\end{pmatrix}$$
where the functions $a$, $b$, $c$, and $d$ satisfy
$$\rho(g) = \begin{pmatrix}a(\rho)& b(\rho)\cr c(\rho)& d(\rho)\cr\end{pmatrix}.$$
We will denote by $X(\Gamma)$ the set of all characters of
representations in $R(\Gamma)$, and by $\tau\co R(\Gamma)\to X(\Gamma)$ the
surjective regular map such that $\tau(\rho)$ is the character of
$\rho$. We give $X(\Gamma)$ the structure of an affine algebraic set
as in \cite{cgls}.

Next suppose that $A$ is an irreducible affine algebraic curve
contained in $X(\Gamma)$ and let $x\in \tilde A$ be an ideal point.
Let $K$ denote the field of rational functions on $A$.  We use the
birational correspondence between $\tilde A$ and $A$ to identify $K$
with the function field of $\tilde A$, and we regard $v_x$ as a
discrete valuation on $K$.  Let $R_0$ denote an irreducible component
of $\tau^{-1}(A)$ which is mapped to a dense subset of $A$ by $\tau$.
We use $\tau$ to identify $K$ with a subfield of the function field
$F$ of $R_0$.  According to \cite[Theorem 1.2.3]{cgls},  we may
extend the valuation $v_x$ to a discrete valuation on $F$ and we will
denote by $F_x$ the resulting discretely valued field.  We consider
the $\gltwo(F_x)$--tree $T_{F_x}$ and we let $T_x$ denote the
$\Gamma$--tree which is the pull-back of $T_{F_x}$ under the
tautological representation $P\co \Gamma \to \sltwo(F_x) < \gltwo(F_x)$
associated to $R_0$.

If $\gamma$ is an element of $\Gamma$ then we shall let $I_\gamma$
denote the rational function on $A$ defined by $I_\gamma(\chi) =
\chi(\gamma)$.  Using the identifications described above, we shall regard
$I_\gamma$ as an element of $F_x$.  We recall from \cite[1.2.4]{cgls}
that $I_\gamma = \trace P(\gamma)$.  If $\mathcal C$ is the
conjugacy class of $\gamma\in\Gamma$ then $I_\gamma = I_{\gamma'}$ for
any element $\gamma'$ of ${\mathcal C}$.  We will write $I_{\mathcal C} =
I_\gamma$.
\EndNumber

\Proposition\label{poleandtranslation}
Let $\Gamma$ be any finitely generated group, let $A\subset
X(\Gamma)$ be any curve and let $x\in \tilde A$ be an ideal point.
Then for any element $\gamma$ of $\Gamma$ we have
  $$2\Pi_x(I_\gamma)=\lambda_{T_x}(\gamma).$$
\EndProposition

\Proof
Using the notation of \ref{sigmadef}, we have
$$\Pi_x(I_\gamma)=\max(0,-v_x(I_\gamma))=\max(0,-v_x(\trace(P(\gamma)))).$$
It therefore follows from Proposition \ref{valuationandtranslation}
that $2\Pi_x(I_\gamma)=\lambda_{T_{F_x}}(P(\gamma))$.  But since the
$\Gamma$--tree $T_x$ is the pull-back of the $\gltwo(F_x)$--tree
$T_{F_x}$ via the representation $P$ we have
$\lambda_{T_{F_x}}(P(\gamma))=\lambda_{T_x}(\gamma)$, and the
assertion follows. 
\EndProof

\Proposition\label{ActionNontrivial}
Let $\Gamma$ be any finitely generated group, let $A\subset
X(\Gamma)$ be any curve and let $x\in \tilde A$ be any ideal point.
Then the tree $T_x$ is a non-trivial bipartite $\Gamma$--tree.
\EndProposition

\Proof
We again use the notation of \ref{sigmadef}.
Since the $\sltwo(F_x)$--tree $T_{F_x}$ is bipartite by
\ref{sltwobipartite}, it follows that the $\Gamma$--tree $T_x$ is
also bipartite.

By definition the ideal point $x$ does not correspond to a point of
the affine curve $A$ under the birational correspondence between
$\tilde A$ and $A$.  Hence there is an element of the coordinate ring
of $A$ which, when regarded as a function on $\tilde A$, has a pole at
$x$.  Since the functions $I_\gamma$ generate the coordinate ring
there exists $\gamma\in \Gamma$ such that $\Pi_x(I_\gamma)>0$.
Proposition \ref{poleandtranslation} implies that $\gamma$ has no
fixed vertex in $T_x$ and hence that $T_x$ is a non-trivial
$\Gamma$--tree.
\EndProof

\Proposition\label{commutatorandidealpoint}
Let $\Gamma$ be any finitely generated group, let $A\subset
X(\Gamma)$ be any curve, let $x\in \tilde A$ be an ideal point.
Suppose that $\Theta$ is a subgroup of $\Gamma$ which fixes an arc of
length $t>0$ in $T_x$. Then for every element $\gamma$ of the
commutator subgroup $[\Theta,\Theta]\le\Theta\le\Gamma$, we have
  $$Z_x(I_\gamma-2)\ge t.$$
\EndProposition

\Proof
Using the notation of \ref{sigmadef}, we have
\Equation\label{commutatorandidealpointone}
Z_x(I_\gamma-2)=\max(0,v_x(I_\gamma)-2)=\max(0,v_x((\trace
P(\gamma))-2))
\EndEquation
for any $\gamma\in\Gamma$.  Now suppose that
$\gamma\in[\Theta,\Theta]$, where $\Theta\le\Gamma$ fixes an arc of
length $t>0$ in $T_x$. Set $J=P(\Theta)$, so that $P(\gamma)\in[J,J]$.

Since the $\Gamma$--tree $T_x$ is the pull-back of the
$\gltwo(F_x)$--tree $T_{F_x}$ via the representation $P$ the group $J$
fixes an arc of length $t$ in $T_{F_x}$.  Hence by Proposition
\ref{valuationandfixedarc} we have
\Equation\label{commutatorandidealpointtwo}
v_x((\trace(P(\gamma)))-2)\ge t .
\EndEquation
The conclusion follows from (\ref{commutatorandidealpointone}) and 
(\ref{commutatorandidealpointtwo}).
\EndProof

\Number\label{principaldef}
If $M$ is a hyperbolic knot manifold, we define a {\it principal
component} $X_0$ of the character variety $X(\pi_1(M))$ to be a
component that contains the character of a discrete, faithful
representation of $\pi_1(M)$.
\EndNumber

\Lemma\label{nonconstant}
Let $M$ be a hyperbolic knot manifold and let $X_0$ be a principal
component of $X(\pi_1(M))$.  If $\gamma$ is a non-trivial peripheral
element of $\pi_1(M))$ then the function $I_\gamma\vert_{X_0}$ is
non-constant.  Furthermore if $\gamma$ is any non-trivial element of
$\pi_1(M)$ then the function $I_\gamma\vert_{X_0}$ cannot be
identically equal to $2$.
\EndLemma

\Proof
If $\gamma$ is a non-trivial peripheral element of
$\pi_1(M)$ then it follows from \cite[Proposition 3.2.1]{boundsep}
that $I_\gamma\vert_{X_0}$ is non-constant.

Now suppose that $\gamma$ is any non-trivial element of $\pi_1(M)$
and that $I_\gamma\vert_{X_0}$ is identically equal to $2$.
Since the principal component $X_0$ contains the character $\chi_0$ of
a discrete, faithful representation $\rho_0$ of $\pi_1(M)$, we
have $\trace\rho_0(\gamma)=I_{\gamma}(\chi_0)=2$. This is possible
only if $\gamma$ is a peripheral element.  But this contradicts the
first part of the statement.
\EndProof

The following result is a strengthened version of Proposition 1.1.2 from
\cite{cgls}.

\Proposition\label{nsurfaces}
Let $M$ be a hyperbolic knot manifold.
Let $X_0$ be a
principal component of $X(\pi_1(M))$ and
let $x_1,\ldots,x_n$ denote the ideal points of $\tilde X_0$.
Then there exists a unique norm $\|\cdot\|$
on the vector space $H_1(\partial M;\mathbb{R})$
such that for any element
$\alpha \in H_1(\partial M;\mathbb{Z})\subset H_1(\partial M;\mathbb{R})$
and any closed curve $c$ representing $\alpha$
we have
$$\|\alpha\|= 2\deg (I_{[c]} \vert_{X_0}) .$$
Moreover there are strict essential surfaces $S_1,\ldots,S_n$ in $M$
(some of which may be closed) such that the following conditions hold.
\Parts
  \Part{(1)} For $i=1,\ldots, n$, the surface $S_i$ is a
  \sigmaxisurface, where $T_{x_i}$ is defined as in
  \ref{sigmadef}, taking $\Gamma=\pi_1(M)$ and $A=X_0$.

  \Part{(2)} If $c$ is any closed curve in $\partial M$ and if
  $\alpha\in H_1(\partial M;\mathbb{Z})\subset H_1(\partial M;\mathbb{R})$
  is the homology class represented by $c$, then $\|\alpha\|=
  \sum_{i=1}^n\Delta_{\partial M}(c,\partial S_i)$. 

  \Part{(3)} For any $k>0$ the set $B_k = \{v \colon \|v\| \le k\}$ in
  $H_1(\partial M;\mathbb{R})$ is a convex polygon.  Furthermore, for each
  vertex $v$ of $B_k$ there is an index $i\le n$ such that $\partial
  S_i \not=\emptyset$ and $v$ is a scalar multiple of the boundary
  class of $S_i$.

  \Part{(4)} If $i\in\{1,\ldots,n\}$ and if $K$ is a $\pi_1$--injective,
  connected, compact subpolyhedron of $S_i$ such that $\chi(K)<0$, and
  if $\Theta\le\pi_1(M)$ is the subgroup defined up to conjugacy by
  $\Theta=\image(\pi_1(K)\to\pi_1(M))$ then, for every element
  $\gamma$ of the commutator subgroup $[\Theta,\Theta]\le
  \Theta\le\pi_1(M)$, we have $Z_{x_i}(I_\gamma-2)\ge t_{S_i}(K)$.
\EndParts
\EndProposition

\Proof
According to Proposition \ref{ActionNontrivial}, $T_i =
T_{x_i}$ is a non-trivial bipartite $\pi_1(M)$--tree for each
$i\in\{1,\ldots,n\}$.  Applying Theorem \ref{dualsurface}, we obtain
an essential $T_{x_i}$--surface $S_i$ satisfying conditions
\ref{dualsurface}(1) and \ref{dualsurface}(2).  (Note that we have
not yet shown that the $S_i$ are strict.)

Let ${\mathcal I}$ denote the unique alternating bilinear form on
$H_1(\partial M; \mathbb{R})$ which restricts to the homological
intersection pairing on $H_1(\partial M;\mathbb{Z})$.  For each $i =
1,\ldots, n$ we define an element $\alpha_i\in H_1(\partial M;\mathbb{Z})\subset H_1(\partial M;\mathbb{R})$ as follows: if $\partial
S_i\not=\emptyset$ we take $\alpha_i$ to be a boundary class for
$S_i$, and if $S_i$ is closed we set $\alpha_i=0$.  We define a linear
functional $l_i$ on $H_1(\partial M; \mathbb{R})$ by $l_i(x) = {\mathcal
I}(x, \alpha_i)$.  Note that if $S_i$ has non-empty boundary then the
kernel of $l_i$ is spanned by $\alpha_i$.

If $c$ is an arbitrary closed curve in $\partial M$ then Proposition
\ref{poleandtranslation} implies for each $i\in\{1,\ldots,n\}$,
that $2\Pi_{x_i}(I_{[c]})$ is equal to $\lambda_{T_i}(c)$, which in
turn is equal to $\Delta_{\partial M}(c,\partial S_i)$ according to
condition \ref{dualsurface}(1). Summing over $i=1,\ldots,n$ we find
that
\Equation\label{degree}
2\deg I_{[c]} = \sum_{i=1}^n 2\Pi_{x_i}(I_{[c]}) =
\sum_{i=1}^n\Delta_{\partial M}(c,\partial S_i) =
\sum_{i=1}^n |l_i(\alpha)| 
\EndEquation
where $\alpha \in H_1(\partial M; \mathbb{R})$ is the class represented
by $c$.  By Lemma \ref{nonconstant}, $\deg I_{[c]}\vert_{X_0}$ is
non-zero for any homotopically non-trivial closed curve $c$ in
$\partial M$.  Thus the $\alpha_i$ span the vector space $H_1(\partial
M; \mathbb{R})$, and we may define a norm on $H_1(\partial M; \mathbb{R})$
by setting
\Equation\label{normdef}
\|v\| = \sum_{i=1}^n |l_i(v)| .
\EndEquation
It follows from \ref{degree} that if $\alpha \in H_1(\partial M;\mathbb{Z})\subset H_1(\partial M;\mathbb{R})$ is represented by a closed curve
$c$ then $\|\alpha\|= 2\deg (I_{[c]} \vert_{X_0})$.  The uniqueness
assertion follows from the observation that, by continuity and
homogeneity, any norm on $H_1(\partial M;\mathbb{R})$ is uniquely
determined by its restriction to the integer lattice $H_1(\partial
M;\mathbb{Z})$.

Conclusion (1) is immediate from the construction of the $S_i$ and
conclusion (2) follows from \ref{degree}.  It follows from
\ref{normdef} that for each vertex $s$ of the convex
polygon $B_k = \{v \;\colon \|v\| \le k\}$, there is an index $i$ such
that the linear functional $l_i$ is not identically 0 and $l_i(s)=0$.
Since $l_i\not\equiv 0$ we have $\partial
S_i\not=\emptyset$.  The kernel of $l_i$ is therefore spanned by the
boundary class $\alpha_i$ of $S_i$.  This implies conclusion (3).

To establish conclusion (4), we suppose that we are given an index
$i\in\{1,\ldots,n\}$ and a $\pi_1$--injective, connected, compact
subpolyhedron $K$ of $S_i$ such that $\chi(K)<0$.  We let
$\Theta\le\pi_1(M)$ denote the subgroup defined up to conjugacy by
$\Theta=\image(\pi_1(K)\to\pi_1(M))$, and suppose that $\gamma$ is an
element of $[\Theta,\Theta]$. If $t_1$ is any positive integer $\le
t_{S_i}(K)$, then according to condition \ref{dualsurface}(2),
$\Theta$ fixes an arc of length $t_1$ in $T_i$.  By
\ref{commutatorandidealpoint} it therefore follows that
$Z_{x_i}(I_\gamma-2)\ge t_1$. As this holds for every positive integer
$t_1\le t_{S_i}(K)$ we conclude that $Z_{x_i}(I_\gamma-2)\ge
t_{S_i}(K)$.

Finally, we must show that each $S_i$
is a strict essential surface. Assume to the contrary that
 $S_i$ is a semi-fiber for some $i$. According to \ref{thicknessdef}
 we then have $t_{S_i}(S_i) = +\infty$. On the other hand, according
 to Proposition \ref{whychiisnegative} we have $\chi(S_i)<0$, and so
$\Theta=\image(\pi_1(S_i)\to\pi_1(M))$ is non-abelian. Let us choose a
non-trivial element $\gamma$ of $[\Theta,\Theta]$. Applying condition
(4) with $K=S_i$ we deduce that $Z_{x_i}(I_\gamma-2)=+\infty$, ie, that
$I_\gamma$ must be the constant function $2$. But since $\gamma$ is
non-trivial, this contradicts Lemma \ref{nonconstant}.
\EndProof

\Remark\label{cyclicremark}
In \cite{cgls} the function $f_\alpha$ in $\mathbb{C}(X_0)$ was defined
by $f_\alpha = I_{[c]}^2 - 4$ where $c$ is a closed curve in $\partial
M$ representing $\alpha$.  The norm referred to in Proposition 1.1.2
of \cite{cgls} satisfies the condition $\|\alpha\| = \deg f_\alpha$
for all $\alpha \in H_1(\partial M; \mathbb{Z})$.  Since the degree of
$f_\alpha$ is twice that of $I_{[c]}$, the norm referred to in
Proposition 1.1.2 of \cite{cgls} is the same as that in the given by
Proposition \ref{nsurfaces}.  We may therefore apply Corollary 1.1.4
of \cite{cgls} to conclude that if $s = \langle\alpha_1\rangle$ is not
the boundary slope of any strict essential surface in $M$, and if the
Dehn filled manifold $M(s)$ has cyclic fundamental group, then then
$\|\alpha\| \le \|\beta\|$ for any non-zero class $\beta \in
H_1(\partial M; \mathbb{Z})$.
\EndRemark

\Corollary\label{twoslopes}
If $M$ is a hyperbolic knot manifold then $M$ has two bounded,
strict, connected essential surfaces with distinct boundary slopes.
\EndCorollary

\Proof
It suffices to show that the surfaces $S_i$ given
by Proposition \ref{nsurfaces} do not all have the same boundary
slope.  If they all did have the same boundary slope, there would be a
non-zero class $\alpha\in H_1(\partial M; \mathbb{Z})$ which is a boundary
class for each $S_i$.  But then the expression given in conclusion (2)
would vanish on the subspace spanned by $\alpha$, contradicting the
fact that this expression defines a norm.
\EndProof

\Proposition\label{poleandintersectionnumber} Suppose that
$M$ is a hyperbolic knot manifold. Let $X_0$ be a principal component
of $X(\pi_1(N))$ and let $x\in \tilde X_0$ be an ideal point.  If $S$ is an
essential $T_x$--surface then for any closed curve $c$ in $M$ we have
$2\Pi_x(I_\gamma)\le\Delta_M(c, S)$.
\EndProposition

\Proof
We have $2\Pi_x(I_{[c]})=\lambda_{T_x}(c)$ by Proposition
\ref{poleandtranslation}, and $\lambda_{T_x}(c)\le\Delta_M(c, S)$ by
Proposition \ref{Deltaandtranslation}. 
\EndProof

\section{Manifolds with few essential surfaces}
\label{fewsurfaces}

The goal of this section is to prove Theorem \ref{AtMostTwoSES},
which gives topological information about an irreducible knot manifold
that has at most two isotopy classes of connected, strict essential
surfaces.  There are a few knot manifolds with this property that
arise as exceptions.  We will discuss these before stating the
theorem.

\Number\label{specialmanifolds}
The solid torus $S^1\times D^2$ and the twisted $I$--bundle $K$ over
the Klein bottle are examples of Seifert fibered knot manifolds which
have no strict essential surfaces at all.  The only connected
essential surface in the solid torus is the meridian disk, which is
obviously a fiber in a fibration over $S^1$ and hence not strict.  The
connected essential surfaces in $K$ are all non-trivial vertical
annuli with respect to the $I$--fibration.  Splitting $K$ along an
essential separating vertical annulus $A$ results in two twisted
$I$--bundles over M\"obius bands for which $A$ is the associated
$\partial I$--bundle, and hence $A$ is a semi-fiber.  Similarly any
non-separating vertical annulus is a fiber in a fibration of $K$ over
$S^1$.

The only 3--manifolds that fiber over the circle with fiber an annulus
are $K$ and $S^1\times S^1\times I$.  Furthermore $K$ is the only
orientable 3--manifold that can be obtained from two twisted
$I$--bundles over M\"obius bands by identifying their $\partial
I$--bundles.  Hence if $M$ is a 3--manifold not homeomorphic to
$K$ or $S^1\times S^1\times I$ then any essential annulus in $M$
is a strict essential surface.

A Seifert fibered manifold $M$ with base surface a disk and two
singular fibers has exactly one isotopy class of connected essential
vertical surfaces, which are annuli.  If the two singular fibers are
both of order $2$ then $M$ is homeomorphic to the twisted $I$--bundle
$K$.  Otherwise a vertical annulus is a strict essential surface;
so $M$ gives an example of an irreducible knot manifold with exactly
one isotopy class of connected strict essential surfaces.
\EndNumber

\Number\label{cabledef}
We define a {\it cable space} to be a Seifert fibered manifold over an
annulus with one singular fiber.  Note that a cable space has exactly
three isotopy classes of essential vertical annuli; one has a
boundary curve on each boundary torus of the cable space and the
other two have both boundary curves on the same boundary torus.
\EndNumber

\Number\label{exceptional}
We will say that an orientable $3$--manifold $M$ is an {\it exceptional
graph manifold} if $M$ is not Seifert fibered and $M$ is
homeomorphic to either
\Parts
  \Part{(1)} a manifold obtained from a disjoint union of a cable
  space $C$ and a twisted $I$--bundle $K$ over a Klein bottle by gluing
  $\partial K$ to a component of $\partial C$ via some homeomorphism; or

  \Part{(2)} a manifold obtained from $P\times S^1$, where $P$ is a
  planar surface with three boundary curves, by gluing two of the
  boundary tori of $P\times S^1$ to each other via some homeomorphism.
\EndParts
\EndNumber

\Proposition\label{exceptionalprop}
 An exceptional graph manifold has exactly two
connected strict essential surfaces up to isotopy.  One of these is a
torus and the other is an annulus. 
\EndProposition

\Proof
Let $M$ denote an exceptional graph manifold.  According to the
definition, $M$ is obtained from a manifold $M'$ by identifying two
torus boundary components of $M'$.  The image of these two tori under
the quotient map $q\co M'\to M$ is a torus $T$ in $M$.  We denote by
$M_0$ the component of $M'$ which contains $q^{-1}(\partial M)$. In
case (1) of the definition $M_0$ is a cable space and the other
component of $M'$, which we shall denote by $M_1$, is a twisted
$I$--bundle over a Klein bottle.  We shall regard $M_1$ as a Seifert
fibration over a disk with two singular fibers of order $2$.  In case
(2) of the definition $M' = M_0$ is homeomorphic to $P\times S^1$,
where $P$ is a planar surface with three boundary components.  In
either case the manifold $M'$ is a Seifert-fibered manifold and, up to
isotopy, there is a unique essential annulus in $M_0$ which has both
boundary components in $q^{-1}(\partial M)$.  We will let $A'$ denote
such an annulus.

Clearly $T$ is an essential surface in $M$, and $T$ is strict by
\ref{meetsall}.  The annulus $A'$ is a strict essential surface
in $M_0$ by \ref{meetsall}, and hence by Proposition \ref{mfldwithtorus}
the annulus $A = q(A')$ is a strict essential surface in $M$.
We will show that any connected strict essential surface in $M$ is
isotopic either to $A$ or to $T$.

Suppose that $F$ is a connected strict essential surface in $M$ which
is not isotopic to $T$.  By Proposition \ref{mfldwithtorus} we may assume
after an isotopy that $F$ is transverse to $T$, that each component of
$F'=q^{-1}(F)$ is essential in the component of $M'$ containing it,
and that some component $S$ of $F'$ is a strict essential surface in
the component of $M'$ containing it.  Since a twisted $I$--bundle over
a Klein bottle has no strict essential surfaces we must have $S\subset
M_0$.  Thus $S$ is a component of $F_0 = F'\cap M_0$.

We claim that $F'$ is isotopic to a union of vertical annuli in the
Seifert-fibered manifold $M'$.  Since the strict essential surface $S$
is a component of $F_0$, it follows from \ref{strictlyvertical} that
$F_0$ cannot be isotopic to a horizontal surface, and hence that it
is isotopic to a vertical surface whose components are all vertical
annuli in $M_0$.  This proves the claim if $M$ satisfies case (2) of
the definition of an exceptional graph manifold.  To complete the
proof in case (1) it is enough to show that $F_1$ is also isotopic to
a vertical surface.  This is true because, in a Seifert fibration over
a disk with two singular fibers of order $2$, every essential surface
is isotopic to a vertical surface.  Thus the claim is proved in both cases.

Next we claim that $\partial F_0 \subset q^{-1}(\partial M)$.  Assume
to the contrary that $\partial F'\cap q^{-1}(T) \not=\emptyset$. Let
$T_0$ and $T_1$ denote the two components of $q^{-1}(T)$, where
$T_0\subset M_0$, and let $h\co T_0 \to T_1$ denote the gluing
homeomorphism.  Since $F'$ is isotopic to a vertical surface and $h$
maps $\partial F'\cap T_0\not=\emptyset$ to $\partial F'\cap T_1$, it
follows that $h$ is isotopic to a fiber-preserving homeomorphism.
Hence $M$ admits a Seifert fibration; this contradicts the definition
of an exceptional graph manifold.

We have now shown that $F_0$ is a vertical surface in $M_0$ and
that $\partial F_0\subset q^{-1}(\partial M)$.  Moreover, $F'$ is
connected by hypothesis.  Thus $F'$ is isotopic to the annulus $A'$,
and $F$ is isotopic to the annulus $A$, as required for the proof of
the proposition.
\EndProof

\Proposition\label{hasSES}
Let $M$ be a compact, irreducible orientable 3--manifold whose boundary
components are all tori, and let $T_0$ be a component of $\partial M$.
If $M$ is not homeomorphic to $S^1\times D^2$, $S^1\times S^1\times I$
or a twisted $I$--bundle $K$ over the Klein bottle, then $M$ contains a
bounded connected strict essential surface $S$ such that $\partial S
\subset T_0$.  Moreover, if $M$ is the compact core of a complete
hyperbolic manifold with finite volume then there are two bounded
connected strict essential surfaces in $M$ which have their boundaries
contained in $T_0$ and which have distinct boundary slopes on $T_0$.
\EndProposition

\Proof
First consider the case where $M$ is the compact core of a complete
hyperbolic manifold with finite volume.  If $\partial M$ is connected,
then the result follows from Corollary \ref{twoslopes}.  If $M$ has
more than one boundary component then \cite[Theorem 3]{boundsep}
implies that $M$ has an essential surface $F$ with $\partial F$
contained in $T_0$; in this case $F$ is disjoint from at least one
component of $\partial M$ and is therefore strict by \ref{meetsall}.

If $M$ is Seifert-fibered and is not homeomorphic to one of the
exceptional manifolds listed in the statement, then we will show that
$M$ has a strict essential vertical annulus $A$ whose boundary is
contained in $T_0$.  This implies the result in this case.  First note
that there is an arc $\alpha$ in the base surface $B$ such that
$\alpha$ is essential (ie, is not the frontier of a disk disjoint
from the image of the singular fibers), and such that $\alpha$ has
both its endpoints on the component of $\partial B$ which is the image
of $T_0$.  Indeed if such an arc $\alpha$ did not exist are where $B$
would be a disk or an annulus, and there would be no singular fibers
in the Seifert fibration of $M$; this would imply that $M$ is
homeomorphic to $S^1\times D^2$ or $S^1\times S^1\times I$, a
contradiction.

The inverse image of $\alpha$ under the Seifert fibration is then an
essential annulus $A$ which has both boundary components contained in
$T_0$.  Since $M$ is not homeomorphic to $K$ or $S^1\times S^1\times
I$, it follows from \ref{specialmanifolds} that the essential
annulus $A$ is a strict essential surface.

To prove the proposition in the general case, let ${\mathcal T}$ be a
maximal collection of disjoint essential tori in $M$, no two of which
are parallel.  We may assume that ${\mathcal T}$ is non-empty since
otherwise, by Thurston's Geometrization Theorem, $M$ would either be
Seifert-fibered or homeomorphic to the compact core of a complete
hyperbolic manifold finite volume.  Let $R$ be a regular neighborhood
of ${\mathcal T}$ and let $N$ be the closure of the component of $M - R$
which contains $T_0$.  It suffices to show that $N$ contains a bounded
strict essential surface which is disjoint from $\partial N - T_0$.
Note that $N$ is not homeomorphic to $S^1\times S^1\times I$, since
the tori in ${\mathcal T}$ are essential.  Also, $N$ cannot be
homeomorphic to $K$ or $S^1\times D^2$ since ${\mathcal T}$ is non-empty.
Since $N$ contains no essential tori it follows from Thurston's
theorem that either $N$ is Seifert-fibered or it is homeomorphic to
the compact core of a complete hyperbolic manifold with finite volume.
Thus $N$ contains a bounded strict essential surface which is disjoint
from $\partial N - T_0$ by the two cases that were handled earlier.
This completes the proof.
\EndProof

\Proposition\label{torusgivesSES}
If a knot manifold $M$ has an essential torus then it also has a
bounded strict essential surface.  Furthermore, a Seifert-fibered knot
manifold which contains an essential torus has infinitely many
distinct isotopy classes of strict essential surfaces.
\EndProposition

\Proof
Consider the manifold $N$ obtained by splitting $M$ along a maximal
family ${\mathcal T}$ of disjoint, non-parallel essential tori.  Let $T_0$
denote the component of $\partial N$ which corresponds to $\partial M$
and let $N_0$ be the component of $N$ containing $T_0$.  Since the
tori in ${\mathcal T}$ are essential, $N_0$ is not homeomorphic to
$S^1\times S^1\times I$.  Thus, since $N_0$ has at least two boundary
components, it is not one of the exceptional manifolds listed in
Proposition \ref{hasSES}.  Hence $N_0$ has a bounded strict essential
surface $F$ with $\partial F \subset T_0$.  Now $F$ is a bounded
essential surface in $M$ which is strict since it is disjoint from an
essential torus.

For the proof of the second assertion, assume that $M$ is
Seifert-fibered and consider an essential vertical torus $T$.  The
image of $T$ under the Seifert fibration map is a simple closed curve
$c$ in the base surface $B$.  Since $T$ is essential, $c$ does not
bound a disk containing fewer than two points which are images of
singular fibers.  Therefore there exists an arc $\alpha$ in $B$ which
meets $c$ transversely in at least one point and has the property that
every disk component of the complement of $c\cup \alpha$ contains the
image of at least one singular fiber.  The images of $\alpha$ under
powers of the Dehn twist about $c$ give an infinite family of
non-isotopic arcs in $B$ whose inverse images in $M$ are strict
essential annuli.
\EndProof

\Theorem\label{AtMostTwoSES}
Let $M$ be an irreducible knot manifold.
\Parts
  \Part{(1)} If $M$ has no strict essential surface then $M$ is
  homeomorphic to either a solid torus or a twisted $I$--bundle over a
  Klein bottle.

  \Part{(2)} If $M$ has exactly one isotopy class of connected strict
  essential surfaces then $M$ is Seifert-fibered over a disk with two
  singular fibers.

  \Part{(3)} If $M$ has exactly two isotopy classes of connected
  strict essential surfaces, represented by surfaces $F_1$ and $F_2$,
  then either

    \Subpart{(3a)} $M$ is an exceptional graph manifold; or

    \Subpart{(3b)} $M$ is a hyperbolic knot manifold, $F_1$ and $F_2$
    are bounded surfaces of negative Euler characteristic, and the
    boundary slopes of $F_1$ and $F_2$ are distinct.
\EndParts
\EndTheorem

The proof of Theorem \ref{AtMostTwoSES} depends on the following
lemma, which is contained in \cite[Lemma 2.3]
{klaffshalen}.

\Lemma\label{cablelemma}
Suppose that $N$ is a cable space (cf \ref{cabledef}) with boundary
tori $T_1$ and $T_2$.  Then there is a bijection $\phi$ from the set
of slopes on $T_1$ to the set of slopes on $T_2$ such that for each
slope $s$ on $T_1$ there exists a connected essential surface in $N$,
having nonempty intersection with both $T_1$ and $T_2$ and having $s$
and $\phi(s)$ as boundary slopes.
\EndLemma

\Proof[Proof of Theorem \ref{AtMostTwoSES}]
If $M$ has no strict essential surface then Proposition \ref{hasSES}
implies that $M$ is homeomorphic to either a solid torus, a twisted
$I$--bundle over a Klein bottle or a product of a torus and an
interval.  Since the latter is not a knot manifold, assertion (1) of
the theorem follows.

As a preliminary to proving assertions (2) and (3) we observe that if
$M$ has no essential torus then, by Thurston's Geometrization Theorem,
$M$ is either a hyperbolic knot manifold or a Seifert-fibered
manifold.  Moreover, all Seifert-fibered knot manifolds with no
essential tori are Seifert-fibered over a disk with two singular
fibers.

To prove assertion (2), suppose that $M$ contains exactly one strict
essential surface $S$.  Proposition \ref{torusgivesSES} implies that
$S$ cannot be a torus, and Proposition \ref{hasSES} implies that $M$
cannot be a hyperbolic knot manifold.  Thus the observation above
implies the conclusion of (2).

We now turn to the proof of assertion (3).  Assume that, up to
isotopy, $M$ contains exactly two strict essential surfaces $F_1$ and
$F_2$.  According to \ref{specialmanifolds} there is only one isotopy
class of strict essential surfaces in a Seifert-fibered manifold over
a disk with at most two singular fibers. Thus by the observation above
we have two cases: either $M$ is hyperbolic, or $M$ has an essential
torus.

We first consider the case that $M$ is hyperbolic. Proposition
\ref{hasSES} implies that $F_1$ and $F_2$ are bounded and have
distinct boundary slopes. Proposition \ref{whychiisnegative} implies
that $\chi(F_1)$ and $\chi(F_2)$ are strictly negative. Thus
conclusion (3b) holds in this case.

The remaining case is that $M$ contains an essential torus $T$.  It
follows from \ref{meetsall} that any closed essential surface in $M$
is a strict essential surface.  In particular $T$ is strict.
Proposition \ref{torusgivesSES} implies that there is also a bounded
strict essential surface $A$ in $M$, and by the hypothesis of
assertion (3) any strict essential surface in $M$ is isotopic to
either $T$ or $A$.  In particular, any closed essential surface in $M$
is isotopic to $T$.  Note that Proposition \ref{torusgivesSES} also
implies that $M$ is not Seifert-fibered.

Let $M'$ be the manifold obtained by splitting $M$ along the torus
$T$, and let $M_0$ denote the component of $M'$ which contains
$q^{-1}(\partial M)$ where $q\co M'\to M$ is the quotient map.  Note that
$M'$ cannot contain a closed essential surface since the image of
such a surface in $M$ would be a closed essential surface but would
not be isotopic to $T$.  In particular, no component of $M'$
contains an essential torus.

We claim that $M_0$ cannot be homeomorphic to the compact core of a
complete hyperbolic manifold with finite volume.  Otherwise by
Proposition \ref{hasSES}, $M_0$ would contain two non-isotopic
bounded strict essential surfaces which are disjoint from $T$.  These
would be strict essential surfaces in $M$ by \ref{meetsall}.  Since
$A$ is the only bounded strict essential surface in $M$ up to isotopy
this is a contradiction. 

It now follows from Thurston's Geometrization Theorem that $M_0$ is a
Seifert fibered manifold.  Moreover, since $T$ is not
boundary-parallel, $M_0$ is not homeomorphic to $S^1\times S^1\times
I$ and must therefore contain an essential annulus with its boundary
contained in $q^{-1}(\partial M)$.  By \ref{meetsall} any such
annulus must be a strict essential surface in $M$.  It follows that
$A$ is an annulus and that, up to isotopy, $A$ is the only essential
annulus in the Seifert fibered manifold $M_0$.

If $T$ is non-separating, then $M_0$ is a Seifert-fibered manifold
with three boundary components and contains only one essential annulus up to
isotopy.  It follows that $M_0$ is homeomorphic to $P\times S^1$ where
$P$ is a planar surface with three boundary curves.  Since $M$ is
not Seifert-fibered, $M$ is an exceptional graph manifold by
\ref{exceptional}(2).  Thus (3a) holds.

Now consider the case where $T$ is separating.  We let $M_1$ denote the
component of $M'$ which does not contain $q^{-1}(\partial M)$, and we
identify $M_0$ and $M_1$ with submanifolds of $M$.  In this case $M_0$
has two boundary components and is a Seifert-fibered manifold which
contains only one essential annulus up to isotopy.  Hence $M_0$ is a
cable space.  It follows from Proposition \ref{hasSES} that either
$M_1$ is a twisted $I$--bundle over a Klein bottle, or $M_1$ has a
bounded strict essential surface.  If $M_1$ is a twisted $I$--bundle
over a Klein bottle then, since $M$ is not Seifert-fibered, $M$ is an
exceptional graph manifold by \ref{exceptional}(2), so (3a) holds.

Finally suppose that $M_1$ contains a bounded strict essential surface
$F$, and let $s$ denote the boundary slope of $F$ in $M_1$.  By Lemma
\ref{cablelemma} there exists a connected essential surface $G$ in
the cable space $M_0$, having boundary slope $s$ on $T$ and boundary
slope $\phi(s)$ on $\partial M_0$. Thus, for suitably chosen positive
integers $m$ and $n$, the surface consisting of $m$ parallel copies of
$F$ in $M_1$ and the surface consisting of $n$ parallel copies of $G$
in $M_0$ have isotopic intersections with $T$.  Hence there exists a
connected surface $\widehat F$ in $M$ which meets $M_1$ in parallel
surfaces isotopic to $F$ and meets $M_0$ in parallel surfaces isotopic
to $G$.  According to Proposition \ref{mfldwithtorus}, $\widehat F$ is a
bounded strict essential surface in $M$.  Since $\widehat F$ is not
isotopic to the annulus $A$, this is a contradiction.  Thus $M_1$
cannot contain a bounded strict essential surface, and the proof is
complete. 
\EndProof

\Definition\label{TwoSurfaceKnotMfld}
We will say that $M$ is a {\it two-surface knot manifold} provided
that $M$ is an irreducible knot manifold and that $M$ has at most two
distinct isotopy classes of strict essential surfaces.  We say that a
two-surface knot manifold $M$ is an {\it exceptional two-surface knot
manifold} if $M$ is Seifert fibered or if $M$ is an exceptional graph
manifold.
\EndDefinition

\Number\label{twosurfprops}
According to Theorem \ref{AtMostTwoSES}, if $M$ is a non-exceptional
two-surface knot manifold then $M$ is a hyperbolic knot manifold
and $M$ has exactly two distinct isotopy classes of connected
strict essential surfaces.  Moreover, if $F_1$ and $F_2$ are
representatives of these two isotopy classes then they have distinct
boundary slopes and both $F_1$ and $F_2$ have negative Euler
characteristic.  These properties of non-exceptional two-surface knot
manifolds will be used throughout the sequel.
\EndNumber

\Definition
A knot $K$ in a closed orientable $3$--manifold $\Sigma$ will be said
to be a {\it non-exceptional two-surface knot} provided that the knot
manifold $M = \Sigma(K)$ is a non-exceptional two-surface knot
manifold.
\EndDefinition

\section{General principles about two-surface knot manifolds}
\label{generalsection}

As we explained in the introduction, there is a combinatorially defined
quantity $\kappa(F_1,F_2)$, associated to two bounded essential
surfaces, which plays a central role in our estimates.  We begin
with the definition.

\Definition\label{kappadef}
Suppose that $F_1$ and $F_2$ are bounded connected essential surfaces
in an irreducible knot manifold $M$.  Let $m_i$ denote the number of
boundary components of $F_i$ for $i=1,2$. We define an element
$\kappa(F_1,F_2)$ of the interval $[0,\infty]$ of the extended real
line by
  $$\kappa(F_1,F_2)=\inf_K \frac{m_2\cdot\#(K\cap F_1)}{m_1\cdot
t_{F_2'}(K)},$$
where $F_2'$ ranges over all surfaces that are isotopic to $F_2$ and
meet $F_1$ transversally, while $K$ ranges over all compact connected
$\pi_1$--injective $1$--dimensional polyhedra of Betti number $2$ which
are contained in $ F_2'$ and meet $F_1$ transversally. (To say that
$K$ meets $F_1$ transversally means in particular that $K\cap F_1$
consists entirely of points at which $K$ is locally Euclidean.
We interpret the quotient in the definition as being $0$ if
$t_{F_2'}=\infty$.)
Note
that $\kappa(F_1,F_2)<\infty$ if and only if $F_2$ contains a
$\pi_1$--injective connected $1$--dimensional polyhedron of Betti number
$2$, ie, if and only if $\chi(F_2)<0$.
\EndDefinition

\Theorem\label{generalinequality}
Suppose that $M$ is a non-exceptional two-surface knot manifold (so
that $M$ is hyperbolic by \ref{twosurfprops}).  Let $X_0$ be a
principal component of $X(\pi_1(M))$ and let $\|\cdot\|$ denote the
norm on $H_1(\partial M,\mathbb{R})$ given by Proposition
\ref{nsurfaces}.  Let $F_1$ and $F_2$ be representatives of the two
isotopy classes of connected strict essential surfaces in $M$, and
for $i=1,2$ let $\alpha_i$ denote a boundary class of $F_i$
(which is a bounded surface by \ref{twosurfprops}).
Then we have
  $$\frac{\|\alpha_1\|}{\|\alpha_2\|}\le{\four}\kappa(F_1,F_2).$$
\EndTheorem

\Proof
Suppose that $F_1$ and $F_2$ satisfy the hypotheses and also meet
transversally.  Let $K$ be a compact connected $\pi_1$--injective
$1$--dimensional polyhedron of Betti number $2$ contained in $ F_2$ and
meeting $F_1$ transversally.  Set $t= t_{F_2}(K)$ and
$\ell=\#(K\cap F_1)$.  In this setting we will show that
$$\frac{\|\alpha_1\|}{\|\alpha_2\|}\le{\four}\frac{m_2\ell}{m_1t}.$$
In view of the definition of $\kappa(F_1, F_2)$, this will
establish the theorem.

Let $x_1,\ldots,x_n$ be the ideal points of $\tilde X_0$.
We fix strict essential surfaces $S_1,\ldots,S_n$ in $M$
satisfying conditions (1)--(4) of Proposition \ref{nsurfaces}.  It
follows from the hypotheses that each component of each $S_i$ is
isotopic to either $F_1$ or $F_2$.
For $i = 1,2$ let $m_i$ denote the number of boundary components of
$F_i$ and let $s_i = \langle\alpha_i\rangle$ denote its boundary
slope.  Since $s_1 \not= s_2$ by \ref{twosurfprops}, there cannot
exist disjoint surfaces isotopic to $F_1$ and $F_2$.  Hence for each
$i\in\{1,\ldots,n\}$, either every component of $S_i$ is isotopic to
$F_1$ or every component of $S_i$ is isotopic to $F_2$. We may
therefore suppose the $S_i$ (and the $x_i$) to be indexed in such a
way that all components of $S_i$ isotopic to $F_1$ when $1\le i\le k$,
and all components of $S_i$ are isotopic to $F_2$ when $k<i\le n$.
Here $k$ is {\it a priori} an integer with $0\le k\le n$. However, if
$k$ were equal to $0$ or $n$ then the $S_i$ would all have the same
boundary slope.  It would then follow from part
(2) of Proposition \ref{nsurfaces} that $\|\alpha\| = 0$,
contradicting the definition of a norm.  Hence $0<k<n$.

We let $\nu_i$ denote the number of components of $S_i$ for
$i=1,\ldots,n$, and we set $N_1=\sum_{i=1}^k\nu_i$ and
$N_2=\sum_{i=k+1}^n\nu_i$. For $1\le i\le k$, the boundary of $S_i$
consists of $\nu_im_1$ simple closed curves of slope $s_1$, and for
$k< i\le n$, the boundary of $S_i$ consists of $\nu_im_2$ simple
closed curves of slope $s_2$.  Thus if $C\subset \partial M$ is a
non-trivial simple closed curve and $s$ denotes its slope, we have
$$\Delta_{\partial M}(C,\partial S_i)=\nu_im_1\Delta(s,s_1)\qquad{\rm
for}\ 1\le i\le k$$
$$\Delta_{\partial M}(C,\partial S_i)=\nu_im_2\Delta(s,s_2)\qquad{\rm
for}\ k< i\le n.\leqno{\hbox{and}}$$
Hence if $\beta$ is the homology class in $H_1(\partial M;\mathbb{Z})\subset H_1(\partial M;\mathbb{R})$ represented by some orientation of
$C$ then \ref{nsurfaces}(2) gives
$$
\begin{aligned}
\|\beta\|=\sum_{i=1}^n\Delta_{\partial M}(C,\partial
S_i)&=\sum_{i=1}^k\nu_im_1\Delta(s,s_1)+\sum_{i=k+1}^n\nu_im_2\Delta(s,s_2)\\
&=N_1m_1\Delta(s,s_1)+N_2m_2\Delta(s,s_2) .
\end{aligned}
$$
In particular, taking $C$ to be a simple closed curve with slope $s_1$
or $s_2$, setting $\Delta=\Delta(s_1,s_2)=\Delta(s_2,s_1)$ and
observing that $\Delta(s_1,s_1)=\Delta(s_2,s_2)=0$, we find that
$$\|\alpha_1\|=N_2m_2\Delta\qquad{\rm and}
\qquad\|\alpha_2\|=N_1m_1\Delta,$$
so that
\Equation\label{generalone}
\frac{\|\alpha_1\|}{\|\alpha_2\|}=\frac{N_2m_2}{N_1m_1} .
\EndEquation
Let us fix collarings $h_1$ and $h_2$
of $F_1$ and $F_2$ in $M$. Since each component of $S_i$ is isotopic
to $F_1$ if $i\le k$ and to $F_2$ if $i>k$, it follows from
Proposition \ref{linemup} that after modifying each $S_i$ within its isotopy
class we may assume that $S_i$ has the form $h_1(F_1\times Y_i)$ (if
$i\le k$) or $h_2(F_2\times Y_i)$ (if $i> k$) where $Y_i\subset
[-1,1]$ is a set of cardinality $\nu_i$. Note also that since, by
\ref{nsurfaces}(1), $S_i$ is a $T_{x_i}$--surface, and since the $T_{x_i}$
are non-trivial $\pi_1(M)$--trees by Proposition \ref{ActionNontrivial},
we have $S_i\ne\emptyset$ for each $i\in\{1,\ldots,n\}$,
and hence $\nu_i>0$ for each $i$. We may therefore assume the isotopic
modifications of the $S_i$ to have been made in such a way that
$0 \in Y_i$ for each $i$. Hence $S_i\supset F_1$ for $i\le k$,
and $S_i\supset F_2$ for $i> k$.

There exist generators $x$ and $y$ of the rank--$2$ free group
$\pi_1(K)$ such that the conjugacy class of the commutator $[x,y]$ is
represented by a map $c\co S^1\to K$ which has the property that
$\#(c^{-1}(p))\le 2$ for all points $p$ at which $K$ is locally
Euclidean.  We regard $c$ as a map of $S^1$ to $M$, ie, a closed
curve in $M$, and denote by $\gamma$ an element of the conjugacy
class in $\pi_1(M)$ which is represented by $c$.  Since $K$ is
$\pi_1$--injective in the essential surface $F\subset M$ we have
$\gamma\ne1$.
 
We consider the function $I_{\gamma}\co X_0\to\mathbb{C}$. Since $K\subset
F_2\subset S_i$ for each $i>k$, and since $\gamma$ is a commutator in
$\image(\pi_1(K) \to \pi_1(M))$, it follows from \ref{nsurfaces}(4)
that $Z_{x_i}(I_{\gamma}-2)\ge t_{S_i}(K)$ for $i=k+1,\ldots,n$. But
since for $i>k$ we have $S_i=h_2(F_2\times Y_i)$, where $0\in
Y_i\subset [-1,1]$ and $\#(Y_i)=\nu_i$, it follows from
\ref{thicknessmultiplies} that $ t_{S_i}(K)\ge\nu_i\cdot
t_{F_2}(K)=\nu_i\cdot t$.  Hence
\Equation\label{generaltwo}
Z_{x_i}(I_{\gamma}-2)\ge\nu_i\cdot t
\EndEquation
for $i=k+1,\ldots,n$.

The function $I_{\gamma}$ is non-constant by Lemma \ref{nonconstant}.
We can therefore estimate its degree by
using (\ref{generaltwo}): we have
$$\deg I_{\gamma}=\deg (I_{\gamma}-2)\ge\sum_{i=k+1}^n
  Z_{x_i}(I_{\gamma}-2)\ge \sum_{i=k+1}^n \nu_i\cdot t,$$
and hence
\Equation\label{generalthree}
\deg I_{\gamma}\ge N_2t .
\EndEquation
We shall compare the lower bound (\ref{generalthree}) for
$\deg I_{\gamma}$ with an upper bound calculated in terms of poles.
The definition of $I_{\gamma}$ shows that it has no poles on the
affine curve $X_0$. For an ideal point $x_i$ with $i>k$, it follows
from (\ref{generaltwo}) that $I_{\gamma}$ takes the value
$2$ at $x_i$ and hence does not have a pole.  Hence
\Equation\label{generalfour}
\deg I_{\gamma}=\sum_{i=1}^k\Pi_{x_i}(I_{\gamma}) .
\EndEquation
For  $i\le k$, it follows from \ref{nsurfaces}(1) and
\ref{poleandintersectionnumber} that
\Equation\label{generalfive}
2\Pi_{x_i}(I_{\gamma})\le\Delta_M(c,S_i) .
\EndEquation
We shall give an upper bound for $\Delta_M(c, S_i)$. Since
$S_i=h_1(F_1\times Y_i)$, where $\#(Y_i)=\nu_i$, and since the
polyhedron $K$ meets $F_1$ transversally, it is clear that $K$ can be
isotoped in $M$ to a polyhedron $K_i$ for which $\#(K_i\cap
S_i)=\nu_i\cdot\#(K_i\cap F_1)=\nu_i\ell$. If we write $K_i=\eta_i
(K)$, where $\eta_i \co M\to M$ is isotopic to the identity, then
$c\co S^1\to M$ is homotopic to $c_i=\eta_i \circ c$.  Since
$\#(c^{-1}(p))\le 2$ for every non-vertex point $p\in K$, we have
\Equation\label{generalsix}
\#(c_i^{-1}(S_i))\le2\nu_i\ell .
\EndEquation
But the definition of geometric intersection number implies that
\Equation\label{generalseven}
\Delta_M(c, S_i)\le \#(c_i^{-1}(S_i)) .
\EndEquation
By combining the inequalities
(\ref{generalfive})--(\ref{generalseven}), summing,
and comparing with (\ref{generalfour}), we find that
$$2\deg I_{\gamma}\le\sum_{i=1}^k2\nu_i\ell$$
ie,
\Equation\label{generaleight}
\deg I_{\gamma}\le N_1\ell .
\EndEquation
 From (\ref{generalthree}) and (\ref{generaleight}) it
follows that
\Equation\label{generalnine} 
N_2t\le N_1\ell .
\EndEquation
Combining (\ref{generalnine}) with (\ref{generalone})
we obtain
$$\frac{\|\alpha_1\|}{\|\alpha_2\|}\le{\four}\frac{m_2\ell}{m_1t}$$
as required.
\EndProof

\Remark\label{weakerhyp}\rm In the proof of Theorem
\ref{generalinequality} we fixed strict essential surfaces
$S_1,\ldots,S_n$ in $M$ satisfying conditions (1)--(4) of Proposition
\ref{nsurfaces} and then concluded from the hypotheses that each
component of each $S_i$ is isotopic to either $F_1$ or $F_2$.  The
proof remains valid as long as $M$ contains two surfaces $F_1$ and
$F_2$ such that each component of each $S_i$ is isotopic to either
$F_1$ or $F_2$.  The results in the rest of this section and the
subsequent sections, would also remain valid under this much weaker,
but much more technical hypothesis.

Nathan Dunfield's computations of A--polynomials suggest that there
are many examples of knot manifolds in lens spaces that have such a pair
of surfaces.

(See \url{http://www.its.caltech.edu/~dunfield/snappea/tables/A-polys}.)
\EndRemark

\Theorem\label{disksandannuli}
Suppose that $M$ is a non-exceptional two-surface knot manifold.  Let
$F_1$ and $F_2$ be representatives of the two isotopy classes of
connected strict essential surfaces.  Suppose that $F_1$ and $F_2$
intersect transversally and that no component of ${\mathcal A} = F_1\cap F_2$ is a
homotopically trivial simple closed curve.  Then for $i=1,2$, every
component of $\inter (F_i-{\mathcal A})$ is an open disk or an open annulus.
\EndTheorem

\Remark\label{standardremark}\rm
  It is a standard observation that if two transversally intersecting
  essential surfaces $F_1$ and $F_2$ are chosen within their
  rel-boundary isotopy classes so as to minimize the number of
  components of $F_1\cap F_2$, then no component of $F_1\cap F_2$ is a
  homotopically trivial simple closed curve.
\EndRemark

\Proof[Proof of Theorem \ref{disksandannuli}]
By symmetry it suffices to prove that every component of
$\inter(F_2-{\mathcal A})$ is a disk or annulus. The hypothesis that
no component of $F_1\cap F_2$ is a homotopically trivial simple closed
curve implies that every component of $\inter(F_2-{\mathcal A})$ is
$\pi_1$--injective in $F_2$, and hence in $M$ since $F_2$ is essential.

By Theorem \ref{generalinequality}, the surfaces $F_1$ and
$F_2$ have non-empty boundaries and, if we let $\alpha_1$ and
$\alpha_2$ denote boundary classes of $F_1$ and $F_2$, we have
$\frac{\|\alpha_1\|}{\|\alpha_2\|}\le{\four}\kappa(F_1,F_2).$ Suppose
that some component $C$ of $\inter(F_2-{\mathcal A})$ is not an open disk or
annulus. Then $\chi(C)<0$. Hence $C$ contains a connected
$1$--dimensional polyhedron $K$ of Betti number $2$ which is
$\pi_1$--injective in $C$ and hence in $M$. Since $K\subset\inter
(F_2-{\mathcal A})$, we have $\#(K\cap F_1)=0$. By the definition of
$\kappa(F_1,F_2)$ it follows that $\kappa(F_1,F_2)=0$, and hence that
$\frac{\|\alpha_1\|}{\|\alpha_2\|}\le0$.  This is impossible, because
$\|\alpha_1\|$ and $\|\alpha_2\|$ are norms of non-zero elements of
$H_1(M;{\partial R})$ and are therefore strictly positive real
numbers.
\EndProof

\Corollary\label{chibound}
Suppose that $M$ is a non-exceptional two-surface knot manifold.  Let
$F_1$ and $F_2$ be representatives of the two isotopy classes of
connected strict essential surfaces.  Let $s_i$ denote the boundary
slope of $F_i$ and let $M_i$ denote the number of boundary components
of $F_i$.  Then for $i=1,2$ we have
$$\chi(F_i) \ge \frac{-m_1m_2\Delta(s_1,s_2)}{2} . $$
\EndCorollary

\proof
Since the number of arc components of $F_1\cap F_2$ is
$m_1m_2\Delta(s_1,s_2)/2$, and since each component of $F_i -
F_1\cap F_2$ has non-negative Euler characteristic by Theorem
\ref{disksandannuli}, we have
$$\chi(F_i) = \chi(F_i - F_1\cap F_2) - \frac{m_1m_2\Delta(s_1,s_2)}{2} \ge
\frac{-m_1m_2\Delta(s_1,s_2)}{2}.\eqno{\qed}$$

\Theorem\label{generalknotinequality} Suppose that $K$ is a
non-exceptional two-surface knot in a closed, orientable $3$--manifold
$\Sigma$ with $\pi_1(\Sigma)$ cyclic.  Set $M = M(K)$ and let $F_1$
and $F_2$ be representatives of the two isotopy classes of connected
strict essential surfaces in the non-exceptional two-surface knot
manifold $M$.  Let $\emm$ denote the meridian slope of $M$, let $s_i$
denote the boundary slope of $F_i$ (which is well-defined by
\ref{twosurfprops}) and assume that $s_2\ne\emm$.  Set
$q_i=\Delta(s_i,\emm)$ (so that $q_i$ is the denominator of $s_i$ in
the sense of \ref{denominator}), and set $\Delta=\Delta(s_1,s_2)$ (so
that $\Delta \not=0$ by \ref{twosurfprops}). Then
$$\frac{q_1^2}{\Delta}\le{\four}2\kappa(F_1,F_2).$$
\EndTheorem

\Proof
The inequality $\frac{q_1^2}{\Delta}\le{\four}2\kappa(F_1,F_2)$ holds
trivially if $s_1=\emm$, since the left hand side is $0$ in this case.
Thus we may assume that $s_1\ne\emm$, and by hypothesis we have
$s_2\ne\emm$. Hence $\emm$ is not the boundary slope of any strict
essential surface.  We choose a meridian class $\mu$ (in the sense of
\ref{denominator}).  Thus $\langle\mu\rangle=\emm$, and $\mu$ is not
a boundary class of any strict essential surface.

If $\|\cdot\|$ denotes the norm on $H_1(\partial M,\mathbb{R})$ associated to a principal component $X_0$ of
$X(\pi_1(M))$, we have
$$\frac{\|\alpha_1\|}{\|\alpha_2\|}\le{\four}\kappa(F_1,F_2),$$
where $\alpha_i$ is a boundary class of $F_i$, so
$s_i = \langle\alpha_i\rangle$.

Let $L$ denote the lattice $H_1(\partial M;\mathbb{Z})$ in the vector
space $V=H_1(\partial M;\mathbb{R})$.  The homological intersection
pairing $L\times L\to\mathbb{Z}$ has a unique extension to an alternating
bilinear form $\omega\co V\times V\to\mathbb{R}$.  Thus for all
$\alpha,\beta\in L$ we have
$\Delta(\langle\alpha\rangle,\langle\beta\rangle)=|\omega(\alpha,\beta)|$.

The alternating form $\omega$ determines an area element on $V$.  If
$v$ and $v'$ are linear independent vectors in $V$ then the
parallelogram with vertices $\{ 0, v, v', v+v'\}$ has area
$|\omega(v,v')|$; in particular, a fundamental parallelogram for $L$
has area $1$.  Furthermore, the parallelogram with vertices $\{ v, v',
-v, -v'\}$ has area $2|\omega(v,v')|$.  (If $e_1$ and $e_2$ form a
basis of $L$ then $|\omega(e_1,e_2)| = 1$ and
$|\omega(xe_1+ye_2,ze_1+we_2)| = |zw-yz|$, so if we use the basis
$(e_1, e_2)$ to identify $L$ with $\mathbb{Z}^2$ and $V$ with $\mathbb{R}^2$
then we recover the standard area element on $\mathbb{R}^2$.)

We set $r=\min_{0\ne\lambda\in L}\|\alpha\|$.
According to \ref{nsurfaces}(3) the set $B_r = \{v\in V | \|v\| \le
r\}$ is a convex polygon.
Since $\|\cdot\|$ is a norm, $B_r$ is {\it
balanced}, ie, invariant under the involution $x\mapsto-x$.
The definition of $r$ implies that $\inter B_r$ contains no non-zero points
of the lattice $L$. It therefore follows from
\cite[Theorem 6.21]{minkowski} that the area of $B_r$ is at most $4$.

According to Proposition \ref{nsurfaces}(3), each vertex of $B_r$,
regarded as a vector, is a scalar multiple of a boundary class of a
strict essential surface.  As $s_1$ and $s_2$ are the only slopes that
arise as boundary slopes of strict essential surfaces, there are at
most two lines through the origin that contain vertices of $B_r$. As
$B_r$ is a balanced convex polygon with non-empty interior, it must be
a parallelogram, in which two opposite vertices are multiples of an
element $\alpha_1 \in L$ such that $\langle\alpha_1\rangle = s_1$,
and the other two vertices are multiples of $\alpha_2\in L$ with
$\langle\alpha_2\rangle = s_2$.

The Dehn-filled manifold $M(\emm)$ is homeomorphic to $\Sigma$, and
hence $\pi_1(M(\emm))$ is cyclic.  As $\emm=\langle\mu\rangle$ is not
the boundary slope of a strict essential surface, it then follows from
Corollary 1.1.4 of \cite{cgls} (see Remark \ref{cyclicremark}) that
$\|\mu\|=r$, so that $\mu$ lies on the boundary of $B_r$. Since $\mu$
is a primitive element of $L$, but is not a boundary class of any
strict essential surface, it cannot be a vertex of $B_r$.  Thus $\mu$
lies on an edge of $B_r$, whose endpoints $v_1$ and $v_2$ must
respectively be scalar multiples of $\alpha_1$ and $\alpha_2$.  We may
suppose the signs of the $\alpha_i$ to be chosen in such a way that
$\alpha_i=\frac{\|\alpha_i\|}{\|\mu\|}v_i$ for $i=1,2$.

As the parallelogram $B_r$ has vertices $\pm v_1$ and $\pm v_2$, its
area is $2|\omega(v_1,v_2)|$. Hence
$$|\omega(v_1,v_2)|\le2.$$
We have $\mu=(1-t)v_1+tv_2$ for some $t\in(0,1)$. Hence
$\omega(v_1,\mu)=t\omega(v_1,v_2)$. This gives
$$q_1=|\omega(\alpha_1,\mu)|=\|\alpha_1\||\omega(v_1,\mu)|
=t\|\alpha_1\||\omega(v_1,v_2)|.$$
On the other hand,
$$\Delta=|\omega(\alpha_1,\alpha_2)|=\|\alpha_1\||\alpha_2\||\omega(v_1,v_2)|.$$
$$\frac{q_1^2}{\Delta}
  =t^2\frac{\|\alpha_1\|}{\|\alpha_2\|}|\omega(v_1,v_2)|
  \le2\frac{\|\alpha_1\|}{\|\alpha_2\|}\le{\four}2\kappa(F_1,F_2),
\leqno{\hbox{Hence}}
$$
 and the proof is complete.
\EndProof

\section{Short subgraphs I}
\label{firstbigirthsection}

This section provides the graph-theoretical background needed for
Theorem \ref{easyboundconsequence}, which is the first of our main
concrete results about two-surface knots in manifolds with cyclic
fundamental group.

\Notation
If $\Gamma_0$ is a subgraph of a graph $\Gamma$, and $r$ is a
non-negative integer, we shall denote by $N_r(\Gamma_0)$ the union of
$\Gamma_0$ with the tracks of all edge paths of length at most $r$
whose initial vertices lie in $\Gamma_0$. Thus $N_r(\Gamma_0)$ is
again a subgraph of $\Gamma$ for each $r\ge0$, and
$N_0(\Gamma_0)=\Gamma_0$. If $v$ is a vertex of $\Gamma$ we shall set
$B_r(v)=N_r(\{v\})$.
\EndNotation

\Lemma\label{tielemma}
Suppose that $\Gamma$ is a finite graph in which every vertex has
valence at least $3$, and that $\Gamma_0$ is a subgraph of
$\Gamma$. Set $\Gamma_1=N_1(\Gamma_0)$. For $i=0,1$, let $n_i$ denote
the number of valence--1 vertices of $\Gamma_i$. Assume that $n_1 <
2n_0$. Then there is a subset $t$ of ${|\Gamma_1|-|\Gamma_0|}$ with
the following properties.
\Parts
  \Part{(1)} The set $t$ is a union of vertices and (open) edges of
  $\Gamma_1$, and is closed in the subspace topology of
  $|\Gamma_1|-|\Gamma_0|$.

  \Part{(2)} If $E_t$ and $V_t$ denote respectively the number of
  edges and the number of vertices contained in $t$, then we have
  $1\le E_t - V_t \le2$.  Furthermore, if $n_1< 2n_0 - 2$ then $E_t -
  V_t = 2$.

  \Part{(3)} If $w$ denotes the number of vertices in $\bar
  t-t\subset|\Gamma_0|$, where $\bar t$ denotes the closure of $t$ in
  $|\Gamma|$, we have $\max(w,E_t)\le2(E_t - V_t)$.
\EndParts
\EndLemma

\Proof
Since $\Gamma_1=N_1(\Gamma_0)$, each edge of $\Gamma_1$
has at least one endpoint in $\Gamma_0$.  Hence if $v$ is a vertex in
$|\Gamma_1|-|\Gamma_0|$, no loop based at $v$ can lie in $|\Gamma_1|$.

We first consider the case in which some vertex
$v_0\in|\Gamma_1|-|\Gamma_0|$ has valence at least $3$ in $\Gamma_1$.
Since no loop based at $v_0$ is contained in $|\Gamma_1|$ we can choose
three distinct edges $e_1$, $e_2$ and $e_3$ having $v_0$ as an
endpoint, and the other endpoints of the $e_i$ must lie in $|\Gamma_0|$.
It follows that $t=\{v_0\}\cup e_1\cup e_2\cup e_3$ satisfies (1).
Furthermore, in the notation of (2) we have $E_t=3$ and $V_t=1$, so
that $E_t-V_t=2$, and both assertions of (2) are automatically
true. If $w$ is defined as in (3), we have $w\le3$ since each $e_i$
has $v\in|\Gamma_0|$ as an endpoint, and so (3) holds as well.

For the rest of the argument we assume that every vertex in
$|\Gamma_1|-|\Gamma_0|$ has valence at most $2$ in $\Gamma_1$.

We denote by ${\mathcal T}$ the set of all connected components of
$|\Gamma_1|-|\Gamma_0|$ which do not contain valence--1 vertices of
$\Gamma_1$. We denote by ${\mathcal T}_0\subset{\mathcal T}$ the set of all
connected components of $|\Gamma_1|-|\Gamma_0|$ which contain no vertices
whatever, and we set ${\mathcal T}_1={\mathcal T}-{\mathcal T}_0$. It is clear
that each element $\tau$ of ${\mathcal T}_0$ consists of a single edge of
$\Gamma_1$ (possibly a loop) whose  endpoints lie in $|\Gamma_0|$.

Now if $\tau$ is an element of ${\mathcal T}_1$, and $v\in\tau$ is a
vertex, then in view of the definition of $\mathcal T$, and the fact that
every vertex in $|\Gamma_1|-|\Gamma_0|$ has valence at most $2$ in
$\Gamma_1$, the valence of $v$ in $\Gamma_1$ must be exactly $2$.
Since no loop based at $v$ can lie in $|\Gamma_1|$, there are exactly
two edges $e_1$ and $e_2$ having $v$ as an endpoint, and the other
endpoints of the $e_i$ (which may or may not coincide with each other)
must lie in $|\Gamma_0|$. As $\{v\}\cup e_1\cup e_2$ is clearly open and
closed in the subspace topology of $|\Gamma_1|-|\Gamma_0|$, we must have
$\tau=\{v\}\cup e_1\cup e_2$.

To summarize, we have shown that each component $\tau\in{\mathcal T}_0$
consists of a single edge whose endpoint lies in $|\Gamma_0|$, and that
each $\tau\in{\mathcal T}_1$ consists of a single vertex $v$ and two
edges, each of which has one endpoint at $v$ and one in $|\Gamma_0|$. It
follows that for each $\tau\in\mathcal T$, if we denote by $E_\tau$ the
number of edges in $\tau$, by $V_\tau$ the number of vertices in
$\tau$, and by $w_\tau$ the number of vertices in
$\bar\tau-\tau\subset|\Gamma_0|$, then we have
\Equation\label{tielemmaone}
E_\tau-V_\tau=1
\EndEquation
and
\Equation\label{tielemmatwo}
\max(w_\tau,E_\tau)\le2 .
\EndEquation
We wish to estimate $m=\#({\mathcal T})$. If $\tau$ is any element of
${\mathcal T}_0$ then $\tau$ consists of a single edge $e$, and each
of the oriented edges with underlying edge $e$ has its initial vertex
contained in $|\Gamma_0|$. If $\tau$ is any element of ${\mathcal
T}_1$ then $\tau$ contains two edges $e_1$ and $e_2$, and each $e_i$
has a unique orientation $\omega_i$ such that $\init(\omega_i)$ lies
in $|\Gamma_0|$.  Hence for every $\tau\in{\mathcal T}$ there are
exactly two oriented edges whose underlying edges lie in $\tau$ and
whose initial points lie in $\tau$.  Thus if we denote by $\Omega_0$
the set of all oriented edges $\omega$ such that $\omega$ lies in
$|\Gamma_1|-|\Gamma_0|$ and such that no endpoint of $|\omega|$ has
valence $1$ in $|\Gamma_1|$, then $\#(\Omega_0)=2m$.

On the other hand, since every vertex of $\Gamma$ has valence at least
$3$, every valence--1 vertex of $\Gamma_0$ is the initial vertex of at
least two oriented edges $\omega_1$, $\omega_2$ in
$|\Gamma|-|\Gamma_0|$; since $\Gamma_1=N_1(\Gamma_0)$, the edge
$|\omega_i|$ lies in $|\Gamma_1|-|\Gamma_0|$ for $i=1,2$. Hence if
$\Omega$ denotes the set of all oriented edges in
$|\Gamma_1|-|\Gamma_0|$ with initial vertices in $\Gamma_0$, we have
$\#(\Omega_0)\ge2n_0$. Now $\Omega - \Omega_0$ consists of all
oriented edges whose initial points lie in $\Gamma_0$ and whose
terminal points are valence--1 vertices of $\Gamma_1$. Since each edge
of $\Gamma_1$ has at least one endpoint in $|\Gamma_0|$, we have
$\#(\Omega-\Omega_0)=n_1$. Hence $\#(\Omega_0)\ge2n_0-n_1$, ie,
\Equation\label{tielemmathree}
2m\ge2n_0-n_1 .
\EndEquation
Since by hypothesis we have $2n_0-n_1>0$, it follows from
(\ref{tielemmathree}) that $m>0$. We set $\beta=\min(2,m)$, so that
$1\le\beta\le2$. We choose a subset ${\mathcal T}'\subset{\mathcal T}$
of cardinality $\beta$. We set $t=\bigcup_{\tau\in{\mathcal
T}'}\tau$. It is clear that $t$ satisfies {(1)}. If we define $E_t$
and $V_t$ as in (2), we have $E_t=\sum_{\tau\in{\mathcal T}'} E_\tau$
and $V_t=\sum_{\tau\in{\mathcal T}'} V_\tau$. Hence it follows from
(\ref{tielemmaone}) that
$$ E_t-V_t=\beta.$$ 
Since $1\le\beta\le2$ it follows that $1\le E_t-V_t
\le2$. Furthermore, if $n_1 < 2n_0 - 2$ then (\ref{tielemmathree})
implies that $m>1$ and hence that $\beta=\max(m,2)=2$. This proves
condition (2) for our choice of $t$. Finally, using
(\ref{tielemmatwo}), we find
$$\max(w,E_t)=\max(\sum_{\tau\in{\mathcal
    T}'}w_\tau,\sum_{\tau\in{\mathcal T}'}E_t)
  \le\sum_{\tau\in{\mathcal T}'}\max(w_\tau,E_\tau)\le2\beta=2(E_t-V_t)$$
This proves condition (3).
\EndProof

\Definition
The {\it bigirth} of a graph $\Gamma$ is the infimum of the lengths of
all finite connected subgraphs $K$ of $\Gamma$ with
$\chi(|K|)<0$. Thus the bigirth of $\Gamma$ is finite if and only if
$\Gamma$ has a component $\Gamma_0$ with
$\chi(|\Gamma_0|)<0$.
\EndDefinition

\Lemma\label{usingties}
Suppose that $\Gamma$ is a graph in which each vertex has valence at
least $3$, and that $v_0$ is a vertex of $\Gamma$.  For each $r\ge1$,
let $m_r$ denote the number of valence--1 vertices of $B_r(v_0)$,
Suppose that $s$ is a positive integer.
\Parts
  \Part{(1)} If $m_{s+1} < 2m_{s}-2$, then $\bigirth(\Gamma)\le4s+4$.

  \Part{(2)} If $m_{s+1} < 2m_{s}$ then either
  $\bigirth(\Gamma)\le4s+4$, or $B_{s+1}(v_0)$ contains a connected subgraph $H$
  such that $\chi(|H|)\le0$, $\length(H)\le2s+2$, and $v_0\in |H|$.

  \Part{(3)} If $m_{s+1}< 2m_{s}$ and $B_s(v_0)$ contains a connected subgraph $H$
  such that $\chi(|H|)\le 0$, $\length(H)\le2s+2$, and $v_0\in |H|$,
  then $\bigirth(\Gamma)\le4s+4$.
\EndParts
\EndLemma

\Proof
Note that $B_{s+1}(v_0)=N_1(B_s(v_0))$. If the hypothesis of any of
the assertions (1)--(3) holds, then in particular $m_{s+1}<
2m_{s}$. Thus the hypotheses of Lemma \ref{tielemma} hold if we set
$\Gamma_0=B_s(v_0)$, $n_0=m_s$, $\Gamma_{1}=B_{s+1}(v_0)$ and
$n_{1}=m_{s+1}$. Let $t\subset|B_{s+1}(v_0)|-|B_s(v_0)|$ be a subset
satisfying conditions (1)--(3) of Lemma \ref{tielemma}. As in Lemma
\ref{tielemma}, we let $E_t$ and $V_t$ denote respectively the number
of edges and the number of vertices contained in $t$.  We set $W =
\bar t-t\subset |B_s(v_0)|$, where $\bar t$ denotes the closure of $t$
in $|\Gamma|$, and as in \ref{tielemma}(3) we let $w$ denote the
number of vertices in $W$.

For every vertex $v\in W\subset B_s(v_0)$ we choose an arc $A_v$ in
the graph $B_s(v_0)$ having length $\le s$ and endpoints $v_0$ and
$v$.  Let $K$ denote the connected subgraph of $B_s(v_0)$ with
$|K|=\bigcup_{v\in W}|A_v|$.  Since by \ref{tielemma}(3) we have
$w\le2(E_t - V_t)$, the length of $K$ is at most $2s(E_t - V_t)$.
Note that the sets $|K|$ and $t$ are disjoint since $|K|\subset
|N_s(v_0)|$.

Now let $H$ denote the subgraph of $N_{s+1}(v_0)$ with $|H|=|K|\cup
t$.  Since the closure of each component of $t$ meets $|K|$, it
follows that $|H|$ is connected.  We have $\length(H)=\length
K+E_t\le2s(E_t - V_t)+E_t$, and since $E_t\le2(E_t - V_t)$ by
\ref{tielemma}(3), we conclude that
\Equation\label{usingtiesone}
\length (H)\le2(s+1)(E_t - V_t) .
\EndEquation
Since $|K|$ and $t$ are disjoint, we have $\chi(|H|)=\chi(|K|)-(E_t -
V_t)$, where $\chi(|K|)\le 1$ since $K$ is a connected graph. Hence
\Equation\label{usingtiestwo}
\chi (|H|)\le 1-(E_t - V_t) .
\EndEquation
We now prove assertion (1) of Lemma \ref{usingties}.  If $m_{s+1}<
2m_{s}-2$, then according to \ref{tielemma}(2) we have $E_t -
V_t=2$. Hence (\ref{usingtiesone}) and (\ref{usingtiestwo}) give
${\length}(H)\le4s+4$ and $\chi(|H|)<0$. From the definition of the
bigirth it follows that $\bigirth(\Gamma)\le4s+4$, and assertion (1)
is proved.

To prove assertion (2) we recall that $E_t - V_t$ is equal to $1$ or
$2$ by \ref{tielemma}(2). If $E_t - V_t=2$ then (\ref{usingtiesone})
and (\ref{usingtiestwo}) again give $\length (H)\le4s+4$ and
$\chi(|H|)<0$, so that $\bigirth (\Gamma)\le4s+4$. If $E_t - V_t=1$
then (\ref{usingtiesone}) and (\ref{usingtiestwo}) give $\length
(H)\le2s+2$ and $\chi(|H|)\le0$. As the construction of $H$ guarantees
that $v_0\in |H|$, this completes the proof of assertion (2).

To prove assertion (3) we assume that $B_s(v_0)$ has a subgraph $H_0$
such that $\chi(|H_0|)\le0$, $\length(H_0)\le2s+2$, and $v_0\in |H_0|$.
We consider the subgraph $H'$ with $|H'|=|H|\cup |H_0|$. Since $H$ has
been shown to have length at most $2s+2$, we have
$${\length}(H')\le\length(H)+\length(H_0)\le4s+4. $$
On the other hand, we may write $|H'|=|K'|\cup t$, where $K'$ is the
subgraph of $B_s(v_0)$ such that $|K'|=|K|\cup |H_0|$. Since $K$ and
$H_0$ are both connected and both have $v_0$, as a vertex, the graph
$K'$ is connected, and hence $\chi(|K'|)\le\chi(|K|)\le0$. Now the
sets $|K'|$ and $t$ are disjoint since $|K'|\subset |N_s(v_0)|$. Since
$E_t-V_t>1$ by \ref{tielemmatwo}, it follows that
$\chi(|H'|)=\chi(|K'|)-(E_t - V_t)\le\chi(|K'|)-1<0$. We have shown
that $K'$ has length at most $4s-4$ and that $\chi(|K'|)<0$; in view
of the definition it follows that $\bigirth(\Gamma)\le4s+4$. This
proves (3).
\EndProof

\Proposition\label{trivalentcase}
Suppose that $\Gamma$ is a non-empty finite graph which has at least
two vertices, and in which every vertex has valence at least $3$.
Then
$$\bigirth(\Gamma)\le4\log_2 V ,$$
where $V$ denotes the number of vertices of $\Gamma$.
\EndProposition

\Proof
We first dispose of some degenerate cases. If some component of
$\Gamma$ has a unique vertex $v$, then since the valence of $v$ is at
least $3$ there must be at least two loops based at $v$.  Hence
$\bigirth(\Gamma)\le2$, and since by hypothesis we have $V\ge2$, the
conclusion holds in this case. If some component contains exactly two
vertices $v_1$ and $v_2$, then there is at least one edge joining
$v_1$ and $v_2$. Furthermore, since each $v_i$ has valence at least
$3$, either there are at least three edges joining $v_1$ and $v_2$, or
there are loops based both at $v_1$ and at $v_2$. In either case it
follows that $\bigirth(\Gamma)\le3$, and again the conclusion
holds. Next suppose that some component of $\Gamma$ has at least three
vertices, but that every vertex of $\Gamma$ lies in a circuit of
length at most $2$. Then there are distinct vertices $v,v_1,v_2$ such
that $v$ is joined to $v_i$ by an edge $e_i$ for $i=1,2$. For $i=1,2$,
choose a circuit $C_i$ of length $\le2$ containing $v_i$, and let $H$
denote the subgraph with $|H|=|C_1|\cup |C_2|\cup\bar e_1\cup\bar
e_2$, so that $\length(H)\le6$. If $C_1\ne C_2$ then there are two
distinct circuits in $H$, and hence $\chi(|H|)<0$. If $C_1=C_2$, then
in $H$ we have the circuit $C_1$ and the arc $A$ with $|A|=\bar
e_1\cup\bar e_2$, which has its endpoints in $C_1$ but contains the
vertex $v\notin C_1$. Hence we have $\chi(|H|)<0$ in this subcase as
well. Hence $\bigirth(\Gamma)\le6$, and since $V\ge3$ in this case,
the conclusion again holds.

Hence we may assume that there is a vertex $v$ of $\Gamma$ which has
valence at least 3 and does not lie in any circuit of length
$\le2$. It follows that the subgraph $B_1(v)$ has at least three
vertices of valence $1$.

For each $r\ge1$, we let $X_r$ denote the set of all valence--1
vertices of $B_r(v)$, and we set $m_r=\#(X_r)$.  Thus $m_1\ge3$. Note
that a valence--1 vertex of $B_r(v)$ has valence at least $3$ in
$B_{r'}(v)$ for every $r'>r$. In particular, the $X_r$ are pairwise
disjoint for $r\ge1$. Since none of the $X_r$ contains $v$, the number
$V$ of vertices of $\Gamma$ is at least $1+\sum_{r=1}^\infty m_r$,
where $m_r=0$ for large enough $r$.

Since $m_r=0$ for large $r$, there is a smallest non-negative integer
$s$ such that $m_{s+1}<2^{s+1}$. Since $m_1\ge3$, we have
$s\ge1$. The minimality of $s$ then implies that $m_{r}\ge2^r$ for
$r=2,\ldots,s$. Hence
$$V\ge1+\sum_{r=1}^s m_r\ge1+3+\sum_{r=2}^s2^r=2^{s+1}.$$
This means that $s+1\le\log_2V$. Hence to prove the proposition it
suffices to show that $\bigirth(\Gamma)\le4s+4$.

We distinguish two cases according to whether $m_s\ge2^s+1$ or
$m_s=2^s$. First suppose that $m_s\ge2^s+1$. Then since
$m_{s+1}<2^{s+1}$, we have $m_{s+1} < 2m_{s}-2$. It therefore follows
from Lemma \ref{usingties}(1) that $\bigirth(\Gamma)\le4s+4$, as
required.

Now suppose that $m_s=2^s$. Note that in this case there must be an
integer $s'$, with $1\le s'<s$, such that $m_{s'+1}<2m_{s'}$. Indeed,
if we had $m_{r+1}\ge2m_r $ for $r=1,\ldots,s-1$, we would have
$m_s\ge2^{s-1}m_1=3\cdot2^{s-1}$, a contradiction to $m_s=2^s$.

If $s'$ satisfies $1\le s'<s$ and $m_{s'+1}<2m_s$, then by Lemma
\ref{usingties}(2), there is a subgraph $H$ of $B_{s'+1}(v)\subset
B_s(v)$ such that $\chi(|H|)\le0$, $\length(H)\le2s'+2<2s+2$, and
$v\in |H|$. On the other hand, we have $2m_{s}-m_{s+1}>2$ since
$m_s=2^s$ and $m_{s+1}<2^{s+1}$. It therefore follows from Lemma
\ref{usingties}(3) that $\bigirth(\Gamma)\le4s+4$ in this case
as well. 
\EndProof

\Lemma\label{roughlyequallengths}
Suppose that $\Gamma$ is a finite graph with $\chi(|\Gamma|)<0$. Set
$\alpha=\length (\Gamma)/|\chi(\Gamma)|$. Then $\Gamma$ has a subgraph
$\Gamma_0$ such that
\Parts
  \Part{(1)} $\chi(|\Gamma_0|)<0$,
  \Part{(2)} every vertex of $\Gamma_0$ has valence at least $2$,
  \Part{(3)} every component of $|\Gamma_0|$ contains a vertex whose
  valence in $\Gamma_0$ is at least $3$,
  \Part{(4)} if $\mathcal V$ denotes the set of all vertices of valence
  $\ge3$ in $\Gamma_0$, every component of $|\Gamma_0|-{\mathcal V}$ is
  homeomorphic to an open interval and contains
  at most $\lfloor\alpha\rfloor$ edges.
\EndParts
\EndLemma

\Proof
Let $\mathcal G$ denote the set of all subgraphs $G$ of $\Gamma$ such
that (a) $\chi(|G|)<0$ and (b) $\length(G)/|\chi(|G|)|\le\alpha$. Note
that $\Gamma\in{\mathcal G}$. Let $\Gamma_0$ be a subgraph in
$\mathcal G$ which is minimal with respect to inclusion. Then
$\Gamma_0$ satisfies (1), and we shall complete the proof by showing
that it satisfies (2)--(4) as well.

If $\Gamma_0$ has an isolated vertex $v$, then the subgraph $G$
defined by $|G|=\Gamma_0-\{v\}$ satisfies
$\chi(|G|)=\chi(|\Gamma|)-1$ and $\length(G)=\length(\Gamma)$. Hence
$G\in{\mathcal G}$, and the minimality of $\Gamma_0$ is
contradicted. If $\Gamma_0$ has a vertex $v$ of valence $1$, and if
$e$ is the edge of $\Gamma_0$ incident to $v$, then the subgraph $G$
defined by $|G|=|\Gamma_0|-(\{v\}\cup e)$ satisfies
$\chi(|G|)=\chi(|\Gamma|)$ and $\length(G)=\length(\Gamma)-1$. Again
it follows that $G\in{\mathcal G}$, in contradiction to the minimality
of $\Gamma_0$. This proves (2). If (3) fails, it now follows that
there is a component $C$ of $\Gamma_0$ whose vertices all have valence
$2$. In this case the subgraph $G$ defined by $|G|=|\Gamma_0|-|C|$
satisfies $\chi(|G|)=\chi(|\Gamma|)$ and $\length(G)<\length(\Gamma)$,
and again we have a contradiction to the minimality of $\Gamma_0$.
This proves (3).

To prove (4), we first note that every vertex in $|\Gamma_0|-{\mathcal
V}$ has valence $2$.  Since every component of $\Gamma_0$ has a vertex
of valence $3$ it follows that every connected component of
$|\Gamma_0|-{\mathcal V}$ is homeomorphic to an open interval.  Now
suppose that some component $C$ of $|\Gamma_0|-{\mathcal V}$ contains
$m>\alpha$ edges, and we consider the subgraph $G$ defined by
$|G|=|\Gamma_0|-C$. Note that since $\Gamma_0\in{\mathcal G}$ we have
$\chi(|\Gamma_0|)<0$; this implies that $\Gamma_0$ has at least one
vertex of valence $\ge3$, and hence that $\length(G)>0$.

Now set $c=-\chi(|\Gamma_0|)$, so that $c-1=-\chi(|G|)$. We have
$c-1\ge0$ since $\chi(|\Gamma_0|)<0$. Since $\Gamma_0\in{\mathcal G}$
we have $\length(\Gamma_0)/|\chi(|\Gamma_0|)|\le\alpha$, ie,
$\length(\Gamma_0)\le c\alpha$. Hence
\Equation\label{roughlyequallengthsone}
\length(G)=\length(\Gamma_0)-m\le c\alpha-m<\alpha(c-1) .
\EndEquation
Since $\length(G)>0$ it follows from (\ref{roughlyequallengthsone})
that $-\chi(|G|)=c-1>0$. From (\ref{roughlyequallengthsone}) we then
conclude that
$$\length(G)<\alpha|\chi(|G|)|.$$
It now follows that $G\in{\mathcal G}$,
and again we have a contradiction to the minimality of $\Gamma_0$.
This establishes (4).
\EndProof

\Lemma\label{subgraphs}
Let $\Gamma$ be a finite graph with no simply connected components.
Then we have $\chi(|\Gamma_0|) \ge \chi(|\Gamma|)$ for any subgraph
$\Gamma_0$ of $\Gamma$.
\EndLemma

\Proof
In this proof we shall denote the set of vertices of a finite graph
$G$ by ${\mathcal V}(G)$, and the set of its edges by ${\mathcal
E}(G)$; thus $\chi(G)=\#({\mathcal V}(G))-\#({\mathcal E}(G))$. The
conclusion of the lemma is equivalent to the assertion that
$\#({\mathcal V}(\Gamma)-{\mathcal V}(\Gamma_0))\le\#({\mathcal
E}(\Gamma)- {\mathcal E}(\Gamma_0))$.

We first consider the special case in which every vertex in the set
${\mathcal V}(\Gamma)- {\mathcal V}(\Gamma_0)$ has valence at least
$2$ in $\Gamma$. Since every edge having an endpoint in the set
${\mathcal V}(\Gamma)- {\mathcal V}(\Gamma_0)$ must belong to the set
${\mathcal E}(\Gamma)- {\mathcal E}(\Gamma_0)$, we have
$$
2\#({\mathcal V}(\Gamma)- {\mathcal
V}(\Gamma_0))\le\sum_{v\in{\mathcal V}(\Gamma)- {\mathcal
V}(\Gamma_0)}{\rm valence}(v)\le2 \#({\mathcal E}(\Gamma)-
{\mathcal E}(\Gamma_0)) ,
$$
which implies the conclusion.

To prove the lemma in general, we use induction on $\#({\mathcal
V}(\Gamma))$. If $\#({\mathcal V}(\Gamma))=0$ the assertion is
trivial. Suppose that $n>0$, that the assertion is true for graphs
with $n-1$ vertices, that $\Gamma$ is a graph with $\#({\mathcal
V}(\Gamma))=n$ and having no simply-connected components, and that
$\Gamma_0$ is a subgraph of $\Gamma$. If every vertex in ${\mathcal
V}(\Gamma)- {\mathcal V}(\Gamma_0)$ has valence at least $2$ in
$\Gamma$ we are in the special case already proved. No vertex of
$\Gamma$ has valence $0$ because no component of $\Gamma$ is simply
connected. Hence we may assume that some vertex $v_1\in{\mathcal
V}(\Gamma)- {\mathcal V}(\Gamma_0)$ has valence $1$ in $\Gamma$. If
$\Gamma'$ denotes the graph obtained from $\Gamma$ by removing $v_1$
and the unique edge incident to $v_1$, then $\Gamma_0$ is a subgraph
of $\Gamma'$. Since $\#({\mathcal V}(\Gamma))=n-1$, the induction
hypothesis gives $\chi(\Gamma_0)\ge\chi(\Gamma')=\chi(\Gamma).$
\EndProof

\Proposition\label{generalbigirthbound}
Suppose that $\Gamma$ is a finite graph such that no component of
$|\Gamma|$ is simply connected and $\chi(|\Gamma|)<0$.  Then
$$
\bigirth(\Gamma)\le4(\log_2|2\chi(|\Gamma|)|)\
\left\lfloor\frac{\length(\Gamma)}{|\chi(|\Gamma|)|}\right\rfloor .
$$
\EndProposition

\Proof
Set $\alpha=\length (\Gamma)/|\chi(|\Gamma|)|$. Let $\Gamma_0$ be a
subgraph of $\Gamma$ satisfying conditions (1)--(4) of Lemma
\ref{roughlyequallengths}.  Conditions (2)--(4) imply that there
exists a graph $\Gamma_0^*$ such that $|\Gamma_0^*|=|\Gamma_0|$ and
every vertex of $\Gamma_0^*$ has valence at least $3$.  (If, as in the
statement of condition (4), we let $\mathcal V$ denote the set of all
vertices of valence at least $3$ in $\Gamma_0$, then the vertices of
$\Gamma_0^*$ are the vertices in $\mathcal V$ and the edges of
$\Gamma_0^*$ are the connected components of $|\Gamma_0|-{\mathcal V}$.)

By condition (1) we have $\chi(|\Gamma_0^*|)<0$, and in particular
$|\Gamma_0^*|\ne\emptyset$. By Proposition \ref{trivalentcase}, if
$V^*$ denotes the number of vertices of $\Gamma_0^*$, we have either
$V=1$ or
$$\bigirth(\Gamma_0^*)\le4\log_2 V^*.$$
Since $\Gamma_0^*$ has $V^*$ vertices, all of valence at least $3$, it
must have at least $3V^*/2$ edges. Hence $\chi(|\Gamma_0^*|)\le
V^*-(3V^*/2)=-V^*/2$. On the other hand, $|\Gamma_0^*|$ is
homeomorphic to $|\Gamma_0|\subset|\Gamma|$, and by Lemma
\ref{subgraphs} $0>\chi(|\Gamma_0^*|)= \chi(|\Gamma_0|)
\ge\chi(|\Gamma|)$.  Hence $V^*\le|2\chi(|\Gamma|)|$, so that
\Equation\label{generalbigirthone}
\bigirth(\Gamma_0^*)\le4\log_2 |2\chi(|\Gamma|)| ,
\EndEquation
provided that $V\ne1$. However, if $V=1$, then since the valence of
the unique vertex $v$ is at least $3$ there must be at least two loops
based at $v$; hence $\bigirth(\Gamma)\le2$, whereas
$\chi(|\Gamma|)\le-1$.  Hence (\ref{generalbigirthone}) holds in all
cases. By definition this means that there is a subgraph $H^*$ of
$\Gamma_0^*$ with $\chi(|H|)<0$ and $\length(H)\le4\log_2
|2\chi(|\Gamma|)|$.

According to condition (4) of Lemma \ref{roughlyequallengths}, every
edge of $\Gamma_0^*$ contains at most $\lfloor\alpha\rfloor$ edges of
the subgraph $\Gamma_0$ of $\Gamma$.  Hence if $H$ denotes the
subgraph of $\Gamma_0\subset\Gamma$ with $|H|=|H^*|$, then
$$\length(H^*)\le\lfloor\alpha\rfloor\length(H)\le
4\lfloor\alpha\rfloor\log_2 |2\chi(|\Gamma|)| .
$$
Since $H$ is in particular a subgraph of $\Gamma$ with $\chi(|H|)<0$,
it follows that
$$\bigirth(\Gamma)\le4\lfloor\alpha\rfloor\log_2 2|\chi(|\Gamma|)| .$$
This is the conclusion of the proposition.
\EndProof

\section{Slopes and genera I}
\label{firstgenussection}

The goal of this section is to prove Theorem \ref{easyboundconsequence}.

\Definition
Suppose that ${\mathcal A}$ is a properly embedded $1$--manifold in a
compact orientable $2$--manifold $F$. A finite graph $\Gamma$ with
$|\Gamma|\subset \inter F$ will be called a {\it dual graph} of
${\mathcal A}$ in $F$ if it has the following properties.
\Parts
  \Part{(1)} Every edge of $\Gamma$ meets ${\mathcal A}$ transversally.
  \Part{(2)} Every component of $F-{\mathcal A}$ contains a unique
  vertex of $\Gamma$.
  \Part{(3)} There is a bijective correspondence $A\mapsto e_A$
  between the components of $\mathcal A$ and the edges of $\Gamma$,
  such that for each component $A$ of $\mathcal A$ we have
  $A\cap|\Gamma|=\{m_A\}$, where $m_A$ denotes the midpoint of $e_A$.
\EndParts
Note that every properly embedded $1$--manifold ${\mathcal A}$ in $F$ has a dual
graph $\Gamma$, and it is unique up to ambient isotopy.  Note also that
$|\Gamma|$ is a retract of $F$, so that in particular $|\Gamma|$ is
$\pi_1$--injective.  Furthermore, $|\Gamma|$ has the same number of connected
components as $F$.
\EndDefinition

\Lemma\label{dualgraphandgenus}
Let $\Gamma$ denote a dual graph of a properly embedded $1$--manifold
$\mathcal A$ in a compact, connected, orientable $2$--manifold
$F$. Suppose that every component of $F-{\mathcal A}$ is a planar
surface. Then the first Betti number of $|\Gamma|$ is greater than or
equal to the genus of $F$.
\EndLemma

\Proof
We shall argue by induction on the number of components of ${\mathcal
A}$. If ${\mathcal A}=\emptyset$ the assertion is true because a
planar surface has genus $0$ by definition. Now suppose that
${\mathcal A}$ has $\nu>0$ components, and that the lemma is true for
all properly embedded $1$--manifolds with fewer than $\nu$ components
in compact, connected, orientable $2$--manifolds. Let $b$ denote the
first Betti number of $\Gamma$, and let $g$ denote the genus of
$F$. Choose a component $A$ of ${\mathcal A}$, let $V$ denote a collar
neighborhood of $A$ in $F$, and let $F'$ denote the closure of $F-V$,
let $e$ denote the unique edge of $\Gamma$ which meets $A$, and let
$\Gamma'$ denote the subgraph of $\Gamma$ with
$|\Gamma'|=|\Gamma|-e$. Then ${\mathcal A}'={\mathcal A}-A$ is a
properly embedded $1$--manifold in $F'$, and $\Gamma'$ is a dual graph
of ${\mathcal A}'$ in $F'$.

If $F'$ is connected then $|\Gamma'|$ is connected and has first Betti
number $b-1$. Furthermore, in this case the genus of $F'$ is at least
$g-1$. Hence the induction hypothesis implies that $b-1\ge g-1$, so
that $b\ge g$. Now suppose that $F'$ is disconnected, let $F_1'$ and
$F_2'$ denote its components. For $i=1,2$, let $\Gamma_i'$ denote the
component of $\Gamma'$ such that $|\Gamma_i'|\subset F_i'$, and let
$g_i'$ and $b_i'$ denote respectively the genus of $F_i'$ and
the first Betti number of $|\Gamma_i'|$. Then $g=g_1'+g_2'$ and
$b=b_1'+b_2'$. But the induction hypothesis gives $b_i'\ge g_i'$ for
$i=1,2$, and hence $b\ge g$.
\EndProof

\Lemma\label{curvesdontmatter}
Suppose that ${\mathcal A}$ is a properly embedded $1$--manifold in a
compact, connected, orientable $2$--manifold $F$. Let ${\mathcal A}_0$
denote the union of all those components of ${\mathcal A}$ that are
arcs.  Assume that a dual graph $\Gamma$ of ${\mathcal A}_0$ in $F$
has finite bigirth. Then there is a compact $\pi_1$--injective
$1$--dimensional polyhedron $K\subset F$ of Betti number $2$ such that
$\#(K\cap {\mathcal A})\le\bigirth(\Gamma)$.
\EndLemma
 
\Proof
First of all, observe that if ${\mathcal A}^*\supset {\mathcal A}_0$
denotes the union of all components of ${\mathcal A}$ that are arcs or
homotopically non-trivial simple closed curves, then ${\mathcal
A}-{\mathcal A}^*$ is contained in a union of disjoint disks
$D_1,\ldots,D_k$ which are disjoint from ${\mathcal A}^*$.  If
$K^*\subset F$ is a compact $\pi_1$--injective $1$--dimensional
polyhedron of Betti number $2$ such that $\#(K^*\cap {\mathcal
A}^*)\le\bigirth(\Gamma)$, there is an isotopy of $F$ supported on
$D_1\cup\cdots\cup D_k$ which carries $K^*$ onto a polyhedron $K$
disjoint from ${\mathcal A}-{\mathcal A}^*$, and it follows that
$\#(K\cap {\mathcal A})\le\bigirth(\Gamma)$.  Hence we may assume
without loss of generality that ${\mathcal A}$ contains no
homotopically trivial simple closed curves.

Since $\bigirth(\Gamma)<\infty$, there is a connected subgraph
${\mathcal H}$ of $\Gamma$ with Betti number $2$ 
such that $\length({\mathcal H}) =
\bigirth(\Gamma)$.  Since the
polyhedron $|\Gamma|\subset F$ is $\pi_1$--injective, $|{\mathcal H}|$
is also $\pi_1$--injective in $F$. It follows that $\pi_1(F)$ is
non-abelian and hence that $F$ is not a torus.  Note also that the
definition of bigirth implies that ${\mathcal H}$ has minimal length
among all connected subgraphs of $\Gamma$ with Betti number $2$; hence
${\mathcal H}$ has no valence--1 vertices.  Since $\chi(|{\mathcal
H}|)=-1$ it follows that every vertex of ${\mathcal H}$ has valence at
least $2$, and that no vertex of ${\mathcal H}$ has valence greater
than $4$.

Since $\length({\mathcal H}) = \bigirth(\Gamma)$, it follows from the
definition of the dual graph that
$\#(|{\mathcal H}|\cap {\mathcal A}_0)=\bigirth(\Gamma) $.

Let $C={\mathcal A}-{\mathcal A}_0$ denote the union of all components
of ${\mathcal A}$ that are simple closed curves. Since $F$ is not a
torus, any two components of $C$ can cobound at most one
annulus. Hence if $Z$ denotes the union of $C$ with all annuli whose
boundaries are contained in $C$, the components of a regular
neighborhood of $Z$ are themselves annuli.  Clearly $Z$ has a regular
neighborhood which is disjoint from ${\mathcal A}_0$.  Define the
complexity of a regular neighborhood $N$ of $Z$ to be the pair $(p,v)$
where $p=\#(\partial N\cap|{\mathcal H}|)$ and $v$ is the number of
vertices of ${\mathcal H}$ which lie in $F-N$.  Among all regular
neighborhoods of $Z$ that are disjoint from ${\mathcal A}_0$ choose
one, $\mathcal Z$, which has minimal complexity with respect to
lexicographical order.  In particular ${\mathcal H}$ and $\partial
{\mathcal Z}$ intersect transversally.

It follows from the construction of $\mathcal Z$ that ${\mathcal
A}\cap\partial{\mathcal Z} =\emptyset$. Another consequence of the
construction is that any annulus cobounded by two curves in
$\partial{\mathcal Z}$ is itself a component of $\mathcal Z$. Note
also that since ${\mathcal A}$ contains no homotopically trivial
simple closed curves, the annuli that make up $\mathcal Z$ are
homotopically non-trivial.

Since $|{\mathcal H}|$ is $\pi_1$--injective and the components of
$\mathcal Z$ are annuli, $|{\mathcal H}|$ cannot be contained in
$\mathcal Z$. If $|{\mathcal H}|\cap{\mathcal Z}=\emptyset$ then
$\#(|{\mathcal H}|\cap {\mathcal A})=\#(|{\mathcal H}|\cap {\mathcal
A}_0)=\bigirth(\Gamma)$, and the conclusion of the lemma follows if we
set $K=|{\mathcal H}|$.  Hence we may assume that there is a component
$X$ of $|{\mathcal H}|\cap\overline{F-{\mathcal Z}}$ such that
$X\cap\partial{\mathcal Z}\ne\emptyset$.

Consider the case in which $X\supset|E|$ for some circuit $E$ of
${\mathcal H}$. Let $\beta\subset X$ be an arc having one endpoint in
$|E|$ and one in $\partial{\mathcal Z}$. Let $c$ denote the component
of $\partial{\mathcal Z}$ containing an endpoint of $\beta$.  In this
case we shall show that the polyhedron $K=|E|\cup\beta\cup c$ has the
properties stated in the lemma.

Since $E$ is a circuit in the dual graph $\Gamma$ of ${\mathcal A}_0$,
at least one edge of $\Gamma$ is contained in $|E|$; hence there is a
component $\alpha$ of ${\mathcal A}_0$ that meets $|E|$ in exactly one
point.  Since $c\cap\alpha\subset{\mathcal Z}\cap{\mathcal
A}_0=\emptyset$ it follows that the homology class $[|E|]\in
H_1(F;\mathbb{Z}/2\mathbb{Z})$ is not a multiple of $[c]$. But $c$ is
homotopically non-trivial in $F$ since $\mathcal Z$ is made up of
non-trivial annuli, and hence $K=|E|\cup\beta\cup c$ is
$\pi_1$--injective. On the other hand we have $K\cap {\mathcal
A}\subset |{\mathcal H}|\cap {\mathcal A}_0$ since $\partial{\mathcal
Z}\cap{\mathcal A}=\emptyset$ and ${\mathcal A}-{\mathcal A}_0\subset
\inter({\mathcal Z})$. Hence $\#(K\cap {\mathcal A})\le\#(|{\mathcal
H}|\cap {\mathcal A}_0)=\bigirth(\Gamma)$, as required.

There remains the case in which there is no circuit $E$ of ${\mathcal
H}$ such that $X\supset|E|$.  Then $X$ is homeomorphic to a finite
tree. In this case we let $\beta$ denote an arc in $X$ that joins two
endpoints of $X$. If $X\cap {\mathcal A}_0$ happens to be non-empty,
we choose $\beta$ to contain a point of $X\cap {\mathcal A}_0$; this
is possible because every point in a finite tree is contained in an
arc joining two endpoints of the tree.

We let $c$ denote the union of all components (there are at most two)
of $\partial{\mathcal Z}$ that contain endpoints of $\beta$. In this
case we shall show that the polyhedron $K=\beta\cup c$ has the
properties stated in the lemma.  We have $K\cap {\mathcal A}\subset
|{\mathcal H}|\cap {\mathcal A}_0$ since $\partial{\mathcal
Z}\cap{\mathcal A}=\emptyset$ and ${\mathcal A}-{\mathcal A}_0\subset
\inter({\mathcal Z})$. Hence $\#(K\cap {\mathcal A})\le\#(|{\mathcal
H}|\cap {\mathcal A}_0)=\bigirth(\Gamma)$, so we need only prove that
$K$ is $\pi_1$--injective in $F$.

There are several subcases. If $c$ has two components which do not
cobound an annulus in $F$ then $K$ is automatically $\pi_1$--injective.
If $c$ has two components $c_1$ and $c_2$ which do cobound an annulus
in $F$, it follows from the construction of $C$ that $c_1$ and $c_2$
are the two boundary components of a single annulus component of
$\mathcal Z$, say ${\mathcal Z}_0$. Thus $K$ is the union of
$\partial{\mathcal Z}_0$ with the properly embedded arc $\beta$ in
$F-\inter{\mathcal Z}_0$. Since $F$ is not a torus, $K$ is
$\pi_1$--injective.

In the remaining subcases, $c$ will be connected, and since $\partial
C$ is made up of non-trivial annuli, $c$ is a homotopically
non-trivial curve. Hence in order to establish $\pi_1$--injectivity it
suffices to show that there is no arc $\delta\subset C$ such that
$\partial\delta=\partial\beta$ and such that $\delta\cup\beta$ bounds
a $D\subset F-\inter{\mathcal Z}$. In these subcases we shall assume
such an arc $\delta$ and disk $D$ exist, and derive a contradiction.

Consider the subcase in which $X\cap {\mathcal A}_0\ne\emptyset$. It
then follows from our choice of $\beta$ that $\beta\cap {\mathcal
A}_0\ne\emptyset$.  Let $\alpha$ denote a component of ${\mathcal
A}_0$ that meets $\beta$. Since $\beta$ is contained in the dual graph
$\Gamma$, it must meet $\alpha$ transversally in exactly one
point. But $\delta\cap\alpha\subset{\mathcal Z}\cap {\mathcal
A}_0=\emptyset$, so that $\delta\cup\beta$ meets $\alpha$
transversally in exactly one point. This gives a contradiction if
$\delta\cup\beta$ bounds a disk $D$.

Now suppose that $X\cap {\mathcal A}_0=\emptyset$. Since $\beta\cup
{\mathcal A}_0\subset X\cap {\mathcal A}_0=\emptyset$ and $\delta\cup
{\mathcal A}_0\subset{\mathcal Z}\cap {\mathcal A}_0=\emptyset$, we
have $D\cap {\mathcal A}_0=\emptyset$.

On the other hand, since $X\cap {\mathcal A}_0=\emptyset$, it follows
in particular that $X$ cannot contain any edge of the subgraph
${\mathcal H}$ of $\Gamma$. Hence either $X$ is an interior arc of an
edge of ${\mathcal H}$, or $X$ is a connected subset of the open star
of some vertex $v$ of ${\mathcal H}$, and $v\in X$. In any event,
since every vertex of ${\mathcal H}$ has valence $2$, $3$ or $4$, $X$
is either a topological arc, or a cone over a three- or four-point
set.  If $X$ is an arc we must have $X=\beta$. If $X$ is a cone over a
three- or four-point set then the cone point is a vertex $v$ of
${\mathcal H}$ which is an interior point of $\beta$.  Furthermore
either $X=\beta\cup\epsilon_1$ or
$X=\beta\cup\epsilon_1\cup\epsilon_2$, where each $\epsilon_i$ is an
arc having $v$ as an endpoint and $\epsilon_i\cap\beta = \{v\}$.  Note
that for each $i$ either $\epsilon_i$ is a properly embedded arc in
the disk $D$, or $\epsilon_i\cap D=\{v\}$.

Let ${\mathcal Z}'$ denote a small regular neighborhood of ${\mathcal
Z}\cup D$ in $F$. Then ${\mathcal Z'}$ is also a regular neighborhood
of $Z$, and it is disjoint from ${\mathcal A}_0$ since $\mathcal Z$
and $D$ are both disjoint from ${\mathcal A}_0$.  If $X=\beta$ we have
$\#(\partial{\mathcal Z}'\cap |{\mathcal H}|)=\#(\partial{\mathcal
Z}\cap |{\mathcal H}|)-2$.  If $X$ is a cone and one of the arcs
$\epsilon_i$ is contained in $D$, then $\#(\partial{\mathcal Z}'\cap
|{\mathcal H}|)\le \#(\partial{\mathcal Z}\cap |{\mathcal H}|)-2$.  If
$X$ is a cone over a three-point set and $\epsilon_1$ meets $D$ only
at the vertex $v$, we may choose the regular neighborhood ${\mathcal
Z}'$ of ${\mathcal Z}\cup D$ so that $\#(\partial{\mathcal Z}'\cap
|{\mathcal H}|)=\#(\partial{\mathcal Z}\cap |{\mathcal H}|)-1$.  If
$X$ is a cone over a four-point set and both of the arcs $\epsilon_1$
and $\epsilon_2$ meet $D$ only at the vertex $v$, then we may choose
the regular neighborhood ${\mathcal Z}'$ of ${\mathcal Z}\cup D$ so
that $\#(\partial{\mathcal Z}'\cap |{\mathcal
H}|)=\#(\partial{\mathcal Z}\cap |{\mathcal H}|)$ but ${\mathcal Z}'$
contains strictly more vertices of ${\mathcal H}$ than ${\mathcal Z}$
does.  In all cases ${\mathcal Z}'$ would have lower complexity than
$\mathcal Z$, contradiction.
\EndProof

\Proposition\label{easykappabound}
Suppose that $F_1$ and $F_2$ are connected essential surfaces in an
irreducible knot manifold $M$. Suppose that $F_1$ and $F_2$ intersect
transversally, and that $\partial F_1$ and $\partial F_2$ are
non-empty and intersect minimally in the sense of
\ref{minimalintersection}. Assume that every component of
$F_2-(F_1\cap F_2)$ is a disk or an annulus. For $i=1,2$, let $g_i$,
$s_i$ and $m_i$ denote, respectively, the genus, boundary slope and
number of boundary components of $F_i$.  Assume that $g_2\ge2$. Then
there is a compact $\pi_1$--injective $1$--dimensional polyhedron
$K\subset F_2$ of Betti number $2$ such that
$$\#(K\cap F_1)\le\frac{2m_1m_2\Delta(s_1,s_2)\log_2(2g_2-2)}{g_2-1}.$$
In particular, according to Definition \ref{kappadef}, we have
$$\kappa(F_1,F_2)\le\frac{2m_2^2\Delta(s_1,s_2)\log_2(2g_2-2)}{g_2-1}.$$
\EndProposition

\Proof
Set ${\mathcal A} =F_1\cap F_2$, so that ${\mathcal A} $ is a properly
embedded $1$--manifold in $F_2$. Let ${\mathcal A} _0$ denote the union
of all those components of ${\mathcal A} $ that are arcs.

We claim that the components of $F_2-{\mathcal A} _0$ are disks and
annuli. Indeed, if $Z$ is any component of $F_2-{\mathcal A} _0$, then
$Z\cap {\mathcal A} $ is a union of simple closed curve components of
${\mathcal A} $. Each component of $Z-(Z\cap {\mathcal A} )$ is a
component of $F_2-{\mathcal A} $, and is therefore a disk or annulus
according to the hypothesis. It follows that $\chi(Z)\ge0$. Since
$\partial F_2\ne\emptyset$, it follows that $Z$ is a disk or an
annulus.

In particular the components of $F_2-{\mathcal A} _0$ are planar
surfaces. Hence if $\Gamma$ denotes a dual graph of ${\mathcal A} _0$
in $F$, it follows from Lemma \ref{dualgraphandgenus} that the first
Betti number of $|\Gamma|$ is greater than or equal $g_2\ge2$. We
therefore have $\chi(|\Gamma|)<0$ and $|\chi(|\Gamma|)|\ge
g_2-1$. Since $g_2-1$ and $\chi(|\Gamma|)$ are positive integers, it
follows that
$$\frac{\log_2|2\chi(|\Gamma|)|}{\chi(|\Gamma|)}\le
  \frac{\log_2(2(g_2-1))}{g_2-1}.$$
Moreover, $\Gamma$ is connected since $F_2$ is connected.

The length of $\Gamma$ is equal to the number of components of
${\mathcal A} _0$.  Since $\partial F_1$ and $\partial F_2$ intersect
minimally, it follows that
$$2\length(\Gamma)=\#((\partial F_1)\cap(\partial F_2))
  =m_1m_2\Delta(s_1,s_2) .$$
Since $\chi(\Gamma)<0$ and $\Gamma$ is connected, it follows from
Proposition \ref{generalbigirthbound} that
\Equation\label{easykappaboundone}
\bigirth(\Gamma) \le
4(\log_2|2\chi(\Gamma)|)\frac{\length(\Gamma)}{|\chi(\Gamma)|} \le
\frac{2\log_2(2(g_2-1))}{g_2-1}\cdot m_1m_2\Delta(s_1,s_2) .
\EndEquation
But according to Lemma \ref{curvesdontmatter}, there is a compact
$\pi_1$--injective $1$--dimensional polyhedron $K\subset F$ of Betti
number $2$ such that $\#(K\cap {\mathcal A} )\le\bigirth(\Gamma)$,
and the first assertion of the proposition therefore follows by
(\ref{easykappaboundone}).  Since we must have $t_{F_2}(K)\ge1$, the
second assertion follows from the first one together with Definition
\ref{kappadef}. 
\EndProof

\Theorem\label{easyboundconsequence}
Suppose that $K$ is a non-exceptional two-surface knot in a closed,
orientable $3$--manifold $\Sigma$ such that $\pi_1(\Sigma)$ is cyclic.
Let $\emm$ denote the meridian slope of $K$ and let $F_1$ and $F_2$ be
representatives of the two isotopy classes of connected strict
essential surfaces in $M(K)$.  Let $s_i$, $g_i$ and $m_i$ denote,
respectively, the boundary slope (well-defined by
\ref{twosurfprops}), the genus and the number of boundary components
of $F_i$. Assume that $s_2\ne\emm$ and that $g_2\ge2$.  Set
$q_i=\Delta(s_i,\emm)$ (so that $q_i$ is the denominator of $s_i$ in
the sense of \ref{denominator}), and set $\Delta=\Delta(s_1,s_2)$
(so that $\Delta\not=0$ by \ref{twosurfprops}).
Then
$$\left(\frac{q_1}{\Delta}\right)^2\le\frac{4 m_2^2\log_2(2g_2-2)}{g_2-1} .$$
\EndTheorem

\Proof
We may assume after an isotopy that $F_1$ and $F_2$ intersect
transversally, and that $\partial F_1$ and $\partial F_2$ intersect
minimally in the sense of \ref{minimalintersection}. Furthermore,
$F_1$ and $F_2$ may be assumed to be chosen within their rel-boundary
isotopy classes so as to minimize the number of components of $F_1\cap
F_2$. Then no component of $F_1\cap F_2$ is a homotopically trivial
simple closed curve (cf  Remark \ref{standardremark}). Set ${\mathcal
A}=F_1\cap F_2$.  It now follows from Theorem \ref{disksandannuli}
that every component of $(\inter F_i)-{\mathcal A}$ is an open disk or
an open annulus.  Hence by Proposition \ref{easykappabound} we have
\Equation\label{easyboundconsequenceone}
\kappa(F_1,F_2)\le\frac{2m_2^2\Delta\log_2(2g_2-2)}{g_2-1} .
\EndEquation
On the other hand, according
to Theorem \ref{generalknotinequality} we have
\Equation\label{easyboundconsequencetwo}
\frac{q_1^2}{\Delta}\le{\four}2\kappa(F_1,F_2) .
\EndEquation
The inequality in the
conclusion of Theorem \ref{easyboundconsequence} follows from
(\ref{easyboundconsequenceone}) and (\ref{easyboundconsequencetwo}).
\EndProof

If $K$ is a knot in a homology $3$--sphere with $M(K)$ irreducible
then it follows from Remark \ref{coho} that $M(K)$ has an essential
surface whose numerical boundary slope with respect to a standard
framing is $0$.

\Corollary\label{easycorollary}
Suppose that $K$ is a knot in a homotopy $3$--sphere $\Sigma$ such that
$M(K)$ is irreducible and has only two essential surfaces up to
isotopy.  Then, with respect to a standard framing, one of these
surfaces has numerical boundary slope $0$ and the other has numerical
boundary slope $r \not=0$.  If $r\not=\infty$ and if the
essential surface with boundary slope $0$ has genus $g\ge 2$ then
$$\frac{g - 1}{4\log_2(2g-2)} \le r^2 .$$
\EndCorollary

\Proof
Any knot in a homology 3--sphere has an essential spanning surface with
connected boundary and numerical boundary slope $0$.  Hence we may denote
the two non-isotopic essential surfaces in $M(K)$ by $F_1$ and $F_2$
where $F_2$ is a spanning surface.

The hypotheses imply in particular that $M(K)$ has at most two strict
essential surfaces, so one of the conclusions of Theorem
\ref{AtMostTwoSES} must hold.  Since $\Sigma$ is a homology sphere, it
cannot contain a non-separating torus or a Klein bottle.  Moreover
$M(K)$ is not a solid torus since it has two non-isotopic essential
surfaces.  This rules out conclusions (1) and (3a) of Theorem
\ref{AtMostTwoSES}.  Therefore $M(K)$ is either Seifert-fibered over
the disk with two singular fibers or a non-exceptional two-surface
knot manifold.  In the first case $\Sigma$ is homeomorphic to $S^3$
and $K$ is a torus knot.  If $K$ is the $(m,n)$--torus knot then $M(K)$
contains an essential annulus of slope $mn\not=0$ and an essential
spanning surface of genus $(m-1)(n-1)/2$.  These must be isotopic to
$F_1$ and $F_2$ respectively.  Thus the conclusions hold in this case.

If conclusion (3b) of Theorem \ref{AtMostTwoSES} holds then $F_1$ and
$F_2$ must be the two non-isotopic strict essential surfaces.  If we
set $r = p/q$, where $p$ and $q$ are relatively prime, then in the
notation of Theorem \ref{easyboundconsequence} we have $\Delta = p$,
$q_1 = q$, $g_2 = g$ and $m_2 = 1$.  In particular $r \not= 0$.
Furthermore if $r\not=\infty$ then $s_2\not=\emm$, and the inequality
in the statement of the corollary is equivalent to the inequality in
the conclusion of Theorem \ref{easyboundconsequence}.
\EndProof

For future reference (see \ref{qualitativediscussion}), we record the
following qualitative consequence of Theorem
\ref{easyboundconsequence}.

\Corollary\label{easyqualitative}
There is a positive-valued function $f_0(x)$ of a positive real
variable $x$ with the following properties.
\Parts
  \Part{(1)} For every $\epsilon>0$ we have
  $$\lim_{x\to\infty}x^{1-\epsilon}f_0(x)=0.$$

  \Part{(2)} If $K$ is any non-exceptional two-surface knot in a
  closed, orientable $3$--manifold $\Sigma$ such that $\pi_1(\Sigma)$
  is cyclic, and if $\emm$, $F_i$, $g_i$, $s_i$, $m_i$, $q_i$ and
  $\Delta$ are defined as in the statement of Theorem
  \ref{easyboundconsequence}, and if $g_2\ge2$, then
  $$
  \left(\frac{q_1}{\Delta}\right)^2 \le {m_2^2} f_0(g_2).
  $$
\EndParts
\EndCorollary

\section{Short subgraphs II}
\label{secondbigirthsection}

This section presents the more subtle combinatorial ideas that are
needed in the proof of Theorem \ref{hardboundconsequence}, our second
main concrete result about two-surface knots in manifolds with cyclic
fundamental group.

\Definition\label{systemdef}
By an {\it essential arc system} in a compact, orientable surface $S$
we mean a non-empty, properly embedded $1$--manifold ${\mathcal
A}\subset S$ such that every component of $\mathcal A$ is a
non-boundary-parallel arc in $S$.  If $\mathcal A$ is an essential arc
system we shall denote by $\partial_{\mathcal A}S$ the union of all
boundary components of $S$ which meet $\mathcal A$.  We shall denote
by $\Gamma_{\mathcal A}$ the trivalent graph with $|\Gamma_{\mathcal
A}|={\mathcal A}\cup\partial_{\mathcal A} S$, in which the vertex set
is $\partial {\mathcal A}$. The edges of $\Gamma_{\mathcal A}$ are the
components of $\inter{\mathcal A}$, which we call {\it interior
edges}, and the components of $(\partial_{\mathcal
A}S)-(\partial{\mathcal A})$, which we call {\it boundary edges}\/.
Note that every vertex of $\Gamma_A$ is an endpoint of a unique
interior edge.

In this section we shall often consider subgraphs of $\Gamma_{\mathcal
A}$. Such a subgraph need not be connected, and it may have vertices
of valence $0$, $1$, $2$ and $3$. If $\Gamma$ is a subgraph of
$\Gamma_{\mathcal A}$ we define the {\it interior edges} and {\it
boundary edges} of $\Gamma$ to be, respectively, the interior edges
and boundary edges of $\Gamma_{\mathcal A}$ that are contained in
$|\Gamma|$.
\EndDefinition

By a {\it reduced arc system} in a compact, orientable surface $S$ we
mean an essential arc system in $S$ such that no two components of
$\mathcal A$ are parallel.

\Proposition\label{pioneinjectivesubgraph}
Suppose that $S$ is a compact, connected, orientable surface which is
not an annulus, that $\mathcal A$ is a reduced arc system in $S$ and
that $\Gamma_0$ is a subgraph of $\Gamma_{\mathcal A}$. Let $\nu$
denote the number of interior edges of $\Gamma_0$. Then $\Gamma_0$ has
a subgraph $\Gamma_1$ such that
\Parts
  \Part{(1)} $|\Gamma_1|$ is $\pi_1$--injective in $S$,
  \Part{(2)} $|\Gamma_1|$ contains every vertex of $\Gamma_0$ and every
  boundary edge of $\Gamma_0$, and
  \Part{(3)} the number of interior edges of $\Gamma_1$ is at least $\nu/3$.
\EndParts
\EndProposition

\Proof
We may assume $\nu>0$. Let ${\mathcal A}_0$ denote the union
of those components of $\mathcal A$ that are contained in $\Gamma_0$, so
that $\nu$ is the number of components of ${\mathcal A}_0$. In particular
${\mathcal A}_0\ne\emptyset$. We let $\mu$ denote the number of components
of $(\inter S)-(\inter{\mathcal A}_0)$ that are topological open disks. We
begin by showing that $\mu\le2\nu/3$. If $N$ is a regular neighborhood of
${\mathcal A}_0$ in $S$, if $D_1,\ldots,D_\mu$ are the components of
$\overline{S-N}$ that are closed disks, and if $c_i$ is the number of
frontier components of $D_i$ then, since each frontier component of $D_1\cup\cdots\cup
D_\mu$ is a component of the frontier of $N$, we have $\sum_{i=1}^\mu
c_i\le2\nu$. We have $c_i\ne0$ for every $i$ because $S$ is connected
and ${\mathcal A}_0\ne\emptyset$. We have $c_i\ne1$ for every $i$ because
no component of ${\mathcal A}_0\subset{\mathcal A}$ is a boundary-parallel
arc. We have $c_i\ne2$ for every $i$ because $S$ is not an annulus and
no two components of ${\mathcal A}_0\subset{\mathcal A}$ are parallel
arcs. Hence $c_i\ge3$ for every $i$, and so $3\mu\le\sum_{i=1}^\mu
c_i\le2\nu$. This proves that $\mu\le2\nu/3$.

To prove the proposition, it suffices to find a submanifold ${\mathcal
A}_1$  which is a union of components of ${\mathcal A}_0$, has at least
$\nu/3$ components, and has the property that no component of $(\inter
S)-(\inter {\mathcal A}_1)$ is an open disk. Indeed, if ${\mathcal A}_1$ has
these properties, and if we define $\Gamma_1$ to be  the subgraph of
$\Gamma_0$ made up of all vertices and boundary edges of $\Gamma_0$
and of those interior edges which are components of $\inter{\mathcal
A}_1$, then $\Gamma_1$ clearly satisfies conditions (1)--(3) of the
proposition.

For any properly embedded $1$--manifold ${\mathcal B}\subset S$ let
$\nu_{\mathcal B}$ and $\mu_{\mathcal B}$ denote respectively the number of
components of $\mathcal B$ and the number of components of $(\inter
S)-(\inter{\mathcal B})$ that are topological open disks. We shall say
that $\mathcal B$ is {\it admissible} if it a union of components of
${\mathcal A}_0$ and if $\nu_{\mathcal B}-\mu_{\mathcal B}=\nu-\mu$. Clearly
${\mathcal A}_0$ is itself admissible. If $\mathcal B$ is any admissible
submanifold, then
$$
\nu_{\mathcal B}\ge\nu_{\mathcal B}-\mu_{\mathcal B}=
\nu-\mu\ge\nu-\frac{2\nu}{3}=\nu/3 ,
$$
ie, $\mathcal B$ has at least $\nu/3$ components.

Among all admissible submanifolds choose one, ${\mathcal A}_1$, for
which the number of components is as small as possible. We shall
complete the proof by showing that no component of $(\inter S)-(\inter
{\mathcal A}_1)$ is an open disk.

Suppose that some component $D$ of $(\inter S)-(\inter {\mathcal
A}_1)$ is an open disk. Choose any component $A$ of the frontier of
$D$ in $S$. Then $A$ is also a component of ${\mathcal A}_1$. Set
${\mathcal B}={\mathcal A}_1-A$. If $D$ is the only component of
$(\inter S)-(\inter {\mathcal A}_1)$ whose closure contains $A$, then
$D\cup A$ is a component of $(\inter S)-(\inter {\mathcal B})$
homeomorphic to an open annulus, and the other components of $(\inter
S)-(\inter {\mathcal B})$ are the components $\ne D$ of $(\inter
S)-(\inter {\mathcal A}_1)$. If $C$ is a second component of $(\inter
S)-(\inter {\mathcal A}_1)$ whose closure contains $A$, then $D\cup
A\cup C$ is a component of $(\inter S)-(\inter {\mathcal B})$
homeomorphic to $C$, and the other components of $(\inter S)-(\inter
{\mathcal B})$ are the components $\ne C,D$ of $(\inter S)-(\inter
{\mathcal A}_1)$. In either case it follows that $\mu_{\mathcal
B}=\mu_{{\mathcal A}_1}-1$. Since $\nu_{\mathcal B}=\nu_{{\mathcal
A}_1}-1$, we have $\nu_{\mathcal B}-\mu_{\mathcal B}=\nu_{{\mathcal
A}_1}-\mu_{{\mathcal A}_1}=\nu-\mu$, so that ${\mathcal B}$ is
admissible. This contradicts the minimality of ${\mathcal A}_1$, and
the proof is complete.
\EndProof

\Definition\label{labelingdef}
Let $S$ be a compact, orientable surface. By a {\it labeling} of an
essential arc system $\mathcal A$ we mean a surjective map $\iota$
from the set of interior edges of ${\mathcal A}$ to some finite set
$I$, called the {\it label set} of $\iota$. If $e$ is an interior edge
of $\Gamma_{\mathcal A}$, we shall refer to $\iota(e)$ as the {\it
label} of $e$. For every vertex $v$ of $\Gamma_{\mathcal A}$ we define
the {\it label} of $v$, denoted $\iota(v)$, to be the label of the
unique interior edge having $v$ as an endpoint. For every $i\in I$ we
define the {\it multiplicity of $i$} with respect to the labeling of
$\iota$, denoted $\theta^\iota_i$, to be the number of arcs in
$\mathcal A$ with label $i$. Thus $\sum_{i\in
I}\theta^\iota_i=\#({\mathcal A})$. For every vertex $v$ and every
edge $e$ we shall write $\theta^\iota(v)=\theta^\iota(\iota(v))$ and
$\theta^\iota(e)=\theta^\iota(\iota(e))$.
 
If $\iota$ is a labeling of an essential arc system $\mathcal A$, then
for every subgraph $\Gamma$ of $\Gamma_{\mathcal A}$ we define
$$
\theta^{\iota}(\Gamma)=\min_e\theta_e^{\mathcal A}
$$
where $e$ ranges over the interior edges of $\Gamma$.

Let $\iota$ be a labelling of an essential arc system $\mathcal A$, with
label set $I$. We define a {\it system of weights} for $\mathcal A$ with
respect to $\iota$ to be a positive real-valued function $\lambda$ on
the set $I$.  If $\lambda$ is a system of weights, and $\Gamma$ is a
subgraph of $\Gamma_{\mathcal A}$, we shall denote by $\lambda(\Gamma)$
the quantity $\sum_{v\in{\mathcal V}}\lambda(\iota(v))$, where $\mathcal V$
denotes the set of vertices of $\Gamma$.
\EndDefinition

\Definition\label{phidef}
For every real number $\tau>1$, we define a positive real-valued
function $\phi_\tau$ on the set of all positive integers by
$$\phi_\tau(n)=\frac{1}{\tau-1}\min_{m}\tau^{2m+2}n^{1/m},$$
where $m$ ranges over all positive integers. Note that $\phi_\tau$ is
monotone increasing, and that for every
$\epsilon>0$ we have
$$
\lim_{n\to\infty}\frac{\phi_\tau(n)}{n^\epsilon}=0 .
$$
\EndDefinition

\Remark
One can show using calculus that for a given $\tau>1$ the function
$\phi_\tau$ grows roughly like $e^{c\sqrt{\ln x}}$ for some constant
$c$.  We will not give a precise statement or proof of this here, but
we will prove and use a more technical result along these lines, Lemma
\ref{precalculus}.
\EndRemark

\Lemma\label{keyinequality}
Suppose that $S$ is a compact, connected, orientable surface which is
not an annulus, that $\mathcal A$ is a reduced arc system in $S$, and
that $\iota$ is a labeling for $\mathcal A$ with label set $I$. Set
$\theta_i=\theta_i^\iota$ for every $i\in I$. Set $\Theta=\#({\mathcal
A})=\sum_{i\in I}\theta^\iota_i$ and $\theta_\infty=\max_{i\in
I}\theta_i$. Let a real number $q>1$ be given, and set
$$
\tau=\frac{7q-1}{q-1} .
$$
Suppose that $E^*$ is a set of interior edges of $\Gamma_{\mathcal
A}$.  For every $i\in I$ set $E_i^*=E^*\cap E_i$ and
$\theta_i^*=\#(E_i^*)$, and suppose that $\theta_i^*\le\theta_i/q$ for
every $i\in I$.  Suppose that $\lambda$ is a system of weights for
$\mathcal A$ with respect to $\iota$. Then there is a subgraph
$\Gamma_1$ of $\Gamma_{\mathcal A}$ such that 
\Parts
  \Part{(1)} $|\Gamma_1|$ is $\pi_1$--injective in $S$,
  \Part{(2)} $|\Gamma_1|$ contains no edge in $E^*$,
  \Part{(3)} the Euler characteristic $\chi(|\Gamma_1|)$ is strictly
  negative, and
  \Part{(4)}$\phantom{99}$\vspace{-15pt}
  $$
  \frac{\lambda(\Gamma_1)}{\theta^\iota(\Gamma_1)\cdot|\chi(|\Gamma_1|)|}<
  \phi_\tau(\theta_\infty)\frac{\sum_{i\in I}\lambda(i)}{\Theta}.
  $$
\EndParts
\EndLemma

\Proof
For every subset $J$ of $I$, let us denote by $E_J$ the set of all
interior edges $e$ of $\Gamma_{\mathcal A}$ such that $\iota(e)\in
J$. We set $\Theta_J=\#(E_J)=\sum_{i\in J}\theta_i$.  We also set
$E_J^*=E_J\cap E^*$ and $\Theta_J^*=\#(E_J)=\sum_{i\in J}\theta_i^*$.

For every  $J\subset I$, let us denote by $V_J$ the set of all
vertices $v$ of $\Gamma_{\mathcal A}$ such that $\iota(v)\in J$.
If a vertex $v$ is an endpoint of an interior edge $e$, the
definition of $\iota(v)$ implies that $\iota(v)=\iota(e)$. Since each
vertex of $\Gamma_{\mathcal A}$ is an endpoint of a unique interior edge,
and each edge has two endpoints, it follows that
\Equation\label{keyone}
\#(V_J)=2\Theta_J
\EndEquation
for every $J\subset I$.

According to the definition of $\phi_\tau$ (\ref{phidef}), we may fix
a positive integer $m$ such that
$$
\phi_\tau(\theta_\infty)=\frac{\tau^{2m+2}}{\tau-1}\theta_\infty^{1/m} .
$$
We set $A=\tau^{m+1}/(\tau-1)$, $\alpha=A^{-1}$, and $\omega=(\sum_{i\in
I}\lambda(i))/\Theta$.

By hypothesis we have $\theta_i^*\le\theta/q$ for every $i\in
I$. Summing this over the labels in any given subset $J$ of $I$, we
find
\Equation\label{keytwo}
\Theta_J^*\le{\Theta_J/q} .
\EndEquation
We decompose the label set $I$ as a disjoint union
$$
I=I_0\sqcup\cdots\sqcup I_m ,
$$
where
$$\begin{aligned}
I_0&=\{i\in I:\lambda(i)>A\omega\theta_i\} ,\\
I_j&=\{i\in I-I_0:\theta_\infty^{(m-j)/m}<\theta_i\le\theta_{\infty}^{(m-j+1)/m}\}
\hbox{ for } j=1,\ldots,m-1, \hbox{ and }\\
I_m&=\{i\in I-I_0:1\le\theta_i\le\theta_{\infty}^{1/m}\} .
\end{aligned}$$
We set $\Theta_j=\Theta_{I_j}$ and $\Theta_j^*=\Theta_{I_j}^*$ for
$j=0,\ldots,m$. Thus we have $\Theta_j=\sum_{i\in I_j}\theta_i$ and
$\Theta=\sum_{j=0}^m\Theta_j$; similarly, $\Theta_j^*=\sum_{i\in
I_j}\theta_i^*$ and $\Theta^*=\sum_{j=0}^m\Theta_j^*$.

Using the definitions of $I_0$ and $\omega$ we find that
$$\begin{aligned}
A\omega\Theta_0&=\sum_{i\in I_0}A\omega\theta_i\le\sum_{i\in I_0}\lambda(i)\\
&\le\sum_{i\in I}\lambda(i)=\omega\Theta ,
\end{aligned}$$
so that
\Equation\label{keythree}
\Theta_0\le\alpha\Theta .
\EndEquation
On the other hand, since
$$
\sum_{j=0}^m\tau^j=\frac{\tau^{m+1}-1}{\tau-1}<A ,
$$
$$
\sum_{j=0}^m\Theta_j=\Theta=A\alpha\Theta>\sum_{j=0}^m\tau^j\alpha\Theta ,
\leqno{\hbox{we have}}$$
and hence $\Theta_j>\tau ^j\alpha\Theta$ for some
$j\in\{0,\ldots,m\}$. Let $k$ denote the smallest index for which
$\Theta_k>\tau ^k\alpha\Theta$. It follows from (\ref{keythree})
that $k>0$.

We set $I^+=I_k\cup\cdots\cup I_m$ and $I^-=I_0\cup\cdots\cup
I_{k-1}$, so that $I=I^+\sqcup I^-$.  We define a subgraph $\Gamma_0$
of $\Gamma_{\mathcal A}$ to consist of all vertices in $V_{I^+}$, all
interior edges in $E_{I_k}-E_{I_k}^*$, and all boundary edges whose
endpoints both lie in $V_{I^+}$.  Since an interior edge in
$E_{I_k}-E_{I_k}^*$ has its endpoints in $V_{I_k}*\subset V_{I^+}$, it
follows that $\Gamma_0$ is indeed a subgraph of $\Gamma_{\mathcal A}$.

The number of interior edges of $\Gamma_0$ is $\Theta_k-\Theta_k^*$.
It therefore follows from Proposition \ref{pioneinjectivesubgraph}
that there is a subgraph $\Gamma_1$ of $\Gamma_0$ such that
$|\Gamma_1|$ is $\pi_1$--injective in $S$, contains every vertex of
$\Gamma_0$ and every boundary edge of $\Gamma_0$, and contains at
least $(\Theta_k-\Theta_k^*)/3$ interior edges. In particular
$\Gamma_1$ satisfies condition (1) of Lemma \ref{keyinequality}. Since
the edge set of $\Gamma_0\supset\Gamma_1$ is $E_{I_k}-E_{I_k}^*$, the
subgraph $\Gamma_1$ satisfies condition (2) as well. We shall complete
the proof by showing that it also satisfies conditions (3) and (4).

For any interior edge $e$ of $\Gamma_1$ we have $e\in E_{I_k}$ and
hence $\iota(e)\in I_k$; since we have shown that $k>0$, the
definition of the $I_j$ then implies that
$\theta_e=\theta_{\iota(e)}\ge\theta_\infty^{(m-k)/m}$. Hence
\Equation\label{keyfour}
\theta^{\iota}(\Gamma_1)\ge\theta_\infty^{(m-k)/m} .
\EndEquation
We next turn to the estimation of $\lambda(\Gamma_1)$. Since
$\mathcal V_{I^+}$ is the vertex set of $\Gamma_1$, we have
$$\begin{aligned}
\lambda(\Gamma_1)&=\sum_{v\in{\mathcal V}_{I^+}}\lambda(\iota(v))\\
&=\sum_{i\in I^+}\lambda(i)\cdot\#(V_i)=2\sum_{i\in I^+}\lambda(i)\theta_i,
\end{aligned}$$
where the last step follows by applying (\ref{keyone}) with $J=I^+$.
Since we have shown that $k>0$, we have $I^+\cap I_0=\emptyset$.
In view of the definition of $I_0$ it follows that
$\lambda(i)\le A\omega\theta_i$ for every $i\in I^+$. Hence
$$
\lambda(\Gamma_1) \le 2A\omega\sum_{i\in I^+}\theta_i^2 .
$$
On the other hand, the definition of the $I_j$ shows that for every
$i\in I^+$ we have $\theta_i\le\theta_\infty^{(m-k+1)/m}$. Hence
$$
\lambda(\Gamma_1)
\le2A\omega\theta_\infty^{(m-k+1)/m}\sum_{i\in I^+}\theta_i
\le2A\omega\theta_\infty^{(m-k+1)/m}\sum_{i\in I}\theta_i ,
$$
ie,
\Equation\label{keyfive}
\lambda(\Gamma_1) \le2A\omega\theta_\infty^{(m-k+1)/m}\Theta.
\EndEquation
Now we turn to the estimation of $\chi(|\Gamma_1|)$. First note that
\Equation\label{keysix}
\chi(|\Gamma_1|)=\beta-\gamma ,
\EndEquation
where $\beta$ denotes the number of simply connected components of the
set ${\mathcal B}=|\Gamma_1|\cap\partial S=|\Gamma_0|\cap\partial S$
and $\gamma$ denotes the number of interior edges of
$\Gamma_1$. According to the construction of $\Gamma_1$ we have
$\gamma\ge(\Theta_k-\Theta_k^*)/3$. By taking $J=I_k$ in
(\ref{keytwo}) we find that $\Theta_k^*\le\Theta_k/q$, and so
$\gamma\ge(q-1)\Theta_k/3q$.  But by our choice of $k$ we have
$\Theta_k>\tau ^k\alpha\Theta$. Hence
\Equation\label{keyseven}
\gamma>\frac{(q-1)\tau ^k}{3q}\alpha\Theta .
\EndEquation
To estimate $\beta$, we let $\beta'$ denote the number of components
of the set ${\mathcal B}'=(\partial S)-{\mathcal B}$, and note that
$\beta\le\beta'$ (with equality holding unless some component of
$\partial S$ is contained in ${\mathcal B}'$). According to our
construction of $\Gamma_0$, the set $\mathcal B$ consists of the
vertices in the set $V_{I^+}$ and the boundary edges of
$\Gamma_{\mathcal A}$ that have both endpoints in $V_{I^+}$. Thus
every component of ${\mathcal B}'$ contains at least one vertex of
$V_{I^-}$, and therefore
$$
\beta\le\beta'\le\#(V_{I^-})=2\Theta_{I^-} ,
$$
where the last step follows by applying (\ref{keyone}) with $J=I^+$.
But since we defined $k$ to be the smallest index for which
$\Theta_k>\tau ^k\alpha\Theta$, we have $\Theta_j\le\tau
^j\alpha\Theta$ for $j=0,\ldots,k-1$. Hence
$$
\Theta_{I^-}=\sum_{j=0}^{k-1}\Theta_j\le\sum_{j=0}^{k-1}\tau
^j\alpha\Theta=
\frac{\tau ^k-1}{\tau-1}\alpha\Theta ,
$$
$$
\beta\le\frac{2(\tau ^k-1)}{\tau-1}\alpha\Theta .
\leqno{\hbox{and so}}
$$
Combining this with (\ref{keysix}) and (\ref{keyseven}) we find that
$$
\chi(|\Gamma_1|)=\beta-\gamma<
\frac{2(\tau ^k-1)}{\tau-1}\alpha\Theta-
\frac{(q-1)\tau^k}{3q}\alpha\Theta .
$$
Since $\tau=(7q-1)/( q-1)$, this last inequality simplifies to
\Equation\label{keyeight}
\chi(|\Gamma_1|)<-\frac{2\alpha\Theta}{\tau-1} .
\EndEquation
It follows from (\ref{keyeight}) that $\Gamma_1$ satisfies condition
(3) of Lemma \ref{keyinequality}. Furthermore, it follows from
(\ref{keyfour}), (\ref{keyfive}) and (\ref{keyeight}) that
$$\begin{aligned}
\frac{\lambda(\Gamma_1)}{\theta^{\mathcal A}(\Gamma_1)\cdot|\chi(|\Gamma_1|)|}
&<2A\frac{\omega\theta_\infty^{(m-k+1)/m}\Theta}
{\theta_\infty^{(m-k)/m}(2\alpha\Theta/(\tau-1))}\\
&=(\tau-1)A^2\omega\theta_\infty^{1/m}\\
&=\frac{\tau ^{2m+2}}{\tau-1}\theta_\infty^{1/m}\omega
=\phi_\tau(\theta_\infty)\omega ,
\end{aligned}$$
which in view of the definition of $\omega$ gives condition (4) of
Lemma \ref{keyinequality}. 
\EndProof

\Remark\label{simplyconnectedremark}
If a graph $\Gamma_1$ satisfies the conclusions of Lemma
\ref{keyinequality}, and if $\Gamma_1'$ denotes the subgraph of
$\Gamma_1$ such that $|\Gamma_1'|$ is union of all non-simply-connected
components of $|\Gamma_1|$, then conclusions (1)--(4) also
hold when $\Gamma_1$ replaced by $\Gamma_1'$.  (To verify condition
(4) observe that $\lambda(\Gamma_1')\le \lambda(\Gamma_1)$,
$\chi(|\Gamma_1'|)\le\chi(|\Gamma_1|)$ and
$\theta^\iota(\Gamma_1')\ge\theta^\iota(\Gamma_1)$.)  Thus Lemma
\ref{keyinequality} remains true if we add the condition
\Parts
{\sl
  \Part{(5)} No component of $|\Gamma_1|$ is simply-connected.
}
\EndParts
\EndRemark

\Definition\label{reductiondef}
If ${\mathcal A}^0$ is any essential arc system in a compact, orientable
surface $S$, we may define an equivalence relation on the set of
components of ${\mathcal A}^0$ in which two components are equivalent if
and only if they are parallel in $S$.  The equivalence classes for
this relation will be called {\it ${\mathcal A}^0$--parallelism classes.}
An essential arc system $\mathcal A$ will be called a {\it reduction} of
${\mathcal A}^0$ if ${\mathcal A}\subset {\mathcal A}^0$ and if every ${\mathcal
A}^0$--parallelism class contains exactly one component of $\mathcal A$. It
is clear that a reduction always exists, that it is unique up to
isotopy, and that it is itself a reduced arc system.

If $\mathcal A$ is a reduction of an essential arc system ${\mathcal A}^0$, we
define the {\it ${\mathcal A}^0$--width} of any interior edge $e$ of
$\Gamma_{\mathcal A}$ to be the cardinality of the ${\mathcal
A}^0$--parallelism class containing the component $\bar e$ of ${\mathcal
A}^0$.

If $\mathcal A$ is a reduction of an essential arc system ${\mathcal A}^0$, we
have $|\Gamma_{\mathcal A}|\subset|\Gamma_{{\mathcal A}^0}|$. Indeed,
$\Gamma_{{\mathcal A}^0}$ has a subgraph $\Gamma'_{\mathcal A}$ which is a
subdivision of a subgraph of $\Gamma_{\mathcal A}$, in the sense that
every vertex of $\Gamma_{{\mathcal A}}$ is a vertex of $\Gamma'_{\mathcal A}$,
and every edge of $\Gamma_{{\mathcal A}}$ is a union of edges and vertices
of $\Gamma'_{\mathcal A}$.  It follows that for every subgraph $\Gamma $
of $\Gamma_{{\mathcal A}}$ there is a unique subgraph $\Gamma ^0$ of
$\Gamma_{{\mathcal A}^0}$ such that $|\Gamma ^0|=|\Gamma |$. We shall
refer to $\Gamma ^0$ as the subgraph of $\Gamma_{{\mathcal A}^0}$ {\it
associated to} $\Gamma $.
\EndDefinition

\Lemma\label{standardweights}
Suppose that $S$ is a compact, connected, orientable surface which is
not an annulus, that ${\mathcal A}^0$ is an essential arc system in
$S$, and that ${\mathcal A}$ is a reduction of ${\mathcal
A}^0$. Suppose that $\iota$ is a labeling for ${\mathcal A}$ with
label set $I$.  For every interior edge $e$ of $\Gamma_{\mathcal A}$,
let $w(e)$ denote the ${\mathcal A}^0$--width of $e$. Define a weight
system for $\mathcal A$ with respect to $\iota$ as follows: for each
label $i\in I$, set $\lambda(i)=\max_e w(e)$, where $e$ ranges over
all interior edges of $\Gamma_{\mathcal A}$ with label $i$. Suppose
that $\Gamma$ is a subgraph of $\Gamma_{\mathcal A}$, and let $\Gamma
^0$ denote the subgraph of $\Gamma_{{\mathcal A}^0}$ associated to
$\Gamma $. Then we have
$$
\length(\Gamma ^0)\le\frac{3}{2}\lambda(\Gamma) .
$$
\EndLemma

\proof
For each vertex $v$ of $\Gamma_{\mathcal A}$ we set $w(v)=w(e)$, where
$e$ is the unique edge of $\Gamma_{\mathcal A}$ having $v$ as an
endpoint. The definition of the $\lambda(i)$ implies that
$w(v)\le\lambda(\iota(v))$ for every vertex $v$.

Since $\mathcal A$ is a reduction of ${\mathcal A}^0$, each arc in
${\mathcal A}^0$ is parallel to $\bar e$ for a unique interior edge
$e$ of $\Gamma_{\mathcal A}$. Thus if for each interior edge $e$ of
$\Gamma_{\mathcal A}$ we let ${\mathcal F}_e$ denote the union of all
arcs in ${\mathcal A}^0$ that are parallel to $\bar e$, then
${\mathcal A}^0$ is the disjoint union of the sets ${\mathcal F}_e $
as $e$ ranges over the interior edges of $\Gamma_{\mathcal A}$. By
definition the number of components of ${\mathcal F}_e$ is the width
$w(e)$. For each $e$, since ${\mathcal F}_e$ is a family of parallel
arcs and $S$ is not an annulus, there is a topological disk or arc
$R_e\subset S$ such that ${\mathcal F}_e\subset R_e$ and $\partial
R_e\subset{\mathcal F_e}\cup\partial S$.  Furthermore,
$R_e\cap\partial S$ consists of two possibly degenerate arcs. If $e$
and $e'$ are distinct interior edges we have $R_e\cap
R_{e'}=\emptyset$. We set ${\mathcal R}=\bigcup_e R_e$, where $e$
ranges over the interior edges of $\Gamma_{\mathcal A}$. We have
${\mathcal A}^0\subset{\mathcal R}$.

If $e$ is an interior edge of $\Gamma_{\mathcal A}$, each of the two arcs
that make up $R_e\cap\partial S$ contains exactly one endpoint of $e$.
Hence if we set ${\mathcal B}={\mathcal R}\cap\partial S$, each component of
$\mathcal B$ contains a unique vertex of $\Gamma_{\mathcal A}$. We shall
denote by $B_v$ the component of $\mathcal B$ containing a given vertex
$v$ of $\Gamma_{\mathcal A}$. Note that since ${\mathcal A}^0\subset{\mathcal R}$,
every vertex of $\Gamma_{{\mathcal A}^0}$ lies in $\mathcal R$. Note also that
since ${\mathcal F}_e$ consists of $w(e)$ arcs for each interior edge $e$
of $\Gamma_{\mathcal A}$, it follows that for each vertex $v$ of
$\Gamma_{\mathcal A}$ the arc $B_v$ contains exactly $w(v)$ vertices of
$\Gamma_{{\mathcal A}^0}$, and hence contains exactly $w(v)-1$ edges of
$\Gamma_{{\mathcal A}^0}$.

Let us fix an orientation of each component of $\partial S$. For every
vertex $v$ of $\Gamma_{\mathcal A}$ we may write $B_v$ in a unique way as
a union of arcs $B_v^+$ and $B_v^{-}$, one or both of which may be
degenerate, such that and $v$ is the negative endpoint of $B_v^+$, and
the positive endpoint of $B_v^-$, with respect to the orientation of
the component of $\partial S$ containing $g$. It follows that
$B_v^+\cap B_v^{-}=\{v\}$. In particular, if $b_v^+$ and $b_v^-$
denote the number of edges of $\Gamma_{{\mathcal A}^0}$ contained in
$B_v^+$ and $B_v^-$ respectively, we have
$$b_v^++b_v^-=w(v)-1.$$
For every boundary edge $g$ of $\Gamma_{\mathcal A}$, let us denote by
$v_+(g)$ and $v_-(g)$ the positive and negative endpoints of $g$ with
respect to the orientation of the component of $\partial S$ containing
$g$. (It may happen that $v_+=v_-$.) We have $\bar g\cap{\mathcal
B}=B_{v^+(g)}^-\cup B_{v^-(g)}^+$. In particular, $B_{v^+(g)}^-\cup
B_{v^-(g)}^+$ contains all the vertices of $\Gamma_{{\mathcal A}^0}$
in $\bar g$. Since each component of $\mathcal B$ is a (possibly
degenerate) topological arc containing exactly one vertex of
$\Gamma_{\mathcal A}$, the arcs $B_{v^+(g)}^-$ and $B_{v^-(g)}^+$ are
distinct, and hence $g-(B_{v^+(g)}^-\cup B_{v^-(g)}^+)$ is an open
topological arc; as it contains no vertices of $\Gamma_{{\mathcal
A}^0}$, it must be an edge of $\Gamma_{{\mathcal A}^0}$. Hence the
number of edges of $\Gamma_{{\mathcal A}^0}$ contained in $g$ is
$b_{v^+(g)}^-+ b_{v^-(g)}^++1$.

It now follows that if $G$ and $G^0$ respectively denote the sets of
boundary edges of the subgraph $\Gamma$ of $\Gamma_{\mathcal A}$ and
of the associated subgraph $\Gamma^0$ of $\Gamma_{{\mathcal A}^0}$,
then
$$
\#(G^0)\le\sum_{g\in G}(b_{v^+(g)}^-+ b_{v^-(g)}^++1)=
\#(G)+\sum_{g\in G}b_{v^+(g)}^-+\sum_{g\in G} b_{v^-(g)}^+ .
$$
Note that $v^+(g)$ and $v^-(g)$ are vertices of $\Gamma$ for every
$g\in G$, and that for a given vertex $v$ of $\Gamma$ there is at most
one edge $g$ such that $v=v^-(g)$, and at most one edge $g'$ such that
$v=v^-(g')$. Hence if $V$ denotes the vertex set of $\Gamma$, we have
$\sum_{g\in G}b_{v^+(g)}^-\le\sum_{v\in V}b_v^-$ and $\sum_{g\in G}
b_{v^-(g)}^+\le\sum_{v\in V}b_v^+$, and therefore
$$\begin{aligned}
\#(G^0)&\le\#(G)+\sum_{v\in V}b_v^-+\sum_{v\in V}b_v^+
=\#(G)+\sum_{v\in V}(b_v^-+b_v^+)\\
&=\#(G)+\sum_{v\in V}(w(v)-1)=(\#(G)-\#(V))+\sum_{v\in V}w(v) .
\end{aligned}$$
But we have observed that $w(v)\le\lambda(\iota(v))$, and since
$\Gamma\cap\partial S$ is a subgraph of a triangulated $1$--manifold we
have $\#(G)\le\#(V)$. Hence
$$
\#(G^0)\le\sum_{v\in V}\lambda(\iota(v))=\lambda(\Gamma) .
$$
Finally, if $E$ denotes the set of interior edges of $\Gamma$, then
since no two edges in $E$ can have a common endpoint, we have
$$
\#(E)\le\frac{1}{2}\#(V)\le\frac{1}{2}\sum_{v\in V}\lambda(\iota(v))
=\frac{1}{2}\lambda(\Gamma) ,
$$
so that
$$
\length(\Gamma^0)=\#(G^0)+\#(E)\le\frac{3}{2}\lambda(\Gamma) .\eqno{\qed}
$$

\Proposition\label{keyconsequence}
Suppose that $S$ is a compact, connected, orientable surface which is
not an annulus, that ${\mathcal A}^0$ is an essential arc system in
$S$, and that ${\mathcal A}$ is a reduction of ${\mathcal A}^0$. Let a
real number $q>1$ be given, and set
$$
\tau=\frac{7q-1}{q-1} .
$$
Suppose that $\iota$ is a labeling for ${\mathcal A}$ with label set
$I$. Set $\theta_i=\theta_i^\iota$ for every $i\in I$. Set
$\Theta=\#({\mathcal A'})=\sum_{i\in I}\theta^\iota_i$ and
$\theta_\infty=\max_{i\in I}\theta_i$.  Suppose that $E^*$ is a set of
interior edges of $\Gamma_{\mathcal A}$. For every $i\in I$ set
$E_i^*=E^*\cap E_i$ and $\theta_i^*=\#(E_i^*)$, and suppose that
$\theta_i^*\le\theta_i/q$ for every $i\in I$. For every interior edge
$e$ of $\Gamma_{\mathcal A}$, let $w(e)$ denote the ${\mathcal
A}^0$--width of $e$. For each label $i\in I$, set $\lambda(i)=\max_e
w(e)$, where $e$ ranges over all interior edges of $\Gamma_{\mathcal
A}$ with label $i$. Then $\Gamma_{\mathcal A}$ has a subgraph $K$ such
that 
\Parts
  \Part{(1)} $|K|$ is $\pi_1$--injective in $S$, 
  \Part{(2)} $|K|$ contains no edge in $E^*$,
  \Part{(3)} the first Betti number of $|K|$ is equal to $2$,
  \Part{(4)} $K$ has no vertices of valence $\le 1$, and
  \Part{(5)} if $K^0$ is the subgraph of $\Gamma_{{\mathcal A}^0}$
   associated to $K$, we have
   $$
   \frac{\length(K^0)}{\theta^\iota(K)}<
   6\phi_\tau(\theta_\infty)(\log_2(2\Theta))\frac{\sum_{i\in I}\lambda_i}{\Theta} .
   $$
\EndParts
\EndProposition

\Proof
If for every $i\in I$ we define $\lambda(i)$ as in the statement of
Proposition \ref{keyconsequence}, then $\lambda$ is a weight system
for $\mathcal A$ with respect to $\iota$. Applying Proposition
\ref{keyinequality} with this choice of the weight system $\lambda$,
we fix a subgraph $\Gamma_1$ of $\Gamma_{\mathcal A}$ such that
conditions (1)--(5) of \ref{keyinequality} and
\ref{simplyconnectedremark} hold.  We denote by $\Gamma_1^0$ the
subgraph of $\Gamma_{{\mathcal A}^0}$ associated to $\Gamma_1$.

According to Lemma \ref{standardweights} we have $\length(\Gamma_1
^0)\le\frac{3}{2}\lambda(\Gamma_1)$. Furthermore, it is apparent that
$\chi(\Gamma_1^0)=\chi(\Gamma_1)$. Hence condition (4) of
\ref{keyinequality} implies that
\Equation\label{keyconsequenceone}
\frac{\length(\Gamma_1^0)}{\theta^\iota(\Gamma_1)\cdot|\chi(|\Gamma_1^0|)|}
<\frac{3}{2}\phi_\tau(\theta_\infty)\frac{\sum_{i\in I}\lambda(i)}{\Theta} .
\EndEquation
On the other hand, since $\chi(|\Gamma_1^0|)=\chi(|\Gamma_1|)<0$ by
condition (3) of \ref{keyinequality}, and since $\Gamma_1^0$ has no
simply connected components by condition (5) of
\ref{simplyconnectedremark}, it follows from Proposition
\ref{generalbigirthbound} that
$$
\bigirth(\Gamma_1^0)\le
4(\log_2|2\chi(\Gamma_1)|)\,\left\lfloor
\frac{\length(\Gamma_1^0)}{|\chi(|\Gamma_1|)|}
\right\rfloor .
$$
Hence $\Gamma_1^0$ has a subgraph $H$ such that $\chi(|H|)<0$ and
\Equation\label{keyconsequencetwo}
\length(H)\le 4(\log_2|2\chi(|\Gamma_1|)|)
\,\left\lfloor\frac{\length(\Gamma_1^0)}{|\chi(|\Gamma_1|)|}\right\rfloor .
\EndEquation
After possibly replacing $H$ by a smaller subgraph (which cannot
increase its length), we may assume that $H$ is connected, has Betti
number $2$ and has no valence--1 vertices. Since $\Gamma_1^0$ is a
subdivision of $\Gamma_1$, the absence of vertices of valence $\le1$
in $H$ implies that $|H|=|K|$ for some subgraph $K$ of $\Gamma_1$.
Thus $H$ is the subgraph of $\Gamma_{{\mathcal A}^0}$ associated to
$K$; for consistency with the statement of Proposition
\ref{keyconsequence} we shall write $K^0=H$. Since $\Gamma_1$
satisfies conditions (1) and (2) of \ref{keyinequality}, it now
follows that $K$ satisfies conditions (1) and (2) of
\ref{keyconsequence}. Our choice of $K$ also guarantees that it
satisfies conditions (3) and (4) of \ref{keyconsequence}.

Since $K$ is a subgraph of $\Gamma_1$, it follows from the definitions
that $\theta^\iota(K)\ge\theta^\iota(\Gamma_1)$. Combining this
observation with the inequalities (\ref{keyconsequenceone}) and
(\ref{keyconsequencetwo}), we deduce that
\Equation\label{keyconsequencethree}
\frac{\length(K_0)}{\theta^\iota(\Gamma_1)}<
6\phi_\tau(\theta_\infty)(\log_2|2\chi(|\Gamma_1|)|)
\frac{\sum_{i\in I}\lambda_i}{\Theta} .
\EndEquation
To estimate the factor $\log_2|\chi(|\Gamma_1|)|$ in
(\ref{keyconsequencethree}), note that since $\Gamma_1$ is a subgraph of
$\Gamma_{\mathcal A}$ it has at most $\Theta$ interior edges. Hence
$\chi(|\Gamma_1|)\ge\chi(|\Gamma_1|\cap\partial S)-\Theta$. But
$|\Gamma_1\cap\partial S|$ is a subpolyhedron of a triangulated $1$--manifold
and must therefore have non-negative Euler characteristic, so that
\Equation\label{keyconsequencefour}
|\chi(|\Gamma_1|)|\le\Theta .
\EndEquation
Condition (5) of  \ref{keyconsequence} follows immediately from
(\ref{keyconsequencethree}) and (\ref{keyconsequencefour}).
\EndProof

\section{Slopes and genera II}
\label{secondgenussection}
The goal of this section is to prove Theorem \ref{hardboundconsequence}.

\Definitions\label{properhomotopydef}
If $M$ is a manifold of arbitrary dimension, we define a {\it proper
path} in $M$ to be a map $\alpha\co I\to M$ such that $\alpha(\partial
I)\subset \partial M$. A {\it proper homotopy (of paths)} in $M$ is a
homotopy $H\co (I\times I)\to M$ such that $H((\partial I)\times
I)\subset\partial M$. Two proper paths $\alpha$ and $\beta$ in $M$ are
{\it properly homotopic} in $M$, or {\it properly $M$--homotopic}, if
there is a proper homotopy $H\co I\times I\to M$ such that $H_0=\alpha$
and $H_1=\beta$.

If $A$ is a properly embedded arc in a manifold $M$, a {\it
parametrization} of $A$, ie, a homeomorphism $\alpha\co I\to A$, is a
proper arc in $M$. We shall say that two properly embedded arcs are
{\it properly homotopic} if they admit properly homotopic
parametrizations.

Now suppose that $F$ is an essential surface in a compact, orientable,
irreducible $3$--manifold $M$. A reduced homotopy $H\co (I\times
I,I\times\partial I)\to(M,F)$ will be termed {\it proper} if the
homotopy $H\co (I\times I)\to M$ is proper. A proper reduced homotopy may
be regarded as a map of triples $H\co (I\times I,I\times\partial
I,\partial I\times I)\to(M,F,\partial M)$.  By a {\it proper reduced
homotopy of length 0} we will mean a homotopy $H\co (I\times
I,I\times\partial I)\to(F,\partial F)$.
\EndDefinitions

\Lemma\label{reduceit}
Suppose that $F$ is an essential surface in an irreducible knot
manifold $M$, and suppose that $\alpha$ and $\beta$ are proper paths
in $F$ which are properly $M$--homotopic. Then there is a proper
reduced homotopy $H$ in $(M,F)$ such that $H_0=\alpha$ and
$H_1=\beta$.
\EndLemma

\Proof
If $M$ is a solid torus then $F$ must be a disk, so that
$\alpha$ and $\beta$ are homotopic in $F$. Thus in this case we may
take $H$ to be a length--$0$ homotopy. We may therefore restrict
attention to the case in which $\partial M$ is $\pi_1$--injective.

Let $J\co I\times I\to M$ be a proper homotopy such that $J_0=\alpha$ and
$J_1=\beta$. For $i=0,1$, let $\gamma_i\co I\to\partial M$ be the path
defined by $\gamma_i(t)=J(i,t)$. We first consider the case in which
$\gamma_1$ is fixed-endpoint homotopic in $\partial M$ to a path in
$\partial F$. In this case we may assume that
$\gamma_1(I)\subset\partial F$. The path $\gamma_0$ is fixed-endpoint
homotopic to a composition $\gamma_0'$ of the paths $\alpha$,
$\gamma_1$ and the inverse path $\bar\beta$ of $\beta$. We have
$\gamma_0'(I)\subset F$ and $\gamma_0(I)\subset\partial M$.  Applying
Lemma \ref{boundaryincompressible}, with $\gamma_0'$ in place of
$\alpha$, we conclude that $\gamma_0'$ is fixed-endpoint homotopic in
$F$ to a path in $\partial F$. Hence $\gamma_0$ is fixed-endpoint
homotopic in $M$ to a path in $\partial F$, and so after modifying the
map $J$, without changing its values on $(I\times\partial
I)\cup(\{1\}\times I)$, we may assume that $J(\partial(I\times
I))\subset F$. Since $F$ is $\pi_1$--injective in $M$, there is a map
$H\co I\times I\to F$ which agrees with $J$ on the boundary of the disk
$I\times I$. Then $H$ is a length--$0$ homotopy from $\alpha$ to
$\beta$, and the conclusion holds in this case. If we assume that
$\gamma_0$ is fixed-endpoint homotopic in $\partial M$ to a path in
$\partial F$, the argument is precisely similar.

Now suppose that neither $\gamma_0$ nor $\gamma_1$ is fixed-endpoint
homotopic in $\partial M$ to a path in $\partial F$.  After modifying
the $\gamma_i$ within their fixed-endpoint homotopy classes in
$\partial M$, we may assume that each $\gamma_i$ is transverse to
$\partial F$, so that $\gamma_i^{-1}(F)$ is a finite set
$\{t_{i,0}\ldots,t_{i,n_i}\}$, where $0=t_{i,0}<\cdots<t_{i,n_i}=1$.
For $i=0,1$ and for $j=1,\ldots,n_i$, let $\delta_{i,j}\co I\to\partial
M$ denote a path which is a reparametrization of
$\gamma_i\vert_{[t_{i,j-1},t_{i,j}]}$. Since $\partial M$ is a torus,
and since in particular no $\gamma_i$ is fixed-endpoint homotopic in
$\partial M$ to a path in $\partial F$, we may assume the $\gamma_i$
to have been chosen within their fixed-endpoint homotopy classes in
$\partial M$ in such a way that no $\delta_{i,j}$ is fixed-endpoint
homotopic in $\partial M$ to a path in $\partial F$.

If some $\delta_{i,j}$ is fixed-endpoint homotopic in $M$ to a path
$\delta'$ in $F$, then by applying Lemma \ref{boundaryincompressible},
with $\delta'$ in place of $\alpha$, we conclude that $\delta'$ is
fixed-endpoint homotopic in $F$ to a path in $\partial F$. Hence
$\delta_{i,j}$ is fixed-endpoint homotopic in $M$ to a path in
$\partial F$. Since $\partial M$ is $\pi_1$--injective, it follows that
$\delta_{i,j}$ is fixed-endpoint homotopic in $\partial M$ to a path
in $\partial F$, a contradiction.  Hence no $\delta_{i,j}$ is
fixed-endpoint homotopic in $M$ to a path in $F$.

After further modifications of the map $J$, which do not change its
values on $\partial(I\times I)$, we may assume that that
$J\vert_{(I\times\inter I)}$ is transverse to $F$. If some component
$C$ of $J^{-1}(F)$ is a simple closed curve, then the
$\pi_1$--injectivity of $F$ implies that $J\vert_C$ is homotopically
trivial in $F$. Hence if $D\subset I\times I$ is the disk bounded by
$C$, we may modify $J$ on a small neighborhood of $D$ to obtain a map
$J'\co  I\times I\to M$, agreeing with $J$ on $\partial(I\times I)$, such
that $J\vert_{(I\times\inter I)}$ is transverse to $F$ and such that
$(J')^{-1}(F)$ has fewer components than $J^{-1}(F)$. After a finite
number of such modifications we may assume that every component of
$J^{-1}(F)$ is an arc.

If some component of $J^{-1}(F)$ has both its endpoints in the same
component $\{i\}\times I$ of $(\partial I)\times I$, then there is a
disk $D\subset I\times I$ such that $D\cap\partial(I\times I)\subset
\{i\}\times I$ and $\frontier D$ is a component of $J^{-1}(F)$. Among
all disks with these properties, we may suppose $D$ to be chosen so as
to be minimal with respect to inclusion. If we set $A=\frontier D$, we
then have $\partial A=\{t_{i,j-1},t_{i,j}\}$ for some $j$ with
$1<j<n_i$. The map $J\vert_D$ defines a fixed-endpoint homotopy from
$\delta_{i,j}$ to a path in $F$ which is a reparametrization of
$J\vert_A$. This is a contradiction. Hence every component of
$J^{-1}(F)$ is an arc which has one endpoint in $\{0\}\times I$ and
one in $\{1\}\times I$.

It follows that by pre-composing $J$ with a self-homeomorphism of
$I\times I$ which is the identity on $I\times\partial I$, we obtain a
homotopy $H\co I\times I\to M$ such that $H^{-1}(F)$ has the form
$I\times Y$ for some finite set $Y\subset I$. Hence $H$ is a
composition of basic homotopies $H^{(1)},\ldots,H^{(n)}$. The fact
that no $\delta_{i,j}$ is fixed-endpoint homotopic in $M$ to a path in
$F$ implies that $H^{(1)},\ldots,H^{(n)}$ are all essential. The fact
that $J\vert_(I\times\inter I)$ is transverse to $F$ implies that,
given a transverse orientation of $F$, for each $i\in\{1,\ldots,n-1\}$
there is an element $\omega$ of $\{-1,+1\}$ such that $H^i$ ends on
the $\omega$ side and $H^{i+1}$ starts on the $-\omega$ side. Thus $H$
is the required reduced homotopy in $(M,F)$ from $\alpha$ to $\beta$.
\EndProof

\Lemma\label{normalize}
Suppose that $F$ is a transversally oriented essential surface in an
irreducible knot manifold $M$. Then there exist a finite set $Y\subset
S^1$ and a homeomorphism $J\co S^1\times S^1\to\partial M$ such that
$J^{-1}(\partial F)=S^1\times Y$.  Furthermore, if $J$ and $Y$ have
this property, and if $H\co (I\times I,I\times\partial I,\partial I\times
I)\to(M,F,\partial M)$ is a proper reduced homotopy of length $k>0$,
then there is a proper reduced homotopy $H'\co (I\times I,I\times\partial
I,\partial I\times I)\to(M,F,\partial M)$ such that
\Parts
  \Part{(1)} $H'$ has length $k$, starts on the same side as $H$ and
  ends on the same side as $H$,
  \Part{(2)} $H'_0=H_0$,
  \Part{(3)} $H'_1$ is properly $F$--homotopic to $H_1$, and
  \Part{(4)} for each $s\in\{0,1\}$, the path $t\mapsto H(s,t)$ is an
  immersion of $I$ in $J(\{c_s\}\times S^1)$ for some $c_s\in S^1$.
\EndParts
\EndLemma

\Proof
The first assertion of the lemma, about the existence of $J$, follows
from the fact that the components of $\partial F$ are disjoint
homotopically non-trivial simple closed curves on the torus $\partial
M$. To prove the second assertion, about the reduced proper homotopy
$H$, we argue by induction on the length $k$ of $H$. If $k=1$ then $H$
is an essential basic homotopy in $(M,F)$; hence for $s=0,1$ the path
$t\mapsto H(s,t)$ is a basic essential path (\ref{reducedhomotopydef})
in the pair $(\partial M,\partial F)$, which is homeomorphic via $J$
to $(S^1\times S^1,S^1\times Y)$. But for any basic essential path
$\alpha$ in $(S^1\times S^1,S^1\times Y)$, there is a homotopy
$A\co I\times I\to M$, constant on $0\in I$, such that for every $s\in I$
the path $A_s$ is a basic essential path in $(\partial M,\partial F)$,
and such that $A_0=\alpha$ and $A_1$ is an immersion of $I$ in
$J(c\times S^1)$ for some $c\in S^1$.  The existence of the required
homotopy $H'$ therefore follows from the homotopy extension property
for polyhedra.

Now assume that $k>1$ and that the assertion is true for reduced
proper homotopies of length $k-1$. We write $H$ as a composition of
$k$ essential basic homotopies $H^1, \ldots, H^k$ in such a way that,
for each $j\in\{1,\ldots,n-1\}$ there is an element $\omega_j$ of
$\{-1,+1\}$ such that $H^j$ ends on the $\omega_j$ side and $H^{j+1}$
starts on the $-\omega_j$ side.  A composition $H^*$ of
$H^1,\ldots,H^{k-1}$ is a reduced proper homotopy of length $k-1$.
After applying the induction hypothesis to $H^*$ and using the
homotopy extension property, we may assume that for each
$s\in\{0,1\}$, the path $t\mapsto H^*(s,t)$ is an immersion of $I$ in
$J(\{c_s\}\times S^1)$ for some $c_s\in S^1$.  According to the case
$k=1$ of our assertion, which has already been proved, there is a
basic essential homotopy $(H^k)'$, starting on the $-\omega_{k-1}$
side and ending on the same side as $H^k$, such that
$(H^k)'_0=H^k_0=H*_1$, $(H^k)'_1$ is properly $F$--homotopic to
$H^k_1=H_1$, and $t\mapsto (H^k)'(s,t)$ is an immersion of $I$ in
$J(\{c_s\}\times S^1)$. We define $H'$ to be a composition of $H^*$
and $(H^k)'$. Since $t\mapsto H^*(s,t)$ and $t\mapsto (H^k)'(s,t)$ are
immersions of $I$ in $J(\{c_s\}\times S^1)$, and since $H^*$ ends on
the $\omega_{k-1}$ side and $H^k$ starts on the $\omega_{k-1}$ side,
it follows that $t\mapsto H'(s,t)$ is also an immersion of $I$ in
$J(\{c_s\}\times S^1)$, and the induction is complete.  \EndProof

\Proposition\label{etawelldefined}
Suppose that $F$ is a transversally oriented essential surface in an
irreducible knot manifold $M$.  Suppose that $H,H'\co (I\times
I,I\times\partial I,\partial I\times I)\to(M,F,\partial M)$ are proper
reduced homotopies with $\length(H)=\length(H')>0$. Suppose that
$H_0$ and $H_0'$ are properly $F$--homotopic and that $H$ and $H'$
start on the same side (in the sense of \ref{reducedhomotopydef}).
Then $H$ and $H'$ end on the same side, and $H_1$ and $H_1'$ are
properly $F$--homotopic.
\EndProposition

\Proof
We may assume without loss of generality that $H_0=H_0'$.  According
to the first assertion of Lemma \ref{normalize}, we may fix a finite
set $Y\subset S^1$ and a homeomorphism $J\co S^1\times S^1\to\partial M$
such that $J^{-1}(\partial F)=S^1\times Y$. According to the second
assertion of \ref{normalize}, the proof of the proposition reduces to
the case in which, for each $s\in\{0,1\}$, the paths $t\mapsto H(s,t)$
and $t\mapsto H'(s,t)$ are immersions of $I$ into submanifolds
$J(\{c_s\}\times S^1)$ and $J(\{c_s'\}\times S^1)$ of $\partial M$ for
some $c_s,c_s'\in S^1$.  Since $H_0=H_0'$, we have $c_s=c'_s$ for
$s=0,1$. In this case, we shall prove that $H_1$ and $H_1'$ are
properly $F$--homotopic, which includes the conclusion of the
proposition.

By hypothesis the proper reduced homotopies $H$ and $H'$ start on the
same side and have the same length. Since $H$ and $H'$ are now assumed
to be immersions of $I$ into $J(\{c_i\}\times S^1)$, it follows that
for each $s\in\{0,1\}$ the path $\gamma_s'\co t\mapsto H'(s,t)$ is a
reparametrization of $\gamma_s\co t\mapsto H(s,t)$.

Let $\bar\gamma_0$ and $\bar\gamma_0'$ denote the inverses of the
paths $\bar\gamma_0$ and $\bar\gamma_0'$. Let $\alpha$ denote a
composition of the paths $\bar\gamma_0$, $H_0$ and $\gamma_1$;
likewise, let $\alpha'$ denote a composition of $\bar\gamma_0'$,
$H_0'=H_0$ and $\gamma_1'$. Then $\alpha'$ is a reparametrization of
$\alpha$. But the existence of the homotopies $H$ and $H'$ imply that
the path $H_1$ is fixed-endpoint homotopic to $\alpha$ in $M$, and
that the path $H_1'$ is fixed-endpoint homotopic to $\alpha'$ in
$M$. Hence the paths $H_1$ and $H_1'$ are fixed-endpoint homotopic to
each other in $M$. Since these paths lie in the $\pi_1$--injective
surface $F\subset M$, they are in fact fixed-endpoint homotopic in
$F$, as required.
\EndProof

\Definition
Let $F$ denote a transversally oriented essential surface in an
irreducible knot manifold $M$. For any proper path $\alpha$ in $F$ and
any $\omega\in\{-1,+1\}$, we define the {\it $\omega$ height} of
$\alpha$, denoted $\height_\omega(\alpha)$, to be the supremum of all
integers $k\ge0$ for which there exists a proper reduced homotopy
$H\co (I\times I,I\times\partial I,\partial I\times I)\to(M,F,\partial
M)$ which starts on the $\omega$ side, has length $k$, and satisfies
$H_0=\alpha$. Thus $\height_\omega(\alpha)$ is either a non-negative
integer or $+\infty$. We define the {\it minheight} of $\alpha$,
denoted $\minheight(\alpha)$, to be
$\min(\height_{+1}(\alpha),\height_{-1}(\alpha))$. It is clear that
$\height_{+1}(\alpha)$, $\height_{-1}(\alpha)$ and
$\minheight(\alpha)$ depend only on the proper homotopy class of
$\alpha$ in $F$.  Furthermore, $\minheight(\alpha)$ is independent of
the choice of a transverse orientation of $F$ in $M$.
\EndDefinition

\Lemma\label{sourceexists}
Suppose that $F$ is a transversally oriented essential surface in an
irreducible knot manifold $M$. For any proper path $\alpha$ in $F$,
any integer $k>0$, and any $\omega\in\{-1,+1\}$, the following
conditions are equivalent:
\Parts
  \Part{(i)} $\height_\omega(\alpha)=k$; and

  \Part{(ii)} there exists a proper reduced homotopy $H\co (I\times
  I,I\times\partial I,\partial I\times I)\to(M,F,\partial M)$ which
  starts on the $\omega$ side and has length $k$, and such that
  $H_0=\alpha$ and $\minheight(H_1)=0$.
\EndParts
\EndLemma

\Proof
First suppose that (i) holds. It follows from the definition of
$\omega$--height that there is a proper reduced homotopy $H$ of length
$k$, starting on the $\omega$ side, such that $H_0=\alpha$. Set
$\beta=H_1$. Define an element $\epsilon$ of $\{-1,+1\}$ by the
condition that $H$ ends on the $\epsilon$ side, and set
$k^*=\height_{-\epsilon}(\beta)>0$.  If $\minheight(\beta)>0$, then in
particular $k^*>0$. Again from the definition, there is a proper
reduced homotopy $H^*$ of length $k^*$, starting on the $-\epsilon$
side, such that $H^*_1=\beta$. Then a composition of $H$ and $H^*$ is
a proper reduced homotopy of length $k+k^*$, starting on the $\omega$
side and having length $k+k^*>k$. This is a contradiction since
$\height_\omega(\alpha)=k$. Hence we must have $\minheight(\beta)=0$,
and (ii) is established.

Conversely, if (ii) holds, then the $\omega$--height $\ell$ of $\alpha$
is by definition $\ge k$.  Assume that $\ell>k$, and fix a proper
reduced homotopy $J$ of length $\ell$, starting on the $\omega$ side,
such that $J_0=\alpha$. Set $k^*=\ell-k>0$.  It follows from the
definition of a reduced homotopy that we may write $J$ as a
composition of a reduced homotopy $H'$ of length $k$ and a reduced
homotopy $H^*$ of length $k^*$, and that for some
$\epsilon\in\{-1,+1\}$, the homotopy $H'$ ends on the $\epsilon$ side
while $H^*$ starts on the $-\epsilon$ side. Set $\beta=H_1$ and
$\beta'=H_1'=H^*_0$. Since $H'$ and the homotopy $H$ given by (ii)
both start on the $\omega$ side, have length $k$ and satisfy
$H_0=H'_0=\alpha$, it follows from Lemma \ref{etawelldefined} that
$\beta$ and $\beta'$ are properly $F$--homotopic.

Since $\beta'=H^*_0$ we have $\height_{-\epsilon}(\beta')\ge k^*>0$.
On the other hand, the inverse homotopy $\bar H'$ of $H'$ starts on
the $\epsilon$ side and has length $k$; hence
$\height_{\epsilon}(\beta')\ge k>0$. It follows that
$\minheight(\beta)>0$. But $\beta=H_1$ has minheight $0$ according to
(ii), and since $\beta$ and $\beta'$ are properly $F$--homotopic they
must have the same minheight. This is a contradiction. It follows that
$\ell=k$, so that (i) holds.
\EndProof

\Proposition\label{atmosttwo}
Suppose that $F$ is an essential surface in an irreducible knot
manifold $M$, that $m$ is a non-negative integer, and that $\mathcal
C$ is a proper $M$--homotopy class of proper paths in $M$. Then
$\mathcal C$ contains at most $2m+2$ proper $F$--homotopy classes of
proper paths in $F$ which are of minheight at most $m$.
\EndProposition

\Proof
We first give the proof in the special case $m=0$.  Suppose that
$\alpha$, $\beta$ and $\beta'$ are proper arcs of minheight $0$ in
$F$, and that no two of them are properly $F$--homotopic. Since
$\minheight(\alpha)=0$, we may assume by symmetry that
$\height_{-1}(\alpha)=0$. Since $\beta$ and $\beta'$ are properly
$M$--homotopic to $\alpha$, it follows from Lemma \ref{reduceit} that
there are proper reduced homotopies $H$ and $H'$ such that
$H_0=H'_0=\alpha$, $H_1=\beta$, and $H'_1=\beta'$. Since neither
$\beta$ nor $\beta'$ is properly $F$--homotopic to $\alpha$, both $k$
and $\ell$ are of strictly positive length. Since
$\height_{-1}(\alpha)=0$, it follows that $H$ and $H'$ must both start
on the $+1$ side. Since $\beta$ and $\beta'$ have minheight $0$, it
follows from Lemma \ref{sourceexists} that
$\length(H)=\height_{+1}(\alpha)=\length(K)$. Hence according to Lemma
\ref{etawelldefined}, the proper paths $\beta$ and $\beta'$ are
properly $F$--homotopic, a contradiction. This completes the proof in
the case $m=0$.

Now suppose that $m>0$. For each proper path $\alpha$ in $F$, let
$[\alpha]$ denote the proper $F$--homotopy class of $\alpha$. Let $Z$
denote the set of all classes $[\alpha]$ such that
$\minheight(\alpha)=0$; by the case of the proposition already proved,
we have $\#(Z )\le2$. For each $[\alpha]\in Z $, and each
$k\in\{0,1,\ldots,m\}$, let $\HH_{([\alpha],k)}$ denote the set of all
classes $[\gamma]$ for which there exists a reduced homotopy $H$ of
length $k$ such that $H_0=\gamma$ and such that $H_1$ is properly
$F$--homotopic to $\alpha$. It follows from Lemma \ref{sourceexists}
that the set
$$
\HH=\bigcup_{([\alpha],k)\in Z\times\{0,\ldots,m\}}\HH_{([\alpha],k)}
$$
contains all classes $[\gamma]$ for which
$
0<\height_{+1}(\gamma)\le m
$ or 
$
0<\height_{-1}(\gamma)\le m
$.
Since
$\HH$ obviously also contains all classes $[\gamma]$ for which
$\minheight(\gamma)=0$, it therefore contains all classes $[\gamma]$
for which $\minheight(\gamma)\le m$.  Thus it suffices to prove that
$\#(\HH)\le2m+2$; and since $\#(Z )\le2$, it suffices to prove that
$\#(\HH_{([\alpha],k)})\le 1$ for each $([\alpha],k)\in Z\times \{0,
\ldots , m\}$.  This is clear for $k=0$, because
$\HH_{([\alpha],0)}=\{[\alpha]\}$ in view of the definition of a
length--$0$ reduced homotopy.

If $k>0$, and if $[\gamma]$ and $[\gamma']$ are elements of
$\HH_{(\alpha,k)}$, we have proper reduced homotopies $H,H'\co (I\times
I,I\times\partial I,\partial I\times I)\to(M,F,\partial M)$ which have
length $k$, and such that $H_0=\gamma$, $H'_0=\gamma'$, and both $H_1$
and $H_1'$ are properly $F$--homotopic to $\alpha$.  The inverse
homotopies $\bar H$ and $\bar H'$ are reduced homotopies of length $k$
such that $\bar H_0$ and $\bar H_0'$ are properly $F$--homotopic to
$\alpha$, while $H_1=\gamma$ and $H'_1=\gamma'$.  If $\bar H$ and
$\bar H'$ start on different sides then the definition of height
implies that both $\height_{-1}(\alpha)$ and $\height_{+1}(\alpha)$
are strictly positive, a contradiction since
$\minheight(\alpha)=0$. Hence $\bar H$ and $\bar H'$ start on the same
side, and it follows from Lemma \ref{etawelldefined} that $\gamma=H_1$
and $\gamma'=H_1'$ are properly $F$--homotopic,
ie, $[\gamma]=[\gamma']$. This shows that $\#(\HH_{([\alpha],k)})\le
1$, as required.
\EndProof
 
\Definition\label{embeddedminheight}
Let $F$ be an essential surface in a compact, orientable, irreducible
$3$--manifold $M$, and let $A$ be a properly embedded arc in $F$. By
the {\it minheight} of $A$ we mean the minheight of any
parametrization of $A$. It is clear that all parametrizations have the
same minheight.
\EndDefinition

\Lemma\label{extendthehomotopy}
Suppose that $F$ is an essential surface in an irreducible knot
manifold $M$.  Suppose that $k$ is a non-negative integer, that
$\mathcal A$ is an essential arc system in $F$ (in the sense of
\ref{systemdef}), and that $K$ is a subgraph of $\Gamma_{\mathcal A}$
such that for every interior edge $e$ of $K$, the arc $\bar e$ has
minheight at least $k$. Then the thickness of $|K|$ (in the sense of
\ref{thicknessdef}) is at least $2k+1$.
\EndLemma
 
\Proof
We may assume that $k>0$, as the case $k=0$ is trivial. We fix a
transverse orientation of $F\subset M$.
 
Let $r\co |K|\to F$ denote the inclusion map.  To establish the
conclusion of the lemma, it suffices to show that there exist reduced
homotopies $H^+\co (|K|\times I,|K|\times\partial I)\to(M,F)$ of length
$k$ such that $H^+$ starting on the $+1$ side, $H^-$ starts on the
$-1$ side, and $H^{+1}_0=H^{-1}_0=r$. Indeed, if $H^+$ and $H^-$ are
such homotopies, then by composing the inverse $\bar H^+$ of $H^+$
with $H^-$, we obtain a reduced homotopy $H$ of length $2k$ such that
$H_t=r$ for some $t\in I$; according to the definition, this implies
that $|K|$ has thickness at least $2k+1$.  By symmetry, it suffices to
construct $H^+$.
 
Let $m$ denote the number of components of $\partial F$, and set
$\zeta=e^{2\pi \sqrt{-1}/m}\in S^1$. It follows from the first
assertion of Lemma \ref{normalize} that there is a homeomorphism
$J\co S^1\times S^1\to\partial M$ such that $J^{-1}(\partial
F)=S^1\times\{1,\zeta,\ldots,\zeta^{m-1}\}$.  Let us denote by
$p\co S^1\times S^1\to S^1$ the projection to the second factor. For each
$j\in\{0,\ldots,m-1\}$ the standard orientation of $S^1$ defines a
transverse orientation of the $0$--manifold $\{\zeta^j\}\subset S^1$,
which in turn pulls back via the submersion $p\circ H^{-1}$ to a
transverse orientation of the component $C_j=J(S^1\times\{\zeta^j\})$
of $\partial F\subset\partial M$. For each $j$, we define $\sigma_j$
to be $+1$ if this pulled back transverse orientation agrees with the
transverse orientation of $C_j\subset\partial M$ induced by the given
transverse orientation of $F\subset M$, and to be $-1$ otherwise.
 
For each interior edge $e$ of $K$ let us fix a parametrization
$\alpha_e$ of $\bar e$. For each $e$ and for each $s\in\{0,1\}$ we
have $\alpha_e(s)\in C_{j(e,s)}$ for a unique
$j(e,s)\in\{0,\ldots,m-1\}$.  Since by hypothesis the arc $\bar e$ has
minheight at least $k$, there is a proper reduced homotopy
$H^e\co (I\times I,I\times\partial I,\partial I\times I)\to(M,F,\partial
M)$ which starts on the $+1$ side and has length $k$, and such that
$H^e_0=\alpha_e$.  According to the second assertion of Lemma
\ref{normalize}, we may suppose the $H^e$ to be chosen so that for
each interior edge $e$ and each $s\in\{0,1\}$, the path $t\mapsto
H^e(s,t)$ is an immersion of $I$ in $J(\{c_s^e\}\times S^1)$ for some
$c_s^e\in S^1$.  In view of our description of $\partial F$ and our
definition of $j(e,s)$ and $\sigma_j$, we may assume after suitable
reparametrization that the immersion $t\mapsto H^e(s,t)$ is given by
$$
t\mapsto J(c_s^e,\zeta^{j(e,s)}\exp({2\pi\sigma_{j(e,s)}kt\sqrt{-1}/m}))
$$
for each interior edge $e$ and each $s\in\{0,1\}$.

We may now define the required reduced homotopy $H^+\co (|K|\times
I,|K|\times\partial I)\to(M,F)$ by setting
$$
H^+(J(x,\zeta^j),t)=J(x,\zeta^j \exp(2\pi \sigma_{j}kt \sqrt{-1}/m))
$$
for every $j\in \{0,\ldots,m-1\}$ and every $x\in S^1$, and
$$
H^+(\alpha_e(s),t)=H^e(s,t)
$$
for every  interior edge $e$ of $K$ and for all $s,t\in I$.
\EndProof

\Lemma\label{nodumbintersections}
Suppose that $M$ is an irreducible knot manifold which contains no
essential annulus. Suppose that $F_1$ and $F_2$ are essential surfaces
in $M$ which intersect transversally, and that $\partial F_1$ and
$\partial F_2$ intersect minimally in the sense of
\ref{minimalintersection}. Suppose that $\alpha$ is a component of
$F_1\cap F_2$.  Then $\alpha$ is not properly $M$--homotopic to a path
in $\partial M$.
\EndLemma

\Proof
If $M$ is a solid torus then the components of $F_1$ and $F_2$ are
disks; the boundary slopes of $F_1$ and $F_2$ must be the same, and
minimality implies that $\partial F_1\cap\partial F_2=\emptyset$.
Thus the statement is vacuously true.  Now suppose that the
irreducible knot manifold $M$ is not a solid torus. Then $\partial M$
is $\pi_1$--injective. Suppose that some component $\alpha$ of $F_1\cap
F_2$ is an arc which is properly $M$--homotopic to a path in $\partial
M$. Then a parametrization of $\alpha$ is fixed-endpoint homotopic to
a path in $\partial F_i$ for $i=1$ and for $i=2$. Hence $\alpha$ is
parallel in each of the $F_i$ to an arc $\beta_i\subset\partial F_i$
with $\partial\beta_i=\partial\alpha_i$. In particular, $\beta_1$ and
$\beta_2$ are fixed-endpoint homotopic arcs in $\partial M$.  Since
$\partial M$ is $\pi_1$--injective, $\beta_1$ and $\beta_2$ are
fixed-endpoint homotopic in $\partial M$. But since $\partial M$ is a
torus and the $1$--manifolds $\partial F_1$ and $\partial F_2$
intersect minimally, no arc in $\partial F_1$ can be fixed-endpoint
homotopic in $\partial M$ to an arc in $\partial F_2$.
\EndProof

The following result was observed by Cameron Gordon, who used it in an
unpublished argument giving a bound for the geometric intersection
number of the boundary slopes of two essential surfaces in a
hyperbolic knot manifold in terms of the intrinsic topological
invariants of the surfaces (cf \cite[Corollary 6.2.5]{veggie}).

\Lemma\label{notbothways}
Suppose that $M$ is an irreducible knot manifold which contains no
essential annulus. Suppose that $F_1$ and $F_2$ are essential surfaces
in $M$ which intersect transversally, and that $\partial F_1$ and
$\partial F_2$ intersect minimally in the sense of
\ref{minimalintersection}. Suppose that $A$ and $A'$ are distinct
components of $F_1\cap F_2$ which are both arcs. Then $A$ and $A'$
cannot be parallel both on $F_1$ and on $F_2$.
\EndLemma

\Proof
As in the proof of Lemma \ref{nodumbintersections} we can
show that the statement is vacuously true if $M$ is a solid torus.
Hence we may assume that $\partial M$ is $\pi_1$--injective. Suppose
that the arcs $A$ and $A'$ are parallel both on $F_1$ and on $F_2$
Then for $i=1,2$ there is a PL disk $R_i\subset F_i$ with $\frontier
R_i=A\cup A'$. Let us write the standard $1$--sphere $S^1$ as a union
of two arcs $r_1$ and $r_2$ with $r_1\cap r_2=\partial r_1=\partial
r_2=\{a,a'\}$, and let $s\co S^1\times I\to M$ be a map which maps
$r_i\times I$ homeomorphically onto $R_i$ for $i=1,2$, and maps
$\{a\}\times I$ and $\{a'\}\times I$ homeomorphically onto $A$ and
$A'$ respectively. Since $\partial F_1$ and $\partial F_2$ intersect
minimally, the arcs $s(r_1\times\{0\})\subset\partial F_1$ and
$s(r_2\times\{0\})\subset\partial F_2$ are not fixed-endpoint
homotopic. Hence $s\vert_{(S^1\times\{0\})}\co (S^1\times\{0\})\to\partial M$
induces an injective homomorphism of fundamental groups. Since
$\partial M$ is $\pi_1$--injective, $s\co S^1\times I\to M$ also induces
an injective homomorphism of fundamental groups. Furthermore, the
proper path $t\mapsto s(a,t)$ in $M$ is a parametrization of the arc
$A$, and hence according to Lemma \ref{nodumbintersections} it is not
properly homotopic to a path in $\partial M$. This shows that $s$,
regarded as a map from $S^1\times I, S^1\times \partial I)$ to
$(M,\partial M)$, is non-degenerate in the sense of \cite{charsub}.
It then follows from the annulus theorem
\cite[Theorem IV.3.1]{charsub} that $M$ contains an essential annulus,
in contradiction to the hypothesis.
\EndProof

\Lemma\label{hardkappabound}
Suppose that $F_1$ and $F_2$ are connected strict essential surfaces
in an irreducible knot manifold $M$.  Suppose that $F_1$ and $F_2$
intersect transversally, and that $\partial F_1$ and $\partial F_2$
are non-empty and intersect minimally in the sense of
\ref{minimalintersection}). Assume that the interior of every
component of $F_2-(F_1\cap F_2)$ is an open disk or an open
annulus. For $i=1,2$, let $g_i$, $s_i$ and $m_i$ denote, respectively,
the genus, boundary slope and number of boundary components of $F_i$,
and let $\chi_i=2-2g_i-m_i\le0$ denote its Euler characteristic.  Then
there exists a positive integer $\Theta\ge|\chi_2|$ with the following
property: if $q$ is any real number greater than $1$, and if we set
$\tau=(7q-1)/ ( q-1)$, then
$$
\kappa(F_1,F_2)\le\frac{36qm_2}{m_1}\cdot
\frac{\phi_\tau(\Theta)\log_2(2\Theta)}{\Theta}|\chi_1| ,
$$
where $\phi_\tau$ is defined by \ref{phidef}.
\EndLemma

\Proof
We set ${\mathcal A}^0=F_1\cap F_2$. It follows from Lemma
\ref{nodumbintersections} that ${\mathcal A}_0$ is an essential arc
system in $F_2$. We choose a reduction of ${\mathcal A}^0$ in $F_2$ in the
sense of \ref{reductiondef}, and denote it by ${\mathcal A}$.

We set $\Theta=\Theta_{\mathcal A}$. By definition, $\Theta$ is the number of
components of the arc system $\mathcal A$. Hence if $N$ denotes a regular
neighborhood of $\mathcal A$ in $F_2$ we have $\chi(\overline{F_2-N})=
\chi(F_2)+\Theta$. But by hypothesis, the interior of every component
of $F_2-{\mathcal A}=F_2-(F_1\cap F_2)$ is an open disk or an open
annulus, and hence $\chi(\overline{F_2-N})\ge0$. It follows that
\Equation\label{hardkappaboundone}
\Theta\ge-\chi(F_2)=|\chi_2|
\EndEquation
Hence in order to prove the lemma it suffices to show that, if we fix
any real number $q>1$, and if we set $\tau=(7q-1)/ ( q-1)$, then
\Equation\tag{$*$}
\kappa(F_1,F_2)\le\frac{36qm_2}{m_1}\cdot
\frac{\phi(\Theta)\log_2(2\Theta)}{\Theta}|\chi_1| .
\EndEquation
The components of $\mathcal A$ may be regarded as properly embedded
arcs in $M$. We denote by $\II$ the set of all proper homotopy classes
of proper paths in $M$, in the sense of \ref{properhomotopydef}, which
are represented by parametrizations of components of $\mathcal A$.  We
define a labeling of $\mathcal A$, in the sense of \ref{labelingdef},
by defining $\iota(e)$ to be the proper $M$--homotopy class of $\bar
e$, for every interior edge $e$ of $\Gamma_{\mathcal A}$.  We have
$\Theta=\Theta_{\mathcal A}=\sum_{i\in \II}\theta^\iota_i$.  As in the
statement of Proposition \ref{keyconsequence}, we set
$\theta_i=\theta_i^\iota$ for every $i\in \II$, and we set
$\theta_\infty=\max_{i\in \II}\theta_i$.

When we wish to think of an element $i$ of $\II$ as a proper
$M$--homotopy class, rather than as a ``label,'' we will denote it by
${\mathcal C}_i$.

For each interior edge of $\Gamma_{\mathcal A}$, the quantity
$\minheight\bar e$ is defined by \ref{embeddedminheight}. Let us
denote by $E^*$ the set of all interior edges of $\Gamma_{\mathcal A}$
such that $\minheight(\bar e)\le(\theta^\iota(e)/(2q))-1$. For every
$i\in \II$ set $E_i^*=E^*\cap E_i$ and $\theta_i^*=\#(E_i^*)$. For
every interior edge $e$ of $\Gamma_{\mathcal A}$, let $w(e)$ denote
the ${\mathcal A}^0$--width of $e$. For each label $i\in \II$, set
$\lambda(i)=\max_e w(e)$, where $e$ ranges over all interior edges of
$\Gamma_{\mathcal A}$ with label $i$.

We wish to apply Proposition \ref{keyconsequence}.  For this purpose
we must verify that $\theta_i^*\le\theta_i/q$ for every $i\in
\II$. According to the definitions we have $\theta_i^*=\#(E_i^*)$,
where $E_i^*$ is a collection of interior edges of $\Gamma_{\mathcal
A}$, whose closures all represent the proper $M$--homotopy class
${\mathcal C}_i$, and all have minheight at most
$(\theta^\iota(e)/(2q))-1$.

Since $\mathcal A$ is a reduced arc system, the closures of the edges
in $E_i^*$ represent pairwise distinct proper $F$--homotopy
classes. Applying Proposition \ref{atmosttwo} with $F=F_2$, ${\mathcal
C}={\mathcal C}_i$ and $m=(\theta^\iota(e)/(2q))-1$, we conclude that
$\theta_i^*\le2m+2=\theta_i/q$, as required for the application of
\ref{keyconsequence}.

Hence $\Gamma_{\mathcal A}$ has a subgraph $K$ satisfying conditions
(1)--(5) of Proposition \ref{keyconsequence}. Since
$\theta_\infty=\max_{i\in \II}\theta_i\le\sum_{i\in
\II}\theta^\iota_i=\Theta$, and since the function $\phi_\tau$ defined
in \ref{phidef} is monotone increasing, it follows from condition (5)
of \ref{keyconsequence} that
\Equation\label{hardkappaboundtwo}
\frac{\length(K^0)}{\theta^\iota(K)}
<\frac{6\phi_\tau(\Theta)\log_2(2\Theta)}{\Theta}\sum_{i\in\II}\lambda(i) ,
\EndEquation
where $K^0$ is the subgraph of $\Gamma_{{\mathcal A}^0}$ associated to
$K$.

In view of the definition of $E^*$, it follows from condition
(2) of \ref{keyconsequence} that for every interior edge $e$ of $K$ we have
$\minheight(\bar e)>(\theta^\iota(e)/(2q))-1$. Since by definition we
have $\theta^\iota(K)=\min_e\theta_e^{\iota}$, where $e$ ranges over
the interior edges of $K$, it follows that
$$
\minheight(\bar e)>(\theta^\iota(K)/(2q))-1
$$
for every interior edge $e$ of $K$. Applying Lemma
\ref{extendthehomotopy}, taking $k$ to be the least integer
$\ge(\theta^\iota(K)/(2q))-1$, we conclude that the thickness
$t_{F_2}(|K|)$ is at least $(\theta^\iota(K)/q)-1$.  Since
$t_{F_2}(|K|)$ is by definition a strictly positive integer, we have
$2t_{F_2}(|K|)\ge t_{F_2}+1 \ge \theta^\iota(K)/q$, and hence
\Equation\label{hardkappaboundthree}
t_{F_2}(|K|)\ge\frac{\theta^\iota(K)}{2q} .
\EndEquation
Combining (\ref{hardkappaboundtwo}) and (\ref{hardkappaboundthree}) we
obtain
\Equation\label{hardkappaboundfour}
\frac{\length(K^0)}{t_{F_2}(|K|)}<
\frac{12q\phi_\tau(\Theta)\log_2(2\Theta)}{\Theta}\sum_{i\in\II}\lambda(i) .
\EndEquation
Let $K^0_\partial$ denote the subgraph of $K^0$ with
$|K^0_\partial|=|K|\cap\partial F_2=|K^0|\cap\partial F_2$.  The
number of vertices of $K_0$ is equal to
$\#(|K^0_\partial|\cap{\mathcal A})$.  It follows from condition (4)
of \ref{keyconsequence} that all vertices of $K^0$ are of valence at
least $2$; hence $K^0$ has at most as many vertices as edges, ie,
\Equation\label{hardkappaboundfive}
\#(|K^0_\partial|\cap {\mathcal A})\le\length(K^0) .
\EndEquation
We may write $|K^0|=|K^0_\partial|\cup{\mathcal B}$, where $\mathcal
B$ is the submanifold of $\mathcal A$ made up of the closures of the
interior edges of $K$.  Let ${\mathcal B}_1$ denote a properly
embedded submanifold of $F_2$ which is topologically parallel to
$\mathcal B$ in $F_2$ and disjoint from ${\mathcal A}$. Then the
$1$--dimensional polyhedron $L=|K^0_\partial|\cup{\mathcal B}_1$ is
isotopic in $F_2$ to $| K^0|=|K| $. Hence
\Equation\label{hardkappaboundsix}
t_{F_2}(L)=t_{F_2}(|K|) .
\EndEquation
The definition of $L$ implies that
$$
L\cap F_1=L\cap{\mathcal A}=|K^0_\partial|\cap{\mathcal A} ,
$$
so that by  (\ref{hardkappaboundfive}) we have
\Equation\label{hardkappaboundseven}
\#(L\cap F_1)\le\length(K^0) .
\EndEquation
Since $L$ is isotopic in $F_2$ to $|K|$, it follows from conditions
(1) and (3) of \ref{keyconsequence} that $L$ is connected, has first
Betti number equal to $2$ and is $\pi_1$--injective in $F_2$.  Hence
according to \ref{kappadef} we have
\Equation\label{hardkappaboundeight}
\kappa(F_1,F_2)\le \frac{m_2\cdot\#(L\cap F_1)}{m_1\cdot t_{F_2}(L)} .
\EndEquation
It follows from (\ref{hardkappaboundfour}), (\ref{hardkappaboundsix}),
(\ref{hardkappaboundseven}) and (\ref{hardkappaboundeight}) that
\Equation\label{hardkappaboundnine}
\kappa(F_1,F_2) \le \frac{12qm_2}{m_1}\cdot
\frac{\phi(\Theta)\log_2(2\Theta)}{\Theta}\sum_{i\in\II}\lambda(i) .
\EndEquation
We need an estimate for the factor $\sum_{i\in \II}\lambda(i)$ which
appears in the right hand side of (\ref{hardkappaboundnine}). By the
definition of the $\lambda(i)$, we may choose for each $i\in{\mathcal
I}$ an edge $e_i$ such that $\iota(e_i)=i$ and $w(e_i)=\lambda_i$. The
${\mathcal A}^0$--width $w(e_i)$ is by definition the cardinality of
the ${\mathcal A}^0$--parallelism class containing the component $\bar
e_i$ of ${\mathcal A}^0$. If we denote this parallelism class by
$C_i$, and set $C=\bigcup_{i\in{\mathcal I}}C_i$, it follows that
\Equation\label{hardkappaboundten}
\sum_{i\in\II}\lambda(i)=\#(C) .
\EndEquation
Suppose that two distinct arcs in $C$, say $A$ and $A'$, are parallel
in $F_1$. We have $A\in C_i$ and $A'\in C_j$ for some $i,j\in\II$.
Then $e_i$ and $e_j$ are respectively parallel in $F_2$ to $A$ and
$A'$, which by our assumption are parallel to each other in $F_1$;
hence $e_i$ and $e_j$ are properly homotopic as properly embedded arcs
in $M$, and from the definition of the labeling $\iota$ it follows
that $i=\iota(e_i)=\iota(e_j)=j$. Thus $A$ and $A'$ both belong to
$C_i$, and are therefore parallel in $F_2$ as well as in $F_1$. Since
by hypothesis the essential surfaces $F_1$ and $F_2$ intersect
transversally, and $\partial F_1$ and $\partial F_2$ intersect
minimally, this contradicts Lemma \ref{notbothways}. It follows that
no two distinct arcs in $C$ can be parallel in $F_1$.

The cardinality of a collection of pairwise non-parallel arcs in $F_1$
is at most $-3\chi(F_1)=3|\chi_1|$. In view of
(\ref{hardkappaboundten} ) it follows that
\Equation\label{hardkappaboundeleven}
\sum_{i\in \II}\lambda(i)\le 3|\chi_1| .
\EndEquation
The inequality ($*$) follows immediately from (\ref{hardkappaboundnine})
and (\ref{hardkappaboundeleven}), and the proof of the lemma is
complete.
\EndProof

The next two lemmas will be needed in order to obtain a concrete estimate
from Lemma \ref{hardkappabound}.

\Lemma\label{precalculus}
For every real number $x>1$ there is a real number $q$ such that
\Equation\tag{$1$}
1<q\le\frac{1}{7}\left({6}\left(\frac{2\ln x}{\ln7}\right)^{1/2}+{31}\right),
\EndEquation
and such that, if we set $\tau=(7q-1)/(q-1)$, we have
\Equation\tag{$2$}
\ln\phi_\tau(x)\le2(2(\ln7)(\ln x))^{1/2}+4\ln7+1 .
\EndEquation
\EndLemma

\Proof
We define a positive integer $\mu$ by
$$
\mu=\left\lfloor\left(\frac{\ln x}{2\ln7}\right)^{1/2}\right\rfloor+1 ,
$$
$$
q=\frac{12\mu+19}{7} .
\leqno{\hbox{and we set}}
$$
Then ($1$) is clear from direct computation. On the other hand, if we
set
$$
\tau=\frac{7q+1}{q+1}=7\left(1+\frac{1}{2m+2}\right) ,
$$
$$
\ln\tau\le\ln7+\frac{1}{2m+2} .
\leqno{\hbox{then}}
$$
The definition of $\phi_\tau$ (\ref{phidef}) implies that
$$
\ln\phi_\tau(x)=\min_m\left((2\mu+2)\ln\tau+\frac{1}{m}\ln x\right) ,
$$
where $m$ ranges over all positive integers. Hence
$$\begin{aligned}
\ln\phi_\tau(x) &\le(2\mu+2)\ln\tau+\frac{1}{\mu}\ln x\\
&\le(2\mu+2)\ln7+1+\frac{1}{\mu}\ln x .
\end{aligned}$$
$$
2\mu+2\le2\left(\frac{\ln x}{2\ln7}\right)^{1/2}+4
\qquad{\rm and }\qquad
\frac{1}{\mu}\le\left(\frac{2\ln7}{\ln x}\right)^{1/2} ,
\leqno{\hbox{Since}}
$$
we now obtain ($2$) by direct computation.
\EndProof

\Lemma\label{calculus}
Let us define a function $f$ on $(1,\infty)$ by
$$
f(x)=2(2(\ln7)(\ln x))^{1/2} +
\ln\left(6\left(\frac{2\ln x}{\ln7}\right)^{1/2}+31\right) +
\ln\left(\frac{\ln x}{\ln2}+1\right)-\ln x .
$$
Then for all integers $m\ge n\ge333$ we have $f(m)\le
f(n)$.
\EndLemma

\Proof
We set $\alpha=2\sqrt{2\ln7}$, $\beta=6\sqrt{2/\ln7}$, $\gamma=31$ and
$\delta=1/\ln2$. Then for $x>1$ we have
$$
f(x)=\alpha(\ln x)^{1/2}+\ln(\beta(\ln x)^{1/2}+\gamma) +
\ln(\delta\ln x+1)-\ln x ,
$$
and hence
\Equation\label{calculusone}
xf'(x)=\frac{\alpha}{2(\ln x)^{1/2}} +
\frac{\beta}{2\beta\ln x+2\gamma(\ln x)^{1/2}} +
\frac{\delta}{\delta\ln x+1}-1
\EndEquation
Since $\alpha$, $\beta$, $\gamma$ and $\delta$ are positive, the right
hand side of (\ref{calculusone}) is obviously monotonically decreasing
for $x>1$, and by direct calculation it is seen to be negative for
$x=334$ (and positive for $x=333$). Hence $f'(x)<0$ for $x\ge334$, and
it follows that the conclusion of the lemma holds in the case $m\ge
n\ge 334$.  Since by direct computation we find that $f(334)<f(333)$,
the conclusion also holds in the case $m\ge n= 333$.
\EndProof

\Proposition\label{explicithardkappabound}
Suppose that $F_1$ and $F_2$ are connected strict essential surfaces
in an irreducible knot manifold $M$.  Suppose that $F_1$ and $F_2$
intersect transversally, and that $\partial F_1$ and $\partial F_2$
are non-empty and intersect minimally in the sense of
\ref{minimalintersection}). Assume that the interior of every
component of $F_2-(F_1\cap F_2)$ is an open disk or an open annulus.
For $i=1,2$, let $g_i$, $s_i$ and $m_i$ denote, respectively, the
genus, boundary slope and number of boundary components of $F_i$, and
let $\chi_i=2-2g_i-m_i\le0$ denote its Euler characteristic.  Assume
that $|\chi_2|\ge333$. Then
$$
\kappa(F_1,F_2)\le g(\chi_2)\frac{m_2|\chi_1|}{m_1|\chi_2|}
$$
where
$$
g(x) = 
12348\
\left(6\left(\frac{2\ln|x|}{\ln7}\right)^{1/2}+31\right)
\left(\frac{\ln|x|}{\ln2}+1\right)\exp(1 + 2(2(\ln7)(\ln |x|))^{1/2}) .
$$
\EndProposition

\Proof
Let $\Theta$ be the integer given by Lemma \ref{hardkappabound}. Thus
$\Theta\ge|\chi_2|$, and if $q$ is any real number greater than $1$,
and if we set $\tau=(7q-1)/ ( q-1)$, then
\Equation\label{explicitone}
\kappa(F_1,F_2)\le\frac{36qm_2}{m_1}\cdot
\frac{\phi_\tau(\Theta)\log_2(2\Theta)}{\Theta}|\chi_1| .
\EndEquation
Applying Lemma \ref{precalculus} with $x=\Theta$, we obtain a
particular value of $q>1$ such that, if we set $\tau=(7q-1)/(q-1)$, we
have
\Equation\label{explicittwo}
q\le\frac{1}{7}\left({6}\left(2\frac{\ln\Theta}{\ln7}\right)^{1/2}+{31}\right) ,
\EndEquation
and
\Equation\label{explicitthree}
\ln\phi_\tau(\Theta )\le2(2(\ln7)(\ln\Theta ))^{1/2}+4\ln7+1 .
\EndEquation
If we define a function $f$ on $(1,\infty)$ by
\Equation\label{explicitfour}
f(x)=2(2(\ln7)(\ln x))^{1/2} +
\ln\left(6\left(2\frac{\ln x}{\ln7}\right)^{1/2}+31\right) +
\ln\left(\frac{\ln x}{\ln2}+1\right)-\ln x ,
\EndEquation
it follows from (\ref{explicittwo}) and (\ref{explicitthree}) that
$$
\ln\left(\frac{q\phi_\tau(\Theta)\log_2(2\Theta)}{\Theta}\right)\le
f(\Theta)+3\ln7+1 .
$$
But since $\Theta\ge|\chi_2|\ge333$, Lemma \ref{calculus} asserts that
$f(\Theta)\le f(|\chi_2|)$. Hence
\Equation\label{explicitfive}
\ln\left(\frac{q\phi_\tau(\Theta)\log_2(2\Theta)}{\Theta}\right)\le
f(|\chi_2|)+3\ln7+1 .
\EndEquation
Combining (\ref{explicitone}) and (\ref{explicitfive}) we obtain
\Equation\label{explicitsix}
\kappa(F_1,F_2)\le\frac{36m_2}{m_1}\cdot\exp(f(|\chi_2|)+3\ln7+1)\cdot|\chi_1|.
\EndEquation
If we set $x=|\chi_2|$ in (\ref{explicitfour}), substitute the
resulting expression for $f(|\chi_2|)$ into (\ref{explicitsix}) and
simplify, we obtain the conclusion of Proposition
\ref{explicithardkappabound}.
\EndProof

\Theorem\label{hardboundconsequence}
Suppose that $K$ is a non-exceptional two-surface knot in a closed,
orientable $3$--manifold $\Sigma$ such that $\pi_1(\Sigma)$ is cyclic.
Let $\emm$ denote the meridian slope of $K$ and let $F_1$ and $F_2$ be
representatives of the two isotopy classes of connected strict
essential surfaces in $M(K)$.  For $i=1,2$, let $g_i$, $s_i$ and $m_i$
denote, respectively, the genus, boundary slope (well-defined by
\ref{twosurfprops}) and number of boundary components of $F_i$, and
let $\chi_i=2-2g_i-m_i < 0$ denote its Euler characteristic.  Assume
that $|\chi_2|\ge333$ and that $s_2\ne\emm$.  Set
$q_i=\Delta(s_i,\emm)$ (so that $q_i$ is the denominator of $s_i$ in
the sense of \ref{denominator}), and set $\Delta=\Delta(s_1,s_2)$ (so
that $\Delta\ne0$ by \ref{twosurfprops}).  Then
$$
\frac{q_1^2}{\Delta}\le 2g(\chi_2)\frac{m_2|\chi_1|}{m_1|\chi_2|}
$$
where
$$
g(x) = 
12348\
\left(6\left(\frac{2\ln|x|}{\ln7}\right)^{1/2}+31\right)
\left(\frac{\ln|x|}{\ln2}+1\right)\exp(1 + 2(2(\ln7)(\ln |x|))^{1/2}) .
$$
\EndTheorem

\Proof
We may assume after an isotopy that $F_1$ and $F_2$ intersect
transversally, and that $\partial F_1$ and $\partial F_2$ intersect
minimally in the sense of \ref{minimalintersection}. Furthermore,
$F_1$ and $F_2$ may be assumed to be chosen within their rel-boundary
isotopy classes so as to minimize the number of components of $F_1\cap
F_2$. Then no component of $F_1\cap F_2$ is a homotopically trivial
simple closed curve (cf  Remark \ref{standardremark}). Set $A=F_1\cap
F_2$.  It now follows from Theorem \ref{disksandannuli} that every
component of $(\inter F_i)-A$ is an open disk or an open annulus.
Hence by Proposition \ref{explicithardkappabound} we have
$$
\kappa(F_1,F_2)\le g(\chi_2)\frac{m_2|\chi_1|}{m_1|\chi_2|} .
$$
On the other hand, according to Theorem \ref{generalknotinequality} we
have
$$
\frac{q_1^2}{\Delta}\le{\four}2\kappa(F_1,F_2) .
$$
The last two inequalities imply the inequality in the conclusion of
Theorem \ref{hardboundconsequence}.
\EndProof

Theorem \ref{hardboundconsequence} has the following qualitative
consequence.

\Corollary\label{hardqualitative}
There is a positive-valued function $f_1(x)$ of a positive real
variable $x$ with the following properties.
\Parts
  \Part{(i)} For every $\epsilon>0$ we have
  $$\lim_{x\to\infty}x^{1-\epsilon}f_1(x)=0.$$

  \Part{(ii)} If $K$ is any non-exceptional two-surface knot in a
  closed, orientable $3$--manifold $\Sigma$ such that $\pi_1(\Sigma)$
  is cyclic, and if $\emm$, $F_i$, $g_i$, $s_i$, $m_i$, $q_i$ and
  $\Delta$ are defined as in the statement of Theorem
  \ref{easyboundconsequence} or \ref{hardboundconsequence},
then
$$
\frac{q_1^2}{\Delta}\le\frac{m_2|\chi_1|}{m_1}f_1(|\chi_2|) .
$$
\EndParts
\EndCorollary

\Number\label{qualitativediscussion}
It is instructive to compare Corollary \ref{hardqualitative} with the
corresponding qualitative consequence of Theorem
\ref{easyboundconsequence}, Corollary \ref{easyqualitative}. Suppose
that $f_1$ is a function with the properties stated in Corollary
\ref{hardqualitative}. After replacing $f_1$ by the function
$x\mapsto\sup_{y\ge x}f_1(y)$, which clearly still has the same
properties, we may assume that $f_1$ is monotone decreasing. Now
suppose that $K$ is any non-exceptional two-surface knot in a closed,
orientable $3$--manifold $\Sigma$ such that $\pi_1(\Sigma)$ is cyclic.
With the notation of \ref{hardqualitative} we have
$$
\frac{q_1^2}{\Delta}\le\frac{m_2|\chi_1|}{m_1}f_1(|\chi_2|) .
$$
On the other hand, according to Corollary \ref{chibound} we have
$|\chi_1| \le m_1m_2\Delta/2$. Hence
$$
\left(\frac{q_1}{\Delta}\right)^2\le {m_2^2} \frac{f_1(|\chi_2|)}{2} .
$$

Since $f$ is monotone decreasing and since $|\chi_2|=2g_2+m_2-2\ge
g_2$ whenever $g_2\ge2$, it follows that
$$\left(\frac{q_1}{\Delta}\right)^2 \le {m_2^2} \frac{f_1(g_2)}{2}$$
provided that $g_2\ge2$.  Hence the conclusions of Corollary
\ref{easyqualitative} hold with $f_0(x)= f_1(x)/2$.
 
This shows that Corollary \ref{hardqualitative}, together with the
relatively easy Corollary \ref{chibound}, implies Corollary
\ref{easyqualitative} by a purely formal argument. This suggests that
Theorem \ref{hardboundconsequence} may be ``stronger'' than Theorem
\ref{easyboundconsequence} in some vague qualitative sense.  Indeed,
the argument given above suggests that when the inequality given by
\ref{chibound} is far from being an equality,
\ref{hardboundconsequence} will give stronger information than
\ref{easyboundconsequence}.
 
On the other hand, note that if we derive Corollary
\ref{easyqualitative} from Theorem \ref{easyboundconsequence}, we
get a function $f_0(x)$ such that $xf_0(x)$ grows like $\ln x$,
whereas if we derive it from Theorem \ref{hardboundconsequence} via
Corollary \ref{hardqualitative}, we get a function
$f_0(x)=f_1(x)/2$ such that $xf_0(x)$ grows like $(\ln
x)^{3/2}\exp(C(\ln x)^{1/2})$, and hence more rapidly than any power
of $\ln x$. In view of the argument given above this suggests that
Theorem \ref{easyboundconsequence} may give stronger information
than Theorem \ref{hardboundconsequence} in the case where the
inequality given by \ref{chibound} is close to being an equality.
\EndNumber

\Addresses\recd

\end{document}
\bye